\documentclass[final,11pt,reqno,twoside]{amsart}

\usepackage[margin=3cm, headsep=24pt, headheight=2cm]{geometry}
\usepackage{amsaddr}

\usepackage{gensymb}

\usepackage{listings}

\usepackage[colorlinks,linkcolor=blue,urlcolor=blue,citecolor=blue]{hyperref}

\usepackage{epsf,epsfig}

\usepackage{graphicx,ulem}

\usepackage{amsmath,amsfonts,amsthm,multirow}

\usepackage{amssymb}
\usepackage[numbers]{natbib}

\usepackage[usenames,dvipsnames,svgnames,table]{xcolor}

\usepackage{pgfplots,pgfplotstable}

\graphicspath{ {FIGURES/} }

\usepackage{caption}
\usepackage{subcaption}

\pgfplotsset{compat=1.13}

\begin{document}


\title[]{Implementation of a Volume-of-Fluid Method in a Finite Element Code with Applications to Thermochemical Convection in a Density Stratified Fluid in the Earth's Mantle}

\author{Jonathan M. Robey}
\address[Jonathan M. Robey]{Department of Mathematics, U.\ C. Davis, Davis, CA 95616, USA}
\email[Corresponding Author]{jmrobey@ucdavis.edu}
\author{Elbridge Gerry Puckett}
\address[Elbridge Gerry Puckett]{Department of Mathematics, U.\ C. Davis, Davis, CA 95616, USA}
\email[Corresponding Author]{egpuckett@ucdavis.edu}



\date{}

\begin{abstract}
    We describe the implementation of a second-order accurate volume-of-fluid interface 
    tracking algorithm in the open source finite element code ASPECT, which is designed to 
    model convection in the Earth's mantle.
    This involves the solution of the incompressible Stokes equations coupled to an
    advection diffusion equation for the temperature, a Boussinesq approximation that governs the dependence of the density on the temperature, and an advection equation for a marker indicating the two initial density states.
    The volume-of-fluid method is fully parallelized and is integrated with the adaptive mesh
    refinement algorithm in ASPECT.
    We present the results of several standard interface tracking benchmarks in order to demonstrate the accuracy of the method as well as the results of several benchmarks commonly used in the computational mantle convection community.
    Finally, we present the results of computations with and without adaptive mesh refinement of a model problem involving thermochemical convection in a computationally stratified fluid designed to provide insight into how thermal plumes, that eventually reach the Earth's surface as ocean island basalts, originate at structures on the core-mantle boundary known as Large Low Shear wave Velocity Provinces.
\end{abstract}


\maketitle




\section{Introduction}
\label{Section:Introduction}

Over more than the past four decades there have been many numerical methods developed to study convection and other processes in the Earth's mantle.
In particular, there have a been a sequence of codes developed over this period of time that are now freely available to any individual who wishes to study mantle dynamics.
They include HC~\cite{BHH-RJO-1981,BHH-RWC:1989,BS:2000}, ConMan~\cite{SDK-AR-BHH:1990}, CitCom~S~\cite{AKM-SZ:2004,ET-EC-PT-MG-MA:2006,SZ-MTZ-LM-MG:2000}, 
Citcom~CU~\cite{LM-MG:1996,SZ:2006} and ASPECT~\cite{TH-JD-RG-WB:2017,MK-TH-WB:2012}.
These codes, as well as others, can be downloaded from the Computational Infrastructure for Geodynamics (CIG) at U.C.~Davis.
\footnote{The CIG is an NSF funded, community driven organization that advances Earth science by developing and disseminating software for geophysics and related fields.}

There are a large number of problems associated with the Earth's mantle that contain one or more interfaces in some form or another.
Although there have been some very specialized computational models of interfaces in the mantle, for example, the dynamics of bubbles and plumes~\cite{MM:1996,MM-HAS:1993,MM-HAS-RJO:1993}, it is only recently that researchers have begun to implement interface tracking algorithms in codes designed to model convection and other processes in the entirety of the Earth's mantle; e.g.,~\cite{HS-ME:2010}.
However, to our knowledge, the Volume-of-Fluid (VOF) method has not yet been implemented in a code designed to model convection in the Earth's mantle or, more generally, used by researchers to model geodynamic flows.

In this article we describe the implementation of a second-order accurate VOF interface tracking algorithm in the open source finite element code ASPECT.
ASPECT is a parallel, extendible finite element code designed to model thermal convection and other processes in the Earth's mantle in two and three dimensions.  
It is built on the deal.II Finite Element Library~\cite{DA-WB-DD-TH-LH-MK-MM-JPP-BT-DW:2017,WB-RH-GK:2007}, which includes adaptive mesh refinement (AMR)~\cite{p4est} and has been shown to scale to thousands of processors~\cite{RG:2016}.
ASPECT has been extended to model other processes that occur in the mantle, such as  modeling grain size evolution in the mantle~\cite{JD-ZE-UF-RG-PM-RM:2017}, melt generation and migration~\cite{JD-TH:2016}, as well as other problems.
There is currently a very active community of researchers extending ASPECT to new problem areas and improving existing algorithms. 
Our VOF algorithm is fully parallelized and is designed to work efficiently with ASPECT's AMR algorithm.

Recent studies utilizing seismic imaging have revealed large regions with anomalous seismic properties in the lower mantle.
In particular, there are two dome-like regions beneath Africa and the Pacific Ocean with low shear-wave velocities that extend some 1000 km above the core-mantle boundary and have horizontal dimensions of several thousand kilometers~\cite{SC-BR:2012,SWF-BR:2015}.
Most interpretations propose that the heterogeneities are compositional in nature, differing from the surrounding mantle, an interpretation that would be consistent with chemical geodynamic models.
Based on geological and geochemical studies it has been argued that LLSVPs have persisted for billions of years~\cite{KB-BS-THT-MAS:2008}.
In this article we compute solutions to a model problem designed to understand the dynamics of plumes that form on the LLSVPs, entrains some of the material in the LLSVP that differs from the surrounding mantle, and brings it to the Earth's surface.
The model problem consists of two horizontal layers, equal in height, in a rectangle, with a density difference of $\Delta \rho = \rho - \rho_0 \ge 0$, where $\rho_0$ is the density of the upper layer.
The initial condition for the temperature is a perturbation from the well-known static temperature field, connecting the temperature boundary conditions $T_0$ at the top of the rectangle and $T_1$ at the bottom of the rectangle~\cite{DLT-GS:2014}.
We study of a range of density differences $\Delta \rho$ that we characterize by the
non-dimensional buoyancy number $\mathrm{B}$, which is the ratio of $\Delta \rho $ to $\rho_0  \, \alpha \, \Delta T$, where $\Delta T = T_1 - T_0$, and $\alpha$ is the volumetric coefficient of thermal expansion.
The temperature perturbation initially drives the convection and, depending on the value of $\mathrm{B}$, determines the dynamics and structure of the resulting flow field.

In Section~\ref{SECTION: THERMOCHEMICAL CONVECTION WITH DENSITY STRATIFICATION} we begin by describing the equations that govern thermochemical convection in the mantle and the modification to these equations that we use to model density stratification in such flows. 
Then, in Section~\ref{SECTION:THE NUMERICAL METHODOLOGY} we describe the numerical methodology, including the underlying Finite Element Method (FEM) and the coupling of our VOF method to this FEM.
In Section~\ref{SECTION:NUMERICAL RESULTS} we begin by presenting two standard interface tracking benchmarks in order to demonstrate the accuracy of the VOF method. 
We then present the results of two benchmarks commonly used by researchers in the computational mantle convection community.
At the end of this section we present computational results of a model problem first proposed in~\cite{EGP-DLT-YH-HL-JM-LHK:2017}, which is designed to provide insight into how thermal plumes, that are thought to eventually reach the Earth's surface as ocean island basalts, originate at structures on the core-mantle boundary known as Large Low Shear wave Velocity Provinces (LLSVPs).
We briefly discuss these latter computational results in Section~\ref{Section:Discussion} and, in Section~\ref{Section:Conclusions}, we present our conclusions.



\vskip 48pt 

\begin{table}[h!]
    \centering
    \caption{A list of symbols used in this paper.}
    \scalebox{0.75}{
        \begin{tabular}{lll || lll}
            \hline
            Symbol     & Quantity                  & Unit             & Symbol        & Quantity                  & Unit \\
            \hline
            $\bf{u}$   &Velocity                   & $m / s$               & $\rho$        &Density                    & $kg\cdot m^{-3}$\\
            $p$        & Dynamic pressure           & Pa               &  $\Delta\rho$ &Density difference   & $kg\cdot m^{-3}$    \\
            $T_0$      &Temperature at the top     & K                & $D$           & Compositional diffusivity  &  $m^2/s$              \\
            $T_1$      &Temperature at the bottom  & K                & $\alpha$      & Thermal expansion coefficient &  1/K               \\
            $T$        &Temperature                & K                & $d$           & Vertical height of fluid layer & m \\
            $\Delta T$ &Temperature difference     & K                & $\mathrm{Pr}$ &Prandtl number    & $\frac{\mu}{\rho\kappa}$    \\
            $C$        &Composition                & -                & $\mathrm{Le}$ &Lewis number     & $\frac{\kappa}{D}$    \\
            $\mu$      &Viscosity                  & Pa $\cdot$ s     & $\mathrm{Ra}$ &Rayleigh number  & $\frac{\rho_0{\bf{g}}\alpha\Delta Td^3}{\mu\kappa}$\\
            $\kappa$   &Thermal diffusivity        & $m^2/s$          & $\mathrm{B}$  &Buoyancy ratio & $\frac{\Delta\rho}{\rho_0\alpha\Delta T}$         \\ 
            $\rho_0$   &Reference density          & $kg\cdot m^{-3}$ &  & & \\    
            \hline
        \end{tabular}
    }
    \label{Table:List of Symbols}
\end{table}


\section{Thermochemical Convection with Density Stratification}
\label{SECTION: THERMOCHEMICAL CONVECTION WITH DENSITY STRATIFICATION}

In this section we present in detail the equations associated with the model problem, which we briefly described above.

\subsection{The Dimensional Form of the Equations}

In order to study the efficacy of our implementation of a VOF algorithm in ASPECT to model processes that occur in the Earth's mantle, we compute a problem that emphasizes the effect of a compositional density difference on thermal convection.
We consider a two-dimensional flow in a horizontal fluid layer with a thickness or height $d$.
Our problem domain $\Omega$ has width $3 \, d$ and height $d$.
At a given reference temperature $T_0$ the region  $d / 2 < y \le d$ has a compositional 
density of $\rho_0$ and the region $0 \le y < d / 2$ has a compositional density of 
$\rho_0 + \Delta \rho$ where $\Delta \rho \, \ll \, \rho_0$.

We also introduce a composition variable $C(x, y, t)$ defined  by
\begin{equation}
  \label{DEF:Composition Variable C}
    C = \frac{\rho - \rho_0}{\Delta \rho} \, .       
\end{equation}
The composition $C$ is the concentration of the dense fluid as a function of space and time.
The initial condition for $C$ is 
\begin{equation}
  \label{EQ:Initial Conditions for C}
    C( x, y, t=0) \, = \,
      \left\{
        \begin{array}{ll}
          1 & \text{for} \quad 0     \leq y \leq d / 2 \, , \\
          0 & \text{for} \quad d / 2   <  y \leq d     \, .
       \end{array}
     \right.
\end{equation}

The upper boundary, at $y = d$, has temperature $T_0$ and the lower boundary at
$y = 0$ has temperature $T_1$.
The fluid is assumed to have a constant viscosity $\mu$ which is large. 
The Prandtl number is assumed to be very large,
\begin{equation}
    \mathrm{Pr} \, = \, \frac{\mu}{\rho_0 \kappa} \, \gg \, 1 \, ,
\end{equation}
where $\kappa$ is the thermal diffusivity so that inertial effects can be neglected.
The fluids in the high density and low density layers are immiscible; i.e., they cannot 
mix by diffusion. 
Similarly the Lewis number is also assumed to be large,
\begin{equation}
    \mathrm{Le} \, = \, \frac{\kappa}{D} \, \gg \, 1 \, ,
.2.1\end{equation}
where $D$ is the diffusion coefficient for the compositional variable $C$. 
Thus, the discontinuous boundary between the high density and low density fluids is 
preserved indefinitely.  

The problem we have posed requires the solution of the standard equations for thermal 
convection with the addition of an equation for the compositional field $C$ that tracks the density difference.
The governing equations are described in detail in~\cite{GS-DLT-PO:2001,DLT-GS:2014}.

We make the assumption that the Boussinesq approximation 
\begin{equation}
\label{EQ:Dimensional Boussinesq approximation}
\rho( x, y, t) \, = \, \rho_0 \, ( 1 - \alpha (T - T_0)) \, + \, \Delta \rho \, C \, .
\end{equation}  
holds; namely, that density differences associated with convection
$\rho_0 \alpha(T_1 - T_0)$ and $\Delta\rho$ are small compared with the reference density $\rho_0$.

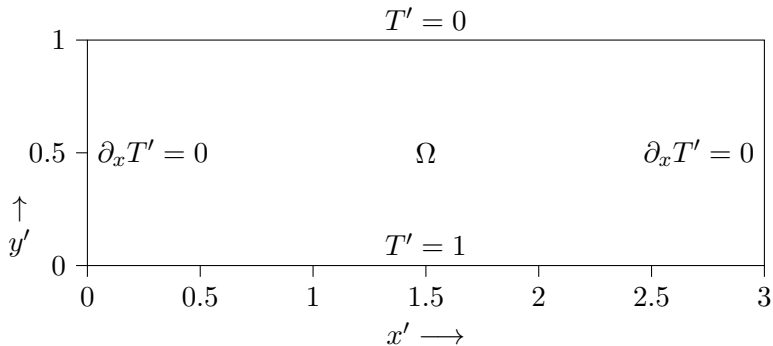
\begin{figure}[h!]
  \centering
    \begin{tikzpicture}[scale=3]
    \draw (0,0) rectangle (3,1);
    \foreach \x in {0, 0.5, 1, 1.5, 2, 2.5, 3}
    \draw (\x, 0) -- (\x, -0.05) node[anchor=north] {$\x$};
    \foreach \y in {0, 0.5, 1}
    \draw (0, \y) -- (-0.05, \y) node[anchor=east] {$\y$};
    \node at (1.5, 0.5) {$\Omega$};
    \node at (1.5, -0.3) {$x^{\prime} \longrightarrow $};
    \node at (-0.3, 0.1) {$y^{\prime}$};
    \node at (-0.3, 0.25) {$\uparrow $};
    \node[anchor=west]  at (0, 0.5) {$\partial_xT^{\prime} = 0$};
    \node[anchor=east]  at (3, 0.5)  {$\partial_xT^{\prime} = 0$};
    \node[anchor=south] at (1.5, 1)  {$T^{\prime} = 0$};
    \node[anchor=south] at (1.5, 0) {$T^{\prime} = 1$};
    \end{tikzpicture}
    \caption{The geometry of the (nondimensional) computational domain $\Omega$ shown with 
        the temperature boundary conditions on the four side walls. 
        The velocity boundary conditions on the side walls are
        $\mathbf{u} \cdot \mathbf{n}=0$ (no flow) and
        $\partial \bf{u} / \partial \boldsymbol{\tau} = 0$ (free slip) where $\mathbf{n}$ 
        and $\boldsymbol{\tau}$ are the unit normal and tangential vectors to the boundary 
        respectively.
            }
  \label{Fig:The Nondimensional Computational Domain}
\end{figure}

Conservation of mass requires
\begin{equation}
  \label{EQ:Dimensional Conservation of Mass}
    \frac{\partial u}{\partial x} + \frac{\partial v}{\partial y} = 0
\end{equation}  
where $x$ and $y$ denote the horizontal and vertical spacial coordinates, oriented as shown 
in Figure~\ref{Fig:The Nondimensional Computational Domain}, and $u$ and $v$ denote the horizontal and 
vertical velocity components, respectively.
We use the Stokes equations
\begin{align}
    \label{EQ:Dimensional Stokes Equation for u}
       0 &= \frac{-\partial P}{\partial x} 
            + \mu(\frac{\partial^2 u}{\partial x^2}
            + \frac{\partial^2 u}{\partial y^2}) \, , \\
    \label{EQ:Dimensional Stokes Equation for v}
      0 &= \frac{-\partial P}{\partial y}
           + \mu(\frac{\partial^2 v}{\partial x^2}
           + \frac{\partial^2 v}{\partial y^2}) 
           + \rho_0 \alpha (T - T_0) g - \Delta\rho C g \, ,
\end{align}
where $\alpha$ is the coefficient of thermal expansion, $g$ is the gravitational acceleration in the negative (downward) $y$ direction as shown in
Figure~\ref{Fig:The Nondimensional Computational Domain}, and
\begin{equation*}
   \label{EQ:P}
    P \, = \, p + \rho_0 \, g \, y
\end{equation*}  
where $p$ is the dynamic pressure and $ \rho_0 \, g \, y$ is the isostatic pressure.
Conservation of energy requires 
\begin{equation}
  \label{EQ:Dimensional Temperature Equation}
    \frac{\partial T}{\partial t} + u\frac{\partial T}{\partial x}
      + v\frac{\partial T}{\partial y} 
      = \kappa(\frac{\partial^2 T}{\partial x^2} + \frac{\partial^2 T}{\partial y^2})\,.
\end{equation}

With no diffusion, i.e., $D = 0$, the composition variable $C$ satisfies the advection equation
\begin{equation}
  \label{EQ:Dimensional Composition Equation}
    \frac{\partial C}{\partial t}
      + u\frac{\partial C}{\partial x} + v\frac{\partial C}{\partial y}
      = 0\,.
\end{equation}

\subsection{The Nondimensional Form of the Equations}
We introduce the non-dimensional variables 
\begin{equation}
  \begin{aligned}
    \label{EQ:Nondimensionalized Variables}
    x' &= \frac{x}{d},               & y' &= \frac{y}{d},            & t' &= \frac{\kappa}{d^2} \, t, \\
    u' &= \frac{d}{\kappa} \, u,     & v' &= \frac{d}{\kappa} \, v,      & \rho' &= \frac{\rho}{\rho_0},   \\
    T' &= \frac{T - T_0}{T_1 - T_0}, & P^\prime &= \frac{d^{\, 2} P}{\mu \kappa},
  \end{aligned}
\end{equation}
and the two nondimensional parameters, the Rayleigh number $\mathrm{Ra}$ and the buoyancy ratio $\mathrm{B}$
\begin{align}
  \label{DEF:Rayleigh Number}
    \mathrm{Ra} \, &= \, \frac{\rho_0 \, g \, \alpha(T_1 - T_0)d^3}{\mu \, \kappa} \, , \\
  \label{DEF:Buoyancy Ratio}
    \mathrm{B}  \, &= \, \frac{\Delta\rho}{\rho_0 \alpha(T_1 - T_0)}.
\end{align}
Substitution of equations~\eqref{EQ:Nondimensionalized Variables}--\eqref{DEF:Buoyancy Ratio} 
into
equations~\eqref{EQ:Dimensional Conservation of Mass}--\eqref{EQ:Dimensional Composition Equation}
gives
\begin{align}
   \label{EQ:Dimensionless Continuity Equation}
    \frac{\partial u'}{\partial x'} + \frac{\partial v'}{\partial y'}
    &= 0                                                                              \, , \\
  \label{EQ:Dimensionless Stokes Equations For u}
    0 &= \frac{-\partial P'}{\partial x'}
       + \frac{\partial^2 u'}{\partial x^{'2}}
       + \frac{\partial^2 u'}{\partial y^{'2}}                                        \, , \\
  \label{EQ:Dimensionless Stokes Equations For v}
    0 &= \frac{-\partial P'}{\partial y'}
       + \frac{\partial^2 v'}{\partial x^{'2}}
       + \frac{\partial^2 v'}{\partial y^{'2}}
       + \mathrm{Ra} T' - \mathrm{Ra} \mathrm{B} C                                    \, , \\
  \label{EQ:Dimensionless Temperature Advection Diffusion Equation}
    \frac{\partial T'}{\partial t'}
      + u' \frac{\partial T'}{\partial x'} 
      + v' \frac{\partial T'}{\partial y'}
     &= \frac{\partial^2 T'}{\partial x^{'2}} + \frac{\partial^2 T'}{\partial y^{'2}} \, , \\
  \label{EQ:Dimensionless Composition Equation}
    \frac{\partial C}{\partial t'} + u'\frac{\partial C}{\partial x'}
                                   + v'\frac{\partial C}{\partial y'}
     &= 0 \, .
\end{align}

This is the superposition of a Rayleigh-Taylor problem and a Rayleigh-B\'{e}nard problem~\cite{SC:1961,DLT-GS:2014}.
In the isothermal limit, $T_0 = T_1$, it is the classic Rayleigh-Taylor problem.
If $C$ is positive, a light fluid is above the heavy fluid and in a downward gravity field the fluid layer is stable. 
If $\Delta \rho$ is negative, a heavy fluid lies over a light fluid and the layer is 
unstable. 
Flows will transfer the heavy fluid to the lower half and the light fluid to the upper 
half and the density layer will overturn. 
If $\Delta \rho = 0$ and hence, $\mathrm{B} = 0$, this is the classic Rayleigh-B\'{e}nard problem for thermal convection. 
The governing parameter is the Rayleigh number $\mathrm{Ra}$. 
If $0 < \mathrm{Ra} < \mathrm{Ra}_{c}$, the critical Rayleigh number, no flow will occur. 
If $\mathrm{Ra}_c < \mathrm{Ra} < \mathrm{Ra}_t$ steady cellular flow will occur. 
If $\mathrm{Ra} > \mathrm{Ra}_t$ the flow becomes unsteady and thermal turbulence 
develops.  

If $\mathrm{Ra} > \mathrm{Ra}_c$ and $\mathrm{B}$ is small, the boundary between the density differences will not block the flow driven by thermal convection.
Kinematic mixing will occur and the composition will homogenize so that the density is 
constant. 
Whole layer convection will occur.
If $\mathrm{B}$ is large, the density difference boundary will block the flow driven by thermal convection.
The compositional boundary will be displaced vertically but will remain intact.  Layered 
convection will occur with the compositional boundary, the boundary between the convecting  layers.  
In this work the Rayleigh number $\mathrm{Ra}$ defined in
equation~\eqref{DEF:Rayleigh Number} is based on the domain thickness $d$ and this is the 
case for which we will show numerical computations.


\section{The Numerical Methodology}
\label{SECTION:THE NUMERICAL METHODOLOGY}

In the following discussion of the numerical methodology, we will only consider the
dimensionless
equations~\eqref{EQ:Dimensionless Continuity Equation}-\eqref{EQ:Dimensionless Composition Equation}
and drop the primes associated with the dimensionless variables.
The vector form of the dimensionless equations on the 2D rectangular domain  
$\Omega=[0,3]\times[0,1]$ shown in Figure~\ref{Fig:The Nondimensional Computational Domain} are given by
\begin{align}
\label{EQ:Nondimensional Vector Form of the Stokes Equations}
-\nabla^2 {\bf{u} } + \nabla  P   
&= (-\mathrm{Ra} \ T + \mathrm{Ra}  \ \mathrm{B} \ C) \, {\bf{g}} \ , \\
\label{EQ:NONDIMENSIONAL VECTOR FORM OF THE CONTINUITY EQUATION}
\nabla\cdot{\mathbf {u}}
&= 0 \ ,                                                            \\
\label{EQ:Nondimensional Vector Form of the Temperature Equation}
\frac{\partial T}{\partial t}+{\bf u}\cdot\nabla T
&= \nabla^2 T \,                                                    \\
\label{EQ:NONDIMENSIONAL VECTOR FORM OF THE COMPOSITION EQUATION}
\frac{\partial C}{\partial t}+{\bf u}\cdot\nabla C
&=0,
\end{align}
where  ${\bf u}=(u,v)$ is the velocity and ${\bf g}= (0, -1)$ is the unit vector pointing 
downward.

Note that the composition 
equation~\eqref{EQ:NONDIMENSIONAL VECTOR FORM OF THE COMPOSITION EQUATION} is equivalent to
\begin{equation}
\label{EQ: MATERIAL DERIVATIVE FORM OF THE COMPOSITION EQUATION}
\frac{D C}{D t} =
\frac{\partial C}{\partial t} + u \, \frac{\partial C}{\partial x}
+ v \, \frac{\partial C}{\partial y} = 0 \, ,
\end{equation}
where
\begin{equation}
\label{DEF: DEFINITION OF THE MATERIAL DERIVATIVE}
\frac{D}{D t} \equiv
\frac{\partial}{\partial t}  + u \, \frac{\partial }{\partial x}  
+ v \, \frac{\partial }{\partial y} 
\end{equation}
is the \textit{material derivative}.
Equation~\eqref{EQ: MATERIAL DERIVATIVE FORM OF THE COMPOSITION EQUATION} implies that the 
composition $C$ is constant on particle paths in the flow~\citep{AJC-JEM:1993}.
Furthermore, since by~\eqref{EQ:Dimensionless Continuity Equation} 
the velocity ${\mathbf {u}}$ is divergence free, the composition
equation~\eqref{EQ:NONDIMENSIONAL VECTOR FORM OF THE COMPOSITION EQUATION} can be written 
in \textit{conservation form}
\begin{equation}
\label{EQ: CONSERVATION FORM OF THE COMPOSITION EQUATION}
\frac{\partial C}{\partial t} + \nabla \cdot \left({\mathbf u} \, C\right) = 0 \, ,
\end{equation}
implying that the composition $C$ is a conserved quantity - it is neither created nor 
destroyed as it is advected in the flow field.

We assume no-flow and free-slip velocity boundary conditions on all boundaries,
\begin{align}
\label{EQ:No-FLow BCs}
\mathbf{u} \cdot \mathbf{n}                         &= 0 
\qquad \text{(no-flow)}   \, , \\
\label{EQ:Free-Slip BCs}
\frac{\partial \bf{u}}{ \partial \boldsymbol{\tau}} &= 0
\qquad \text{(free slip)} \, ,
\end{align}
where $\mathbf{n}$ and $\boldsymbol{\tau}$ are the unit normal and tangential vectors to 
the boundary respectively.
We impose Dirichlet boundary conditions for the temperature on the top and bottom of the 
computational domain and Neumann boundary conditions (no heat flux)  on the sides of the 
computational domain,
\begin{eqnarray}
  \label{EQ:Temperature BCs at y = 0} 
    T ( x, 0, t) = 1             \, ,   \\
  \label{EQ:Temperature BCs at y = 1}
    T ( x, 1, t) = 0             \, ,  \\
  \label{EQ:Temperature BCs at x = 0}
    \partial_{x} T ( 0, y, t) = 0  \, , \\
  \label{EQ:Temperature BCs at x = 3}
    \partial_{x} T ( 0, y, t) = 0 \, .
\end{eqnarray}
The geometry of the computational domain together with the boundary conditions 
on the temperature are shown in Figure~\ref{Fig:The Nondimensional Computational Domain}.

\subsection{Decoupling of the Nonlinear System}
\label{Cubsection:Nonlinear Decoupling}
The incompressible Stokes equations can be considered as a constraint on the temperature and 
composition at any given time leading to a highly nonlinear system of equations.
To solve this nonlinear system, we apply the Implicit Pressure Explicit Saturation (IMPES) approach, originally developed for computing solutions of equations for modeling problems in porous media flow~\citep{RH-RH:1999,JWS-WTC:1959}, to decouple the incompressible Stokes 
equations~\eqref{EQ:Dimensionless Continuity Equation}--\eqref{EQ:Dimensionless Stokes Equations For v}
from the temperature and compositional equations
\eqref{EQ:Dimensionless Temperature Advection Diffusion Equation}--\eqref{EQ:Dimensionless Composition Equation}. 
This leads to three discrete systems of linear equations, the Stokes equations, the temperature equation, and the composition equation, thereby allowing each equation to be solved easily and efficiently. 

\subsection{Discretization of the Stokes Equations}
Let $t^k$ denote the discretized time at the k{\it th} time step with a time step size of
$\Delta t^k =t^k - t_{k-1}$, $k = 0, \, 1, \, \ldots$ 
Given the temperature $T^k$ and composition $C^k$ at time $t = t^k$, we first solve for our approximation to the Stokes
equations
\eqref{EQ:Dimensionless Continuity Equation}--\eqref{EQ:Dimensionless Stokes Equations For v} to obtain the 
velocity  $\mathbf{u}^k = (u^k , \, v^k)$ and pressure $P^k$
\begin{align}
  \label{EQ:Stokes Discretization a}
    -\nabla^2 {\bf{u} }^k +\nabla  P^k
     &= (-\mathrm{Ra} \ T^k + \mathrm{Ra} \ \mathrm{B} \, C^k) \ {\mathbf{g}} \ , \\
  \label{EQ:Stokes Discretization b}
     \nabla \cdot {\mathbf{u}}^k 
    &= 0 \ . 
\end{align}
For the incompressible Stokes
equations~\eqref{EQ:Stokes Discretization a}--\eqref{EQ:Stokes Discretization b}, we use the standard mixed FEM method with a Taylor-Hood element \cite{JD-AH:2005} for the spatial approximation.
We refer the interested reader to \cite{MK-TH-WB:2012} for a more detailed discussion of the spatial discretization and the choice of Stokes preconditioners and solvers.

\subsection{Discretization of the Temperature Equation}
In all of the computations presented here we use the algorithm currently implemented in ASPECT~\cite{ASPECT-v1_5_0:2017} to approximate the spatial and temporal terms in the temperature equation~\eqref{EQ:Nondimensional Vector Form of the Temperature Equation}. 
This algorithm includes an entropy viscosity stabilization technique described 
in~\cite{JLG-RP-BP:2011,MK-TH-WB:2012}.
If we introduce the inner product of two scalar functions $u$ and $v$ on $\Omega$
\begin{align}
(u,v)_{\Omega} = \int\limits_\Omega \, u \,v \, dx \, dy
\end{align}
and $\Gamma_D=\{y=0\}.$
By multiplying the test function $\psi(x,y)$ and taking the integration, the weak form of this spatial discretization is
\begin{align}
\label{EQ:WEAK FORM OF THE SPATIAL DISCRETIZATION OF THE TEMPERATURE EQUATION}
(\frac{\partial T}{\partial t},\psi)_{\Omega}+({\bf u}\cdot\nabla T,\psi)_{\Omega}
&=-(\nabla T,\nabla\psi)_{\Omega}-(\nu_h(T)\nabla T,\nabla\psi)_{\Omega}  
+ (\frac{\partial T}{\partial {\bf n}},\psi)_{\Gamma_D}
\end{align}
where  $\nu_h(T)$ is the entropy viscosity function as defined in~\cite{MK-TH-WB:2012}, 
except here we do not use a second-order extrapolation to treat the advection term
$({\bf u}\cdot \nabla T,\psi)$ and the entropy viscosity term  
$(\nu_h(T)\nabla T,\nabla\psi)_{\Omega}$ explicitly.
We use the fully implicit adaptive Backward Differentiation Formula of order 2 (BDF2) 
\citep{GW-EH:1991, MK-TH-WB:2012} to discretize the temperature equation in time.
Thus, the full discretization of the temperature equation is
\begin{align}
\label{EQ:Weak Form of the Discretization of the Temperature Equation}
&(\frac{1}{\Delta t^{k+1}}\left (\frac{2\Delta t^{k+1}+\Delta t^{k}}{\Delta t^{k+1}+\Delta t^{k}}T^{k+1}-\frac{\Delta t^{k+1}+\Delta t^{k}}{\Delta t^{k}}T^k+\frac{(\Delta t^{k+1})^2}{\Delta t^{k}(\Delta t^{k+1}+\Delta t^{k})}T^{k-1}\right ),\psi)_{\Omega}
\nonumber \\
=&-({\bf u}^{k}\cdot\nabla T^{k+1},\psi)_{\Omega}-(\nabla T^{k+1},\nabla\psi)_{\Omega}-(\nu_h^k(T)\nabla T^{k+1},\nabla\psi)_{\Omega}   
+(\frac{\partial T^{k+1}}{\partial {\bf n}},\psi)_{\Gamma_D}\,.
\end{align}
The entropy-viscosity function $\nu_h^k(T)$ is a non-negative constant within each cell that only 
adds artificial diffusion in cells for which the local P\'{e}clet number 
$\mathrm{Pe} = \mathrm{Ra} \cdot \mathrm{Pr}$ is large and the solution is not smooth.

\subsection{Discretization of the Composition Equation}
\label{SUBSECTION:DISCRETIZATION OF THE COMPOSITION EQUATION}

In this article we use the Volume-of-Fluid (VOF) interface tracking algorithm described in Section~\ref{SUBSECTION:THE VOLUME-OF-FLUID INTERFACE TRACKING ALGORITHM} below to discretize the composition 
equation~\eqref{EQ:NONDIMENSIONAL VECTOR FORM OF THE COMPOSITION EQUATION}.
Prior to our implementation of the VOF method in ASPECT the only algorithms one could use to model the solution
of~\eqref{EQ:NONDIMENSIONAL VECTOR FORM OF THE COMPOSITION EQUATION} were based on a spatial discretization of the weak form of the of the composition equation,
\begin{align}
  \label{EQ:WEAK FORMULATION OF THE COMPOSITION EQUATION}
    (\frac{\partial C}{\partial t},\psi)_{\Omega} + ({\bf u}\cdot\nabla C,\psi)_{\Omega} 
      \, = \, 0 \, .
\end{align}

The first advection algorithm that was implemented in ASPECT is based on the same spatial discretization as in
equation~\eqref{EQ:WEAK FORM OF THE SPATIAL DISCRETIZATION OF THE TEMPERATURE EQUATION}.
However, the entropy-viscosity stabilization term on the right-hand side in 
\begin{align}
  \label{EQ:WEAK FORMULATION OF THE COMPOSITION EQUATION WITH ENTROPY VISCOSITY}
    (\frac{\partial C}{\partial t},\psi)_{\Omega} + ({\bf u}\cdot\nabla C,\psi)_{\Omega} 
       \, = \, -(\nu_h(C) \nabla C, \nabla \psi)_{\Omega}
\end{align}
is computed separately for the composition field; i.e, it does \textit{not} have the same value in each cell as does the entropy viscosity function $\nu_h(T)$ for the temperature  field.
The adaptive BDF2 algorithm is also used for the time discretization of the composition equation, leading to the following FEM Entropy Viscosity (FEM-EV) discretization
of equation~\eqref{EQ:NONDIMENSIONAL VECTOR FORM OF THE COMPOSITION EQUATION},
\begin{align}
  \label{EQ:THE ENTIRE FEM-EV DISCRETIZATION}
    \frac{1}{\Delta t^{k+1}} ( & \frac{2\Delta t^{k+1}+\Delta t^{k}}{\Delta 
      t^{k+1}+\Delta t^{k}}C^{k+1}-\frac{\Delta t^{k+1}+\Delta t^{k}}{\Delta t^{k}}C^k+\frac{(\Delta t^{k+1})^2}{\Delta t^{k}(\Delta t^{k+1}+\Delta t^{k})}C^{k-1},\psi )_{\Omega} \nonumber\\
       =& - ({\bf u}^{k}\cdot\nabla C^{k+1},\psi)_{\Omega}-(\nu_h^k(C)\nabla C^{k+1},\nabla\psi)_{\Omega} \, .
\end{align}
In equation~\eqref{EQ:THE ENTIRE FEM-EV DISCRETIZATION} the entropy viscosity 
function $\nu_h^k(C)$ has the same purpose as $\nu_h^k(T)$.

The other algorithm for modeling solutions 
of~\eqref{EQ:NONDIMENSIONAL VECTOR FORM OF THE COMPOSITION EQUATION} that is implemented in ASPECT is a Discontinuous Galerkin method with a Bound Preserving limiter.
See~\cite{YH-EGP-MIB:2017} for a description of this algorithm in ASPECT and a comparison with the advection method with entropy viscosity described above.
Also see~\cite{EGP-DLT-YH-HL-JM-LHK:2017} for a comparison of these two methods with the VOF method described here and a particle method for modeling the solution of the composition equation~\eqref{EQ:NONDIMENSIONAL VECTOR FORM OF THE COMPOSITION EQUATION}.

\subsection{The Volume-of-Fluid Interface Tracking Algorithm}
\label{SUBSECTION:THE VOLUME-OF-FLUID INTERFACE TRACKING ALGORITHM}

The Volume-of-Fluid (VOF) method is an interface tracking method in which, at each time step, the interface between the two compositions, one ($C \, = \, 1$) with density 
$\rho = \rho_0 + \Delta \rho$ and the other ($C \, = \, 0$) with density $\rho = \rho_0$, is explicitly reconstructed in every cell that contains a portion of the interface.
Given this explicit (approximate) location of the interface at the current time step one then uses this information to advance the interface in time. 
In this sense the VOF method approximates the compositional interface on a subgrid scale.
In addition, in an incompressible flow both the VOF \textit{interface reconstruction algorithm} and the VOF \textit{advection algorithm} presented here conserve the
volume of each of the two compositions throughout the course of the computation.

\subsubsection{Background}
\label{SUBSUBSECTION: BACKGROUND}

There are a wide variety of possible VOF interface reconstruction and advection
algorithms; e.g., see ~\cite{JEP-EGP:2004} and the references there.
The VOF method was first developed at the U.S.~National Labs in the 1970s~\cite{WFN-PRW:1976} and have continued to be used and developed by researchers at the National
Labs~\cite{CWH-BDN:1981,BDN-CWH-RSH:1980,MDT-RCM-LRS:1987,MDT-LDC-RCM-CWH:1985} as well as around the world.
The advantage that VOF methods have over other interface tracking algorithms is that
they are designed to naturally satisfy a conservation law; namely,
equation~\eqref{EQ:CONSERVATION EQUATION FOR THE VOLUME FRACTION f} below.
Thus, materials that should be conserved as they move with the flow are conserved,
without the need to resort to additional numerical algorithms such as the
redistancing step in a Level Set method~\cite{JAS:1999}.
VOF methods can and have been used effectively to model a wide variety of moving
interface problems, including interfaces in compressible flow with shock
waves~\cite{LFH-PC-EGP:1991}, interfaces with shock waves in materials in the limit of no 
strength effects~\cite{GHM-EGP:1994,GHM-EGP:1996}, jetting in meteorite impacts~\cite{EGP-GHM:1996}, \textit{nonconservative} interface
motion such as photolithography~\cite{JJH-PC-EGP:1997,JJH-PC-EGP-MD:1996}, the transition 
from deflagration to detonation~\cite{JEP-EGP:1998} and more than two materials; i.e., 
more than one interface in a cell~\cite{HRA-KM:2011,RNH-MS:2013}.

\subsubsection{Description}
\label{SUBSUBSECTION:DESCRIPTION}
In this article we use a two-dimensional VOF algorithm to discretize the conservation equation
\begin{equation}
  \label{EQ:CONSERVATION EQUATION FOR THE VOLUME FRACTION f}
    \frac{\partial f}{\partial t} \, + \, \nabla \cdot \mathbf{F} \left( f \right) \,
     \, = \, 0 \, ,
\end{equation}
where $\mathbf{u} = (u, \, v)$ is the velocity field, $f$ is the
\textit{volume fraction} of one of the compositional fields, say
$C \, = \, 1$, the field with density $\rho_0 \, + \, \Delta \rho$, which we will refer to as `Composition~1' or $C_1$ for short,\footnote{Throughout this section and beyond we will use the terms ``volume'' and ``volume fraction'' of $C_1$, etc., although it is to be understood that in two dimensions the quantity in question is an area.} and
\begin{equation}
  \label{DEF:THE VOLUME FRACTION FLUX}
    \mathbf{F}\left( f \right)
     \, = \, \left( F \left( f \right) , G (f) \right)
     \, = \, \left( u \, f , \, v \, f \, \right)
\end{equation}
is the volume fraction flux associated with $C_1$.
In our VOF implementation in ASPECT we use the `Efficient Least Squares VOF Interface
Reconstruction Algorithm' (ELVIRA), which is described in detail in~\cite{JEP-EGP:2004}
and is based on the ideas in~ \cite{JEP:1992} and~\cite{EGP:1991}.
The ELVIRA interface reconstruction algorithm reconstructs lines on a uniform grid with
square cells \textit{exactly}.
We will explain this in more detail in
Section~\ref{SUBSUBSECTION:THE ELVIRA INTERFACE RECONSTRUCTION ALGORITHM} and give a
example in Section~\ref{SUBSUBSECTION:ADVECTION OF A LINEAR INTERFACE} below.
Since the ELVIRA algorithm reconstructs lines in square cells exactly it is natural to
assume that the algorithm is second-order accurate on a uniform grid with identical
square cells.
This turns out to be true \cite{EGP:2010a,EGP:2010b,EGP:2014}.
We use a second-order accurate operator splitting advection
method~\cite{JEP-EGP:2004,EGP-DLT-YH-HL-JM-LHK:2017,WGS:1968} to update the values of the volume fractions in time.

For simplicity of exposition we will assume the finite element grid consists entirely of
square cells $\Omega_e$, of side $h$, indexed by the variable $e$, and aligned parallel to the $x$ and $y$ axes.
The discretization of equation~\eqref{EQ:CONSERVATION EQUATION FOR THE VOLUME FRACTION f}
proceeds as follows.
Let $\Omega_e$ denote an arbitrary finite element cell in our domain $\Omega$
and let $f^k_e$ denote the discretized volume fraction in $\Omega_e$ at time $t^k$.
The variable $f^k_e$ is a scalar that satisfies $0 \le f^k_e \le 1$ such that
\begin{equation}
  \label{EQ: THE VOLUME OF COMPOSITION 01 IN Omega_e}
    f^k_e \, \approxeq \,
         \frac{1}{h^2} \int \limits_{\Omega_e} \, f(x, \,  y, \, t^k \, ) \, dx \, dy \, .
\end{equation}
Thus, the discretized volume, $V^k_e$, of $C_1$ in $\Omega_e$ at time $t^k$ is
\begin{equation}
  \label{EQ:The Volume of C_1 in Omega_e}
    V^k_e \,
      = \, \int\limits_{\Omega_e} \, f^k_e \; dx \, dy \,
      = h^2 \, f^k_e \, .
\end{equation}
Note that for an incompressible velocity field $\mathbf{u} = ( u,v)$ we have
$\nabla \cdot \mathbf{u} = 0$ and hence, the volume of `parcels' or regions of $C_1$ are
constant as they evolve in time.

From a mathematical point of view the variable $f(x,y)$ is the characteristic function associated with $C_1$.
In other words,
\begin{equation}
  \label{DEF: f AS A CHARACTERISTIC FUNCTION}
  f(x,y) \, = \,
    \begin{cases}
      f(x,y) \, &= \, 1 \quad \text{if $(x,y)$ is occupied by Composition~1}     \, , \\
      f(x,y) \, &= \, 0 \quad \text{if $(x,y)$ is not occupied by Composition~1} \, .
    \end{cases}
\end{equation}
This implies $1 - f(x,y)$ is the characteristic function associated with $C_2$, the composition with density $\rho \, = \, \rho_0$.
In this article we restrict ourselves to modeling the interface between two compositions.
However, there is currently a great deal of research into modeling three or more
interfaces in one cell with a VOF method; e.g., see~\cite{MJ-MMS-MS:2015}.

In its simplest form our implementation of the VOF algorithm in ASPECT proceeds as follows.
  \begin{figure}[b!]
      \begin{center}
          \includegraphics[height=3.0in]{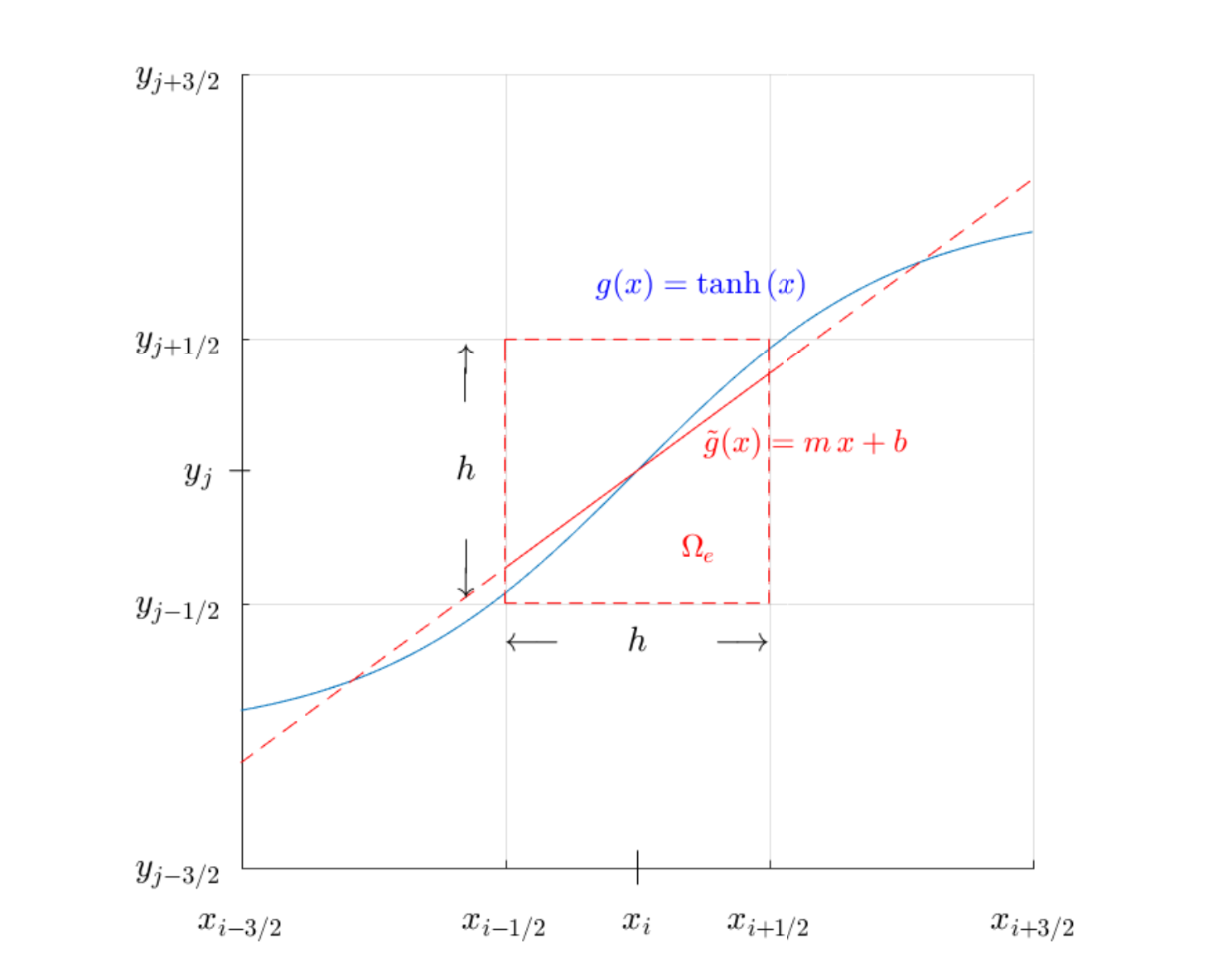}
          \caption{In our implementation of the VOF interface reconstruction algorithm
              the true interface, which in this example is $g(x) = \tanh (x)$, is
              approximated as a line segment $\tilde{g}_{e} (x) = m_{e} \, x + b_{e}$
              in each cell \textcolor{Red}{$\Omega_e$} that has a volume fraction
              $f_{e}$ with $0 < f_{e} < 1$.
              The approximate interface in \textcolor{Red}{$\Omega_e$} is depicted as the
              solid red line segment in the center cell \textcolor{Red}{$\Omega_e$}.
              In this example, as with all VOF methods, the volume $h^2 f_e^{true}$
              beneath the true interface in \textcolor{Red}{$\Omega_e$} is
              \textit{exactly equal} to the volume $h^2 f_e$ beneath the approximate
              interface $\tilde{g}$ in \textcolor{Red}{$\Omega_e$};
              i.e., $f_e^{true} = f_e$.  Note that, for convenience, we have used the notation $(x_i, y_j)$ to denote the center of the cell $\Omega_e$,
             $[x_{i-1/2}, \, x_{i+1/2}] \times [y_{i-1/2}, \, y_{i+1/2}]$ to denote the cell $\Omega_e$, etc.
            }
            \label{FIG: ELVIRA APPROXIMATION TO TANH}
        \end{center}
    \end{figure}

Given the values $f^k_e$ at time $t^k$ and the velocity field at time $t^k$ we do the
following to obtain the volume fractions $f^{k+1}_e$ at time $t^{k+1}$.
\begin{enumerate}

    \item \label{ITEM: THE INTERFACE RECONSTRUCTION STEP} THE INTERFACE RECONSTRUCTION
    STEP:
    Given a cell $\Omega_e$ that contains a portion of the interface, so
    $0 < f^k_e < 1$ where $f^k_e$ is the volume fraction in $\Omega_e$ at time
    $t^k$, use the volume fractions  $f^k_{e'}$ in the $3 \times 3$ block of cells
    $\Omega_{e'}$ centered on the cell $\Omega_e$ to \textit{reconstruct} the interface
    in $\Omega_e$.
    The reconstructed interface will be a piecewise linear approximation to the
    true interface as shown in Figure~\ref{FIG: ELVIRA APPROXIMATION TO TANH} that
    preserves the given volume $h^2 \, f^k_e$ of $C_1$ in $\Omega_e$.
    We give a brief description of how we determine the linear approximation
    $\tilde{g}_{e} (x) = m_{e} \, x \, + \, b_{e}$, to the true interface in cells $\Omega_e$ for which $0 < f^k_e < 1$ in
    Section~\ref{SUBSUBSECTION:THE ELVIRA INTERFACE RECONSTRUCTION ALGORITHM} below.

    \begin{figure}[h!]
        \centering
        \includegraphics[height=3.0in]{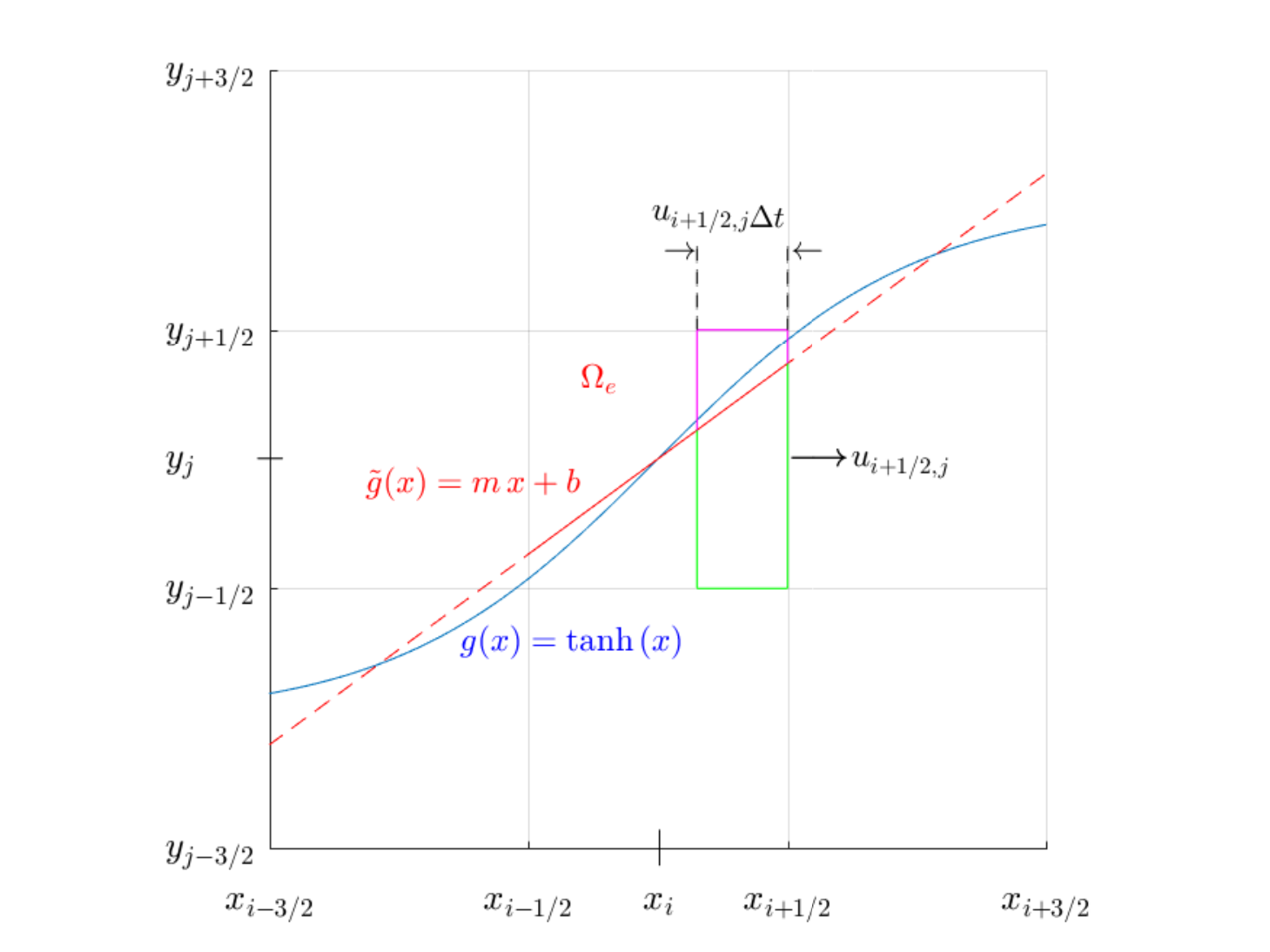}
        \caption{The volume $V^k_{i + 1/2, j}$ of $C_1$ in the quadrilateral 
            outlined in green on three sides and by a portion of the solid red line on top is the flux of of $C_1$ that will cross the right-hand edge of \textcolor{Red}{$\Omega_e$} during the time step from time $t^k$ to $t^{k+1}$.
            Here $\Delta t^k = t^{k+1} - t^k$ and we have dropped the superscript $k$
            from $u^k_{i+i/2,j}$ and $\Delta t^k$ for clarity.
            The solid red line in \textcolor{Red}{$\Omega_e$} is the reconstructed
            interface that approximates the true interface $g(x) = \tanh \, (x)$
            in \textcolor{Red}{$\Omega_e$} as shown in
            Figure~\ref{FIG: ELVIRA APPROXIMATION TO TANH}.
                }
        \label{FIG:COMPUTATION OF VOLUME FLUX ACROSS THE RH EDGE}
    \end{figure}

        \vskip 12pt

        \item \label{ITEM: COMPUTATION OF THE FLUXES} COMPUTATION OF THE FLUXES: In the
        computations presented in this article we use a second-order accurate operator
        split algorithm, often referred to as Strang Splitting~\cite{WGS:1968}, in order
        to advance the interface in time.
        However, as mentioned above, for clarity and simplicity of exposition we will  only describe a first-order accurate operator split VOF advection algorithm here.
        See~\cite{JEP-EGP:2004} for the details of a second-order accurate operator split VOF advection algorithm.

        \vskip 06pt

        \hskip 12pt For convenience and clarity of exposition, for the remainder of this
        section we will use the index notation $(i,j)$, as shown in
        Figures~\ref{FIG: ELVIRA APPROXIMATION TO TANH}--\ref{FIG: A LINE WITH TWO EXACT COLUMN SUMS}.
        Thus, we have nine cells with centers $(x_{i'}, y_{j'})$ for
        $i' = i-1, \, i, \, i+1$ and $j' = j-1, \, j, \, j+1$ with edges indexed as shown
        in the figure.
        In the ELVIRA interface reconstruction algorithm we use the information in the
        $3 \times 3$ block of cells $\Omega_{i'j'}$ immediately adjacent to the cell
        $\Omega_{ij}$ in which we wish to reconstruct the interface.
        Given the reconstructed interface $ \tilde{g}_{e} \, = \, \tilde{g}_{ij} (x)$ in
        \begin{equation}
            \Omega_{e} \equiv \Omega_{ij} = [x_{i-1/2},x_{i+1/2}] \times [y_{j-1/2},y_{i+1/2} ]
        \end{equation}
        as shown in Figure~\ref{FIG:COMPUTATION OF VOLUME FLUX ACROSS THE RH EDGE} and the
        velocity $u^k_{i \pm 1/2, j}$ normal to the right and left edges of
        $\Omega_{ij}$ at time $t^{k}$, we wish to determine the \textit{volumes}
        $V^k_{i \pm 1/2, j}$ of $C_1$ that cross the right and left edges of $\Omega_e$
        in the time interval $[t^{k}, t^{k+1}]$.
        These volumes are determined \textit{geometrically}.
        A diagram for how to determine the \textit{volume} $V^k_{i + 1/2, j}$ of $C_1$
        that crosses the right-hand edge of $\Omega_{ij}$ in the time interval
        $[t^{k}, t^{k+1}]$, given the assumption that $u^k_{i \pm 1/2, j} > 0$, is
        outlined in green on three sides and by a portion of the solid red line on top in Figure~\ref{FIG:COMPUTATION OF VOLUME FLUX ACROSS THE RH EDGE}.

        \vskip 12pt

        \item \label{ITEM: THE CONSERVATIVE UPDATE} THE CONSERVATIVE UPDATE:
        Given the \textit{volumes} $V^k_{i \pm 1/2, j}$ of $C_1$ that
        cross the left and right-hand edges of $\Omega_{ij}$ in the time interval
        $[t^{k}, t^{k+1}]$ we use the following equation to determine an
        \textit{intermediate} volume $\hat{V}^{k}_{ij}$ in $\Omega_{ij}$ for the first
        part of the two part operator split algorithm:
        \begin{equation}
            \label{EQ: CONSERVATIVE VOF UPDATE}
            \hat{V}^{k}_{ij} \,
            = \,  V^k_{ij} \, + \, V^k_{i - 1/2, j} - V^k_{i + 1/2, j} \, ,
        \end{equation}
        where  $V^k_{ij} = h^2 f^k_{ij}$ and $\hat{V}^{k}_{ij}$ denotes the
        `intermediate' volume in $\Omega_{ij}$ after the first part of the operator split
        advection step from time $t^{k}$ to time $t^{k+1}$.

        \vskip 12pt Given the nine intermediate volume fractions
        $\hat{f}^k_{i'j'} \equiv \hat{V}^k_{i'j'} / h^2$
        in $\Omega_{ij}$ and the $3 \times 3$ block of cells $\Omega_{i'j'}$ surrounding
        $\Omega_{ij}$, together with all of the intermediate volume fractions in the
        $3 \times 3$ block of cells surrounding each of the cells $\Omega_{i'j'}$,
        reconstruct an intermediate interface $\hat{g}_{i'j'}(x)$ in each cell
        $\Omega_{i'j'}$ and use it to geometrically determine the volumes
        $\hat{V}^{k}_{i, j \pm 1/2}$~of $C_1$ that cross the \textit{top} and
        \textit{bottom} edges of $\Omega_{ij}$ in the time interval $[t^{k}, t^{k+1}]$ in
        the same manner as illustrated in
        Figure~\ref{FIG:COMPUTATION OF VOLUME FLUX ACROSS THE RH EDGE}, but this time in the $y$-direction.
        The volume of $ V^{k+1}_{ij}$ in $\Omega_{ij}$ at the new time $t^{k+1}]$ is thus,
        \begin{equation}
          \label{EQ: COMPUTATION OF V^{k+1}_{ij}}
            V^{k+1}_{ij} \,
              = \, \hat{V}^{k}_{ij} \,
              + \, \hat{V}^k_{i, j - 1/2} - \hat{V}^k_{i, j + 1/2}
        \end{equation}
        The new volume fraction in $\Omega_{ij}$ is now
        \begin{equation}
            f^{k+1}_{ij} \, = \, V^{k+1}_{ij} / h^2
        \end{equation}

        There are also \textit{unsplit} VOF advection algorithms; e.g., see \cite{JEP-EGP:2004,EGP-ASA-JBB-DLM-WJR:1997}.

    \end{enumerate}

   \begin{figure}[h!]
        \centering
        \includegraphics[height=3.0in]{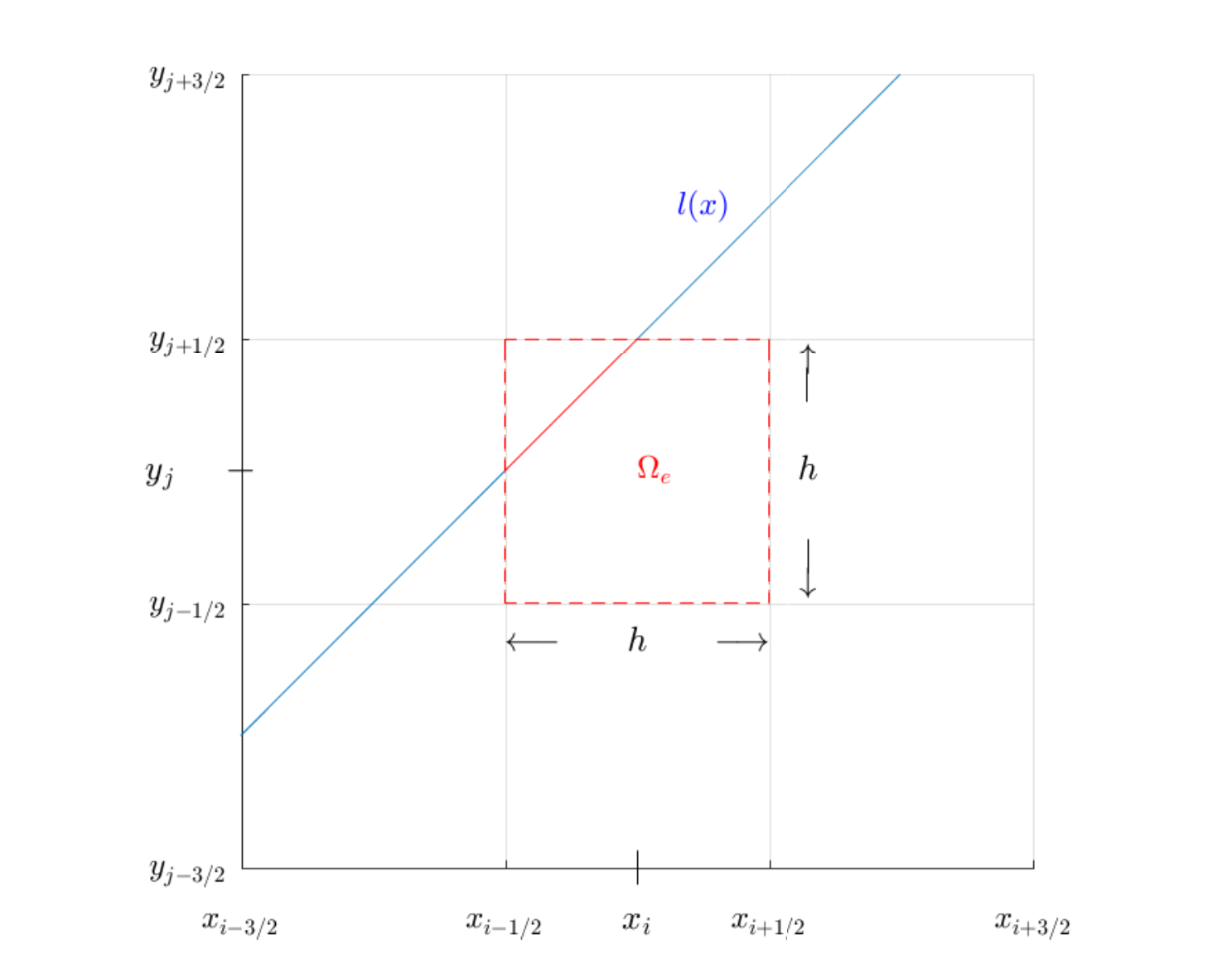}
        \caption{In this example the true interface is the line $l(x) = m \, x + b$
             Note that that the volumes $V_{i-1}$ and $V_{i}$ under the line in
            the first two columns $i-1$ and $i$ are exactly equal to the volumes due to
            the column sums $\tilde{V}_{i-1} = h^2 \,S_{i-1}$ and
            $\tilde{V}_{i} = h^2 \, S_{i}$ in the first and second columns of the
            $3 \times 3$ block of cells $B_{ij}$ centered on the center cell
            $\Omega_{e} (= \Omega_{ij})$.
            In this case the slope $\tilde{m} = S_{i} - S_{i-1}$ is \textit{exactly}
            equal to the slope $m$ of the interface as shown
            in~\eqref{EQ: m = S_i -_S_(i-1)}.
            It is always the case that if the true interface is a line, then one of the
            four standard rotations of $B_{ij}$ by a multiple of $90$ degrees about its
            center will orient the block so at least one of the divided differences
            of the column sums in~\eqref{EQ: THREE DIFFERENCES OF COLUMN SUMS IN x} or
            \eqref{EQ: THREE DIFFERENCES OF COLUMN SUMS IN y} is exact and hence, one of
            the linear approximations to the interface in the center cell $\Omega_{e}$
            defined in~\eqref{DEF: THE SIX LINES FROM THE COLUMN AND ROW SUMS}
            will always equal the interface in that cell, \textit{exactly},
            $\tilde{g}_{ij} (x) = m_{ij} \, x \, + \, b_{ij} = m \, x + b = l(x)$.
            In other words, the piecewise linear VOF approximation to $l(x)$ will
            \textit{always} reconstruct  the linear interface exactly.
            }
        \label{FIG: A LINE WITH TWO EXACT COLUMN SUMS}
    \end{figure}

\subsubsection{The ELVIRA Interface Reconstruction Algorithm}
\label{SUBSUBSECTION:THE ELVIRA INTERFACE RECONSTRUCTION ALGORITHM}

Here we briefly describe the ELVIRA interface reconstruction algorithm
\cite{JEP-EGP:2004} we used in this article.
In this example we present the simplest possible case; namely, when the true interface is
a line that passes through the center cell of the $3 \times 3$ block $B_{ij}$ of cells
$\Omega_{i'j'}$ centered on the cell $\Omega_{ij}$ as shown in
Figure~\ref{FIG: A LINE WITH TWO EXACT COLUMN SUMS}.
The following description is intended to be easy to understand.
However, the reader should be aware that there are many VOF interface reconstruction
algorithms in both two \cite{MDT-LDC-RCM-CWH:1985} and three dimensions
\cite{MDT-RCM-LRS:1987} and on every conceivable grid;
e.g., \cite{DBK-EGP-MWW-3:1999}, as well as numerous hybrid VOF/ Level Set algorithms
\cite{MMS-EGP:2000}.
See \cite{GT-RS-SZ:2007} and the references therein for a more complete overview VOF
methods.

In the ELVIRA algorithm the approximate interface will be a \textit{piecewise linear}
approximation $\tilde{g}_{ij} (x) = m_{ij} \, x \, + \, b_{ij}$ to the true interface in
$\Omega_{ij}$ as depicted in Figure~\ref{FIG: ELVIRA APPROXIMATION TO TANH}.
Furthermore \textit{the approximate interface is subject to the constraint that} the
volume fraction in the center cell due to the true interface  $g(x)$ and the approximate
interface $\tilde{g}_{ij}$ are equal; i.e., $f_{ij}^{true} = f_{ij}$.

Consider the example shown in Figure~\ref{FIG: A LINE WITH TWO EXACT COLUMN SUMS}.
In this example the true interface is a line $l(x) = m \, x + b$.
Assume we are given the exact volume fractions $f_{i'j'}$ associated with the
line $l(x)$, which is the true interface, in each cell $\Omega_{i'j'}$ of the
$3 \times 3$ block.
Then in this example the first two column sums
\begin{equation}
    S_{i-1} \equiv \sum_{j' =j-1}^{j+1} f_{{i-1},j'}
      \qquad \text{and} \qquad
    S_{i} \equiv \sum_{j' = j-1}^{j+1} f_{i,j'}
\end{equation}
are \textit{exact} in the sense that
\begin{equation}
  \label{EQ: AN EXACT COLUMN SUM}
    S_{i} \, = \,
     \frac{1}{h^2} \,
     \int\limits_{x_{i-1/2}}^{x_{i+1/2}} \, \left( l(x) - y_{j-3/2} \right) \, dx
\end{equation}
and similarly for $S_{i-1}$, but \textit{not} for $S_{i+1}$.
Thus, using~\eqref{EQ: AN EXACT COLUMN SUM} we find the difference in the column sums $S_{i}$ and
$S_{i-1}$ is
\begin{align}
  \label{EQ: DIFFERENCE IN COLUMN SUMS}
    h^2 \, \left( S_{i} \, - \, S_{i-1} \right)
     &= \, \int\limits_{x_{i-1/2}}^{x_{i+1/2}} \, \left(m \, x - b \right) \,
                                             - \,  y_{j-3/2} \; dx \,
      - \, \int\limits_{x_{i-3/2}}^{x_{i-1/2}} \, \left(m \, x - b \right)
                                             - \, y_{j-3/2} \; dx  \nonumber \\
      &= \, \int_{x_{i-1/2}}^{x_{i+1/2}} \, m \, x  \, dx
       - \, \int_{x_{i-3/2}}^{x_{i-1/2}} \, m \, x  \, dx  \nonumber \\
      &= \,  m \; \frac{x^2}{2} \; \bigg|_{x_{i-1/2}}^{x_{i+1/2}} \,
       - \,  m \; \frac{x^2}{2} \; \bigg|_{x_{i-3/2}}^{x_{i-1/2}}     \\
      & = \frac{m}{2} \left[ \left( {x_{i+1/2}} \right)^2
        - \left( {x_{i-3/2}} \right)^2  \right]
        - \, \frac{m}{2} \left[ \left( {x_{i+1/2}} \right)^2
                              - \left( {x_{i-1/2}} \right)^2 \right]  \nonumber \\
         &= \, \frac{m}{2} \, h \, \left( {x_{i+1/2}} - {x_{i-3/2}} \right)   \nonumber \\
         &= \, m \, \; h^2 \, .  \nonumber
\end{align}
Thus,
\begin{equation}
  \label{EQ: m = S_i -_S_(i-1)}
     m = S_i \, - \, S_{i-1}
\end{equation}
and we have recovered the \textit{exact} slope $m$ of the true interface $l(x)$ in the
center cell simply by differencing the correct pair of column sums of volume
fractions.
A little thought will show that the constraint
\begin{equation}
    f_{ij} = f^{true}_{ij}
\end{equation}
determines $b$ uniquely, thus determining the linear approximation
\begin{equation}
    g_{ij}(x) \, = \, m \, x \, + b
\end{equation}
which is exactly equal to the true interface $l(x)$.
In actual fact one needs to know whether the region containing the composition $C_1$ is
above, below, or to the left or right of $C_2$.
However, there are a variety of algorithms for doing this; e.g., see \cite{AJC:1985}
or~\cite{JEP-EGP:2004}.
This always works on a uniform grid of square cells, each of side $h$.

However, there are a few caveats:
There are three ways to difference the column sums,
    \begin{align}
      \label{EQ: THREE DIFFERENCES OF COLUMN SUMS IN x}
        m^{x,l} &= \left( S_{i} \, - \, S_{i-1} \right)   \nonumber \\
        m^{x,c} &= \frac{\left( S_{i+1} \, - \, S_{i-1} \right)}{2} \\
        m^{x,r} &= \left( S_{i+1} \, - \, S_{i} \right)   \nonumber
    \end{align}
    and three ways to difference the row sums
    \begin{align}
        \label{EQ: THREE DIFFERENCES OF COLUMN SUMS IN y}
        m^y_l &= \left( R_{j} \, - \, R_{j-1} \right)   \nonumber \\
        m^y_c &= \frac{\left( R_{j+1} \, - \, R_{j-1} \right)}{2} \\
        m^y_l &= \left( R_{j+1} \, - \, R_{j} \right)   \nonumber
      \end{align}
where the \textit{row sums} are defined by
\begin{equation}
  \label{DEF: THE ROW SUMS}
    R_{j-1} \equiv \sum_{i'=i-1}^{i+1} f_{{i'},j-1} \, ,
    \quad
    R_{j} \equiv \sum_{i' = i-1}^{i+1} f_{i',j}
    \quad \text{and} \quad
    R_{j+1} \equiv \sum_{i' =i-1}^{i+1} f_{{i'},j+1}
\end{equation}
In order to determine the best linear approximation to the true interface we compare the
the volume fractions $f^{x,l}_{i'j'}, \, f^{x,c}_{i'j'}, \, f^{x,r}_{i'j'}, \dots \,
f^{y,r}_{i'j'}$ due to each of the six lines
\begin{align}
  \label{DEF: THE SIX LINES FROM THE COLUMN AND ROW SUMS}
    g^x_l &= m^x_l \; x \, + \,  b^x_l \quad
   &g^y_l  = m^y_l \; x \, + \,  b^y_l \nonumber \\
    g^x_c &= m^x_c \; x \, + \,  b^x_c
   &g^y_c  = m^y_c \; x \, + \,  b^y_c           \\
    g^y_r &= m^x_r \; x \, + \,  b^x_r
   &g^y_r  = m^y_r \; x \, + \,  b^y_r \nonumber
\end{align}
we obtain from each of the six slopes in~\eqref{EQ: THREE DIFFERENCES OF COLUMN SUMS IN x}
and~\eqref{EQ: THREE DIFFERENCES OF COLUMN SUMS IN y} in the $3 \times 3$ block $B_{ij}$
centered on the cell of interest $\Omega_{ij}$ and use the line that minimizes the
difference between the given volume fractions and the volume fractions due to the lines
in~\eqref{DEF: THE SIX LINES FROM THE COLUMN AND ROW SUMS}.
We now explain this procedure in a bit more detail.

\subsubsection{Approximating an Unknown Interface}
\label{PARAGRAPH: APPROXIMATING AN UNKNOWN INTERFACE}

Suppose $g(x)$ is an unknown interface that passes through the center cell $\Omega_{ij}$
of a $3 \times 3$ block of cells $B_{ij}$ containing nine square cells $\Omega_{i'j'}$,
each of side $h$, centered on $\Omega_{ij}$.
Furthermore, assume the only information we have are the nine \textit{exact} volume
fractions $f_{i'j'}$ in the cells $\Omega_{i'j'}$ due to $g(x)$.
For example, in Figure~\ref{FIG: ELVIRA APPROXIMATION TO TANH} the `unknown' interface is
$g(x) = \tanh(x)$, which is the blue curve, and the volume fractions are nonzero only in
cells that either contain the curve or are below it.
We want to find a line segment $\tilde{g}_{ij} (x) = m_{ij} \, x + b_{ij}$ that is a
second-order accurate approximation to $g(x)$, in the following sense,
\begin{equation}
  \label{DEF: SECOND ORDER IN THE MAX NORM}
   \max \left| g(x) - \tilde{g}_{ij} (x) \right|
     \le \, \tilde{C} \, h^2
     \qquad \text{for all } x \in [x_{i-1/2}, \, x_{i+1/2}] \, ,
\end{equation}
 where $\tilde{C}$ is a constant that is independent of $h$.

First we define a way to measure the error $E(\tilde{m})$ between the volume fractions
$f_{i'j'}$  we are given that are due to the unknown interface and the approximate volume
fractions $\tilde{f}_{i'j'}$ due to a line segment
$\tilde{g} (x) = \tilde{m} \, x + \tilde{b}$ that passes through the center cell
$\Omega_{ij}$ and the $3 \times 3$ block $B_{ij}$ centered on $\Omega_{ij}$,
\begin{equation}
  \label{DEF: TWO NORM SQUARED}
    E(\tilde{m}) \,
      = \, \sum\limits_{i' = i-1}^{i+1}
        \, \sum\limits_{j' = j-1}^{j+1} \left(f_{i'j'} - \tilde{f}_{i'j'} \right)^2 \, .
\end{equation}
Note that this is the square of the \textit{two norm} on vector spaces $R^n$ from linear
algebra, where in our case $n =9$, \cite{WGS:2016}.

Now take the volume fractions we are given, namely $f_{i'j'}$, and form all six of the
slopes in~\eqref{EQ: THREE DIFFERENCES OF COLUMN SUMS IN x}
and~\eqref{EQ: THREE DIFFERENCES OF COLUMN SUMS IN y} and the six candidate lines
in~\eqref{DEF: THE SIX LINES FROM THE COLUMN AND ROW SUMS} from these slopes.
Remember that the `$y$ intercept' $b$ for each of the lines
in~\eqref{DEF: THE SIX LINES FROM THE COLUMN AND ROW SUMS} is determined by the constraint
$f_{ij}^{true} = f_{ij}$.
Each of the six line produces nine volume fractions in the $3 \times 3$ block $B_{ij}$.
For example, given the slope $m^{x,c}$ defined
in~\eqref{EQ: THREE DIFFERENCES OF COLUMN SUMS IN x} we obtain the line
$g^x_c = m^x_c \; x \, + \,  b^x_c$ defined
in~\eqref{DEF: THE SIX LINES FROM THE COLUMN AND ROW SUMS}, which in turn gives us nine
volume fractions $f^{x,c}_{i'j'}$ for $i' = i-1, \, i, \, i+1$ and
$j' = j-1, \, j, \, j+1$.
Now compute $E(m^x_c)$ and repeat this procedure for each of the other lines
in~\eqref{DEF: THE SIX LINES FROM THE COLUMN AND ROW SUMS} with slopes computed as in
in~\eqref{EQ: THREE DIFFERENCES OF COLUMN SUMS IN x}
and~\eqref{EQ: THREE DIFFERENCES OF COLUMN SUMS IN y}.
Finally, take the line from~\eqref{DEF: THE SIX LINES FROM THE COLUMN AND ROW SUMS}
that minimizes the error defined in~\eqref{DEF: TWO NORM SQUARED}; i.e., pick the slope
from~\eqref{EQ: THREE DIFFERENCES OF COLUMN SUMS IN x}
and~\eqref{EQ: THREE DIFFERENCES OF COLUMN SUMS IN y}, call it $\tilde{m}$, that satisfies
\begin{equation}
  \label{EQ: MINIMIZE THE ERROR}
    E(\tilde{m}) \, = \,
      \min \left\{ E(m^x_l), \, E(m^x_c), \ldots, \, E(m^y_r) \right\} \, .
\end{equation}
The line
\begin{equation}
  \label{EQ: ELVIRA LINEAR APPROXIMATION}
    \tilde{g} \, = \, \tilde{m} \; x \, + \, \tilde{b}
\end{equation}
is the linear approximation to the true interface $g(x)$ in $\Omega_{ij}$ that we use in
the VOF algorithm in this article.
In~\cite{EGP:2010a} and \cite{EGP:2014} it is proven that this algorithm produces a
second-order accurate approximation in the sense
of~\eqref{DEF: SECOND ORDER IN THE MAX NORM} to the interface provided that
\begin{equation}
  \label{DEF:VOF EGP Constraint}
    h \, \le \, \frac{2}{33 \, \sigma_{max}}
\end{equation}
where  $\sigma_{max}$ denotes the maximum curvature of the interface, $h$ is the grid
size, and the volume fractions due to the true interface are exact.

\subsubsection{Implementation of the VOF method in ASPECT}
\label{SUBSUBSECTION:Implementation of the VOF method in ASPECT}

In the work described in this article we implemented the VOF algorithm described above in
ASPECT. 
We will now describe our implementation for \textit{square}, two dimensional, cells $\Omega_e$ in \textit{physical space} (often referred to as the `real' cell in ASPECT) and show computational results for such cells.
First note that the VOF method is essentially a specialized version of a Finite Volume Method, which is equivalent to a Discontinuous Galerkin (DG) method with values $f_e$ that are constant on each cell $\Omega_e$.
Approaching the VOF algorithm from this point of view, we note that both methods require the computation of the flux of the volume of $C_1$, or, equivalently, the volume fraction of $C_1$, across each of the edges of $\Omega_e$.

In a VOF method it is natural to use the method of characteristics to calculate the flux of $C_1$ through each of the cell edges.
This is done by tracing backward in time along a linear approximation to each  characteristic in order to identify the \textit{total} volume that will cross a given edge and then computing that portion of the volume associated with the fluid that is being tracked as shown in
Figure~\ref{FIG:COMPUTATION OF VOLUME FLUX ACROSS THE RH EDGE}; i.e., by computing the volume fraction of $C_1$ in the total volume that will cross that edge. 
See~\cite{PC:1990} and~\cite{RLV:1996} for examples of computing a second-order accurate  flux in this manner in a finite volume discretization
of~\eqref{EQ:CONSERVATION EQUATION FOR THE VOLUME FRACTION f}, rather than a VOF discretization of~\eqref{EQ:CONSERVATION EQUATION FOR THE VOLUME FRACTION f}, as well as higher resolution versions of these algorithms.
In our computation of the volume fraction flux we make use of several algorithms that we developed for the interface reconstruction step.
We will describe these algorithms is more detail below.

There are a number of approaches one can consider for obtaining the velocities on the $k$th edge from the approximate FEM solution of the incompressible Stokes equations.
Two approaches are (1) a point sample of the normal velocity on the $k$th edge and (2) the velocity integrated 
\begin{equation}
  \label{EQ:INTEGRATED VELOCITY ON THE KTH EDGE}
    \int \limits_{\partial \Omega_{e,k}} \mathbf{u} \cdot \mathbf{n}_k \; ds
\end{equation}
along the $k$th edge of $\Omega_e$, where $\mathbf{n}_k$ denotes the unit normal to the  $k$th edge of $\Omega_e$.
For a finite volume method, (1) and (2) are both reasonable approximations to the edge velocities.
However, the latter method~\eqref{EQ:INTEGRATED VELOCITY ON THE KTH EDGE} is a closer analogue to the type of procedure one would choose for a finite element method.

\begin{figure}
    \begin{center}
        \begin{tikzpicture}
            \draw (0,1) rectangle (1,2);
            \draw[fill=blue, opacity=0.5] (0,1) -- (0,1.25) -- (1,1.75) -- (1,1) -- cycle;
            \draw[blue] (0,1.25) -- (1,1.75);
            \draw[pink] (0.75,1) rectangle (1,2);
            \draw[green] (0.75,1) -- (0.75,1.625) -- (1,1.75) -- (1,1) -- cycle;
            \draw[red] (0.75,1.625) -- (1,1.75);
            \draw[->] (1,1.5) -- (1.25,1.5) node[above] {$u$};
            \draw[<->] (0,0.9) -- node[below] {$h$} (1,0.9);
            \draw[<->] (-0.1,1) -- node[left] {$h$} (-0.1,2);
            \draw (0.5,0) node {$\Omega_e$};
            \draw[densely dotted] (1.1,1) -- (1.9,0);
            \draw[densely dotted] (1.1,2) -- (1.9,4);
            \draw (2,0) rectangle (6,4);
            \draw[->] (6,2) -- (6.25,2) node[above] {$\mathbf{n}_{k}$};
            \draw[blue] (2,1) -- (6,3);
            \fill[opacity=0.5] (4,2) circle[radius=0.05];
            \draw (4,2) node[below] {$x^c_{\hat{e}}$};
            \draw (3,2) node {$d_{\hat{e}}=0$};
            \draw[->] (4,2) -- (3.875,2.5) node[above] {$\mathbf{n}_{\hat{e}}$};
            \draw[pink] (5,0) rectangle (6,4);
            \draw[green] (5,0) -- (5,2.5) -- (6,3) -- (6,0) -- cycle;
            \draw[red] (5,2.5) -- (6,3);
            \draw[<->] (2,-0.1) -- node[below] {$1$} (6,-0.1);
            \draw[<->] (1.9,0) -- node[left] {$1$} (1.9,4);
            \draw (4,-1) node {$\hat{\Omega}_e$};
            \draw[black,densely dotted] (6.1,4) -- (6.9,4);
            \draw[black,densely dotted] (6.1,0) -- (6.9,0);
            \draw[pink] (7,0) rectangle (11,4);
            \draw[green] (7,0) -- (7,2.5) -- (11,3) -- (11,0) -- cycle;
            \draw[red] (7,2.5) -- (11,3);
            \fill[opacity=0.5] (9,2) circle[radius=0.05];
            \draw (9,2) node[below] {$x^c_I$};
            \draw[<->] (9,2) -- node[right] {$d_I$} (8.908,2.738);
            \draw[->] (8.908,2.738) -- (8.8712,3.0332) node[above] {$\mathbf{n}_I$};
            \draw[<->] (7,-0.1) -- node[below] {$1$} (11,-0.1);
            \draw[<->] (6.9,0) -- node[left] {$1$} (6.9,4);
            \draw (9,-1) node {$\hat{\Omega}_I$};
        \end{tikzpicture}
        \caption{
            A diagram of the mapping of the region (in purple) containing the compositional field $C_1$ in the  \textit{real} cell $\Omega_e$ to its associated unit cell $\hat{\Omega}_e$.
            In this diagram we have assumed a constant velocity field $u$ lying only in the $x$ direction in the unit cell $\hat{\Omega}_e$ so that the flux of $C_1$ across the RH edge of $\hat{\Omega}_e$ is a rectangular region.
            This allows us to compute the \textit{total} volume $V_F$ of $C_1$ and $C_2$ that is fluxed across the RH edge of $\hat{\Omega}_e$; namely the rectangle on the RH edge of $\hat{\Omega}_e$.
            We then map this rectangle to another unit cell $\hat{\Omega}_I$ in order to compute the volume \textit{fraction} $f_k$ of the (mapped) rectangle that contains the composition $C_1$. 
            Since in this article, linear interfaces map to linear interfaces, we can use the unit normal $\mathbf{n}_I$ and distance $d_I$ to calculate the volume fraction $f_k$ of $C_1$ in this rectangle.
            Note that, since this diagram has been chosen to correspond exactly to the one in Figure~\ref{FIG:COMPUTATION OF VOLUME FLUX ACROSS THE RH EDGE}, $d_{\hat{e}} = 0$. 
            However, in general, $d_{\hat{e}} \ne 0$.
        }
        \label{FIG:SKETCH OF HOW TO COMPUTE V_F}
    \end{center}
\end{figure}
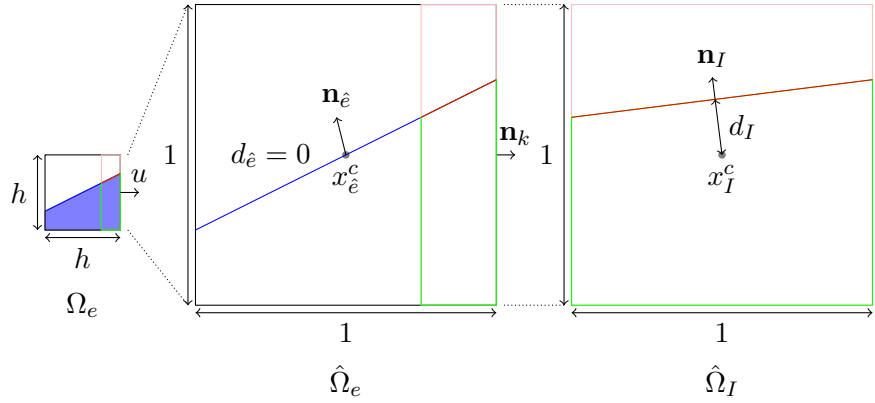

We now describe our implementation of the computation of the volume flux of $C_1$, on a square cell $\Omega_e$ of side $h$.
(When we employ AMR, $h$ denotes the length of each side of the most finely resolved cells in the FEM grid.)
All of the information that we use to describe the interface in a real cell $\Omega_e$; namely, its distance $d_e$ to the center of the cell and the unit normal $\mathbf{n}_e$ to the interface, is stored with respect to the interface's location relative to the center of the \textit{unit} cell $\hat{\Omega}_e$ associated with $\Omega_e$ when $\Omega_e$ is mapped to $\hat{\Omega}_e$ as depicted in
Figure~\ref{FIG:SKETCH OF HOW TO COMPUTE V_F}.
In particular, the interface in the unit cell $\hat{\Omega}_e$  is given by
\begin{equation}
  \label{EQ:EQUATION DESCRIBING THE INTERFACE IN THE UNIT CELL}
	\mathbf{n}_{\hat{e}} \cdot (\mathbf{x} \, - \, \mathbf{x}^{c}_{\hat{e}})
      \, = \, d_{\hat{e}}
\end{equation}
where $\mathbf{x}^{c}_{\hat{e}}$ is the center of  $\hat{\Omega}_e$, $d_{\hat{e}}$ is the distance of the interface from the center $\mathbf{\mathbf{x}}^{c}_{\hat{e}}$ of $\hat{\Omega}_e$, and $\mathbf{n}_{\hat{e}}$ is a unit vector that is perpendicular to the approximate (linear) interface in $\hat{\Omega}_e$, with the convention that $\mathbf{n}_{\hat{e}}$ always points away from the region containing Composition~1.
The location of the interface is stored by recording $\mathbf{n}_{\hat{e}}$ and $d_{\hat{e}}$ for each cell $\Omega_e$ that contains a portion of the interface.
For the case when the velocity field is perpendicular to a cell edge, say 
$\partial \Omega_{e,k}$, for some $k  \, =  \, 1, \, 2, \,  3,  \, 4$, let $\hat{\mathbf{n}}_k$ be the unit normal vector to the $k$th edge
$\partial \hat{\Omega}_{e,k}$ of the unit cell $\hat{\Omega}_e$,  and let $V_F$ denote the \textit{total} volume flux (i.e., the volume of $C_1$ \textit{and} $C_2$) that will flux / advect across $\partial \Omega_{e,k}$.

As shown in Figure~\ref{FIG:SKETCH OF HOW TO COMPUTE V_F}, with only a few computationally inexpensive transformations we can use the same algorithm we used to compute the volume fraction on a cell $\Omega_e$ in the reconstruction step to compute the volume flux of $C_1$ across each of the edges of $\Omega_e$.
If we map $V_F$ from $\hat{\Omega}_e$ to another unit cell $\hat{\Omega}_I$ and assuming the velocity is perpendicular to the $k$th cell edge $\partial \hat{\Omega}_{e,k}$ of $\Omega_e$, we find that the interface within the unit cell $\hat{\Omega}_I$ is given by
\begin{equation*}
     \mathbf{n}_I \cdot(\mathbf{x} \, - \, \mathbf{x}^{c}_I) \, = \, d_I
\end{equation*}
where  $\mathbf{x}^{c}_I$ is the center of $\hat{\Omega}_I$ as shown in Figure~\ref{FIG:SKETCH OF HOW TO COMPUTE V_F}.
The values of $\mathbf{n}_I$ and $d_I$ in terms of $\mathbf{n}_e$, $\mathbf{n}_{k}$, and $d_e$ are given by
\begin{align}
    \mathbf{n}_I &= \mathbf{n}_{\hat{e}} +(\frac{V_F}{V_e}   
    -\mathbf{n}_{\hat{e}} \cdot \mathbf{n}_{k} )\mathbf{n}_{k}           \, , \\
    d_I         &= d_{\hat{e}}-(\frac{1}{2}+\frac{V_F}{2V_e})(\mathbf{n}_{\hat{e}}\cdot\mathbf{n}_{k}) \, .
\end{align}
where $V_e$ is the volume of $\Omega_e$ (the upwind cell for this edge), and $\mathbf{n}_{k}$ is the outward pointing normal of the cell edge $\partial \hat{\Omega}_{e,k}$.

Given the assumptions we have made regarding a uniform square grid, we have a constant Jacobian, so the volumes on the unit cell and the volume in physical space are related by
a constant multiple.
For a given interface, there is a simple formula to calculate the volume of $C_1$ on the side opposite the unit normal $\mathbf{n}$; e.g., see \cite{RS-SZ:2000}.
In our notation this formula is
\begin{equation}
    \label{EQ:VOLUME FRACTION OF C_1}
    f(\mathbf{n},d) = \begin{cases}
        1 & \frac{1}{2}\leq \bar{d}\\
        1-\frac{(\bar{d}-\frac{1}{2})^2}{2m(1-m)} & \frac{1}{2}-m<\bar{d}<\frac{1}{2}\\
        \frac{1}{2}+\frac{\bar{d}}{(1-m)} & m-\frac{1}{2}\leq \bar{d}\leq \frac{1}{2}-m \\
        \frac{(\bar{d}+\frac{1}{2})^2}{2m(1-m)} & -\frac{1}{2} < \bar{d} <m-\frac{1}{2} \\
        0 & \bar{d} \leq -\frac{1}{2}\\
    \end{cases}
\end{equation}
Where $m = 1 - \frac{\|\mathbf{n}\|_\infty}{\| \mathbf{n} \|_{1 \,}}$ and the components of $\mathbf{n}$ are parallel to sides of the unit cell $\hat{\Omega}_e$, and $\bar{d}=\frac{d}{\|\mathbf{n}\|_{1 \,}}$.
We use~\eqref{EQ:VOLUME FRACTION OF C_1} to compute the flux of $C_1$ across the RH edge, which is $f(\mathbf{n}_I,d_I) \, V_F$.
In general we use an analogous procedure to compute the (volume) flux of $C_1$ across the other edges of $\Omega_e$.

\subsubsection{Volume Correction}
\label{SECTION:Volume Correction}


In our current implementation of the VOF advection algorithm we use a dimensionally split
algorithm as described in Section~\ref{SUBSUBSECTION:DESCRIPTION} above.
Consequently, we cannot assume that the velocity at the intermediate step is divergence free due to the decoupling of the cell edges from one spatial dimension to the other.
This decoupling removes the guarantee that the volume fractions retain the bound
$0 \leq f \leq 1$, since there is no guarantee that the velocity is divergence free during the intermediate step.
Since the reconstruction algorithm requires $0 \leq f \leq 1$, it is therefore necessary to modify the advection algorithm.

First, note that the equation that governs the advection of the characteristic function $f$ is
\begin{equation}
  \label{EQ:ADVECTION EQUATION FOR THE CHARACTERISTIC FUNCTION F}
    \frac{\partial}{\partial t} f + \mathbf{u} \cdot \nabla f = 0 \, .
\end{equation}
Using~\eqref{EQ:NONDIMENSIONAL VECTOR FORM OF THE CONTINUITY EQUATION} we obtain a modified form of~\eqref{EQ:ADVECTION EQUATION FOR THE CHARACTERISTIC FUNCTION F},
\begin{equation}
  \label{EQ:MODIFIED ADVECTION EQUATION FOR THE CHARACTERISTIC FUNCTION F}
    \overbrace{\frac{\partial}{\partial t} f + \nabla\cdot(\mathbf{u}f)}^{\text{Advection}} -
    \overbrace{f(\nabla\cdot\mathbf{u})}^{\textrm{Correction}} = 0
\end{equation}

Note that in~\eqref{EQ:MODIFIED ADVECTION EQUATION FOR THE CHARACTERISTIC FUNCTION F} the first term is a conservation law for $f$.
If~\eqref{EQ:NONDIMENSIONAL VECTOR FORM OF THE CONTINUITY EQUATION} is satisfied \textit{exactly} then the correction term 
in~\eqref{EQ:MODIFIED ADVECTION EQUATION FOR THE CHARACTERISTIC FUNCTION F} will be zero. However, in the case of a dimensionally split algorithm the assumption that the velocity $\mathbf{u}$ is divergence free, even to $O(h^q)$ for some integer $q \ge 2$, breaks down.
One expects this to add a small error to the computation.
However, since the VOF interface reconstruction algorithm requires $f$ to satisfy
$0 \leq f \leq 1$ it is necessary to retain this correction term.

We approximate~\eqref{EQ:MODIFIED ADVECTION EQUATION FOR THE CHARACTERISTIC FUNCTION F} by
\begin{equation}
  \label{EQ:DISCRETIZED ADVECTION EQUATION FOR THE CHARACTERISTIC FUNCTION F}
    f^{n+1}_e \, V_e \, = \, f^n_e \, V_e \, + \, \sum\limits_k \, f_k \, U_k 
       \, - \, \bar{f}_e \, \sum\limits_k U_k \, ,
\end{equation}
where $e$ is an index that ranges over all cells $\Omega_e$, $V_e$ is the volume of $\Omega_e$, $k$ is an index that ranges over the cell edges of $\Omega_e$, $f_k$ is the volume \textit{fraction} of $C_1$ that will be fluxed across the $k$th edge as described  in the caption to Figure~\ref{FIG:SKETCH OF HOW TO COMPUTE V_F}, and 
\begin{equation}
  \label{EQ:DEFINITION OF U_k}
    U_k \, = \, \Delta t \, 
      \int\limits_{\partial\Omega_{e,k}} \, \tilde{\mathbf{u}}_k \cdot \mathbf{n}_k \, ds \, .
\end{equation}
where $ \tilde{\mathbf{u}}_k$ is a time centered approximation to the velocity  $\mathbf{u}$ on the $k$th edge, 
\begin{equation*}
    \tilde{\mathbf{u}}_k \, = \, \frac{\mathbf{u}^{n+1} + \mathbf{u}^n}{2} \, .
\end{equation*}
The term $\bar{f}_e$
in~\eqref{EQ:DISCRETIZED ADVECTION EQUATION FOR THE CHARACTERISTIC FUNCTION F} can be one of several approximations to the volume fraction $f_e$ in $\Omega_e$.
The two simplest cases are (1) $\bar{f}_e=f^{old}_e$, i.e., an explicit correction term, and (2) $\bar{f}_e=f^{new}_e$, i.e., an implicit correction term.
In the results shown here we use the explicit term $\bar{f}_e \, = \, f^{old}_e$.

Since we are using Strang
splitting,~\eqref{EQ:DISCRETIZED ADVECTION EQUATION FOR THE CHARACTERISTIC FUNCTION F} is evaluated once for each spatial dimension in the problem at each time step, alternating the order of the dimensions in the subsequent time step.

\subsubsection{Model Coupling Procedure}
\label{SUBSUBSECTION:MODEL COUPLING PROCEDURE}

Having now described our implementation the VOF method in ASPECT, it is necessary to establish how the computed $C$ field may be used by the other portions of the model in cases where the tracked fluid is not a passive tracer.

There is a significant reduction in the complexity of the implementation and duplication of work if the $C$ field can be discretized in the same manner as what are known as  ``compositional fields'' in ASPECT~\cite{ASPECT-MANUAL:2017}.
Furthermore, in order to avoid interfering with the values on neighboring cells, we prefer to use a discontinuous element.
For a number of reasons, often relating to the physical interpretation of the quantity $C$, it is also desirable to ensure that it will always be bounded; 
e.g., $0 \leq C \leq 1$.

A basic implementation can be done by directly discretizing the volume fraction data on a discontinuous $P_0$ element, which is equivalent to a least squares approximation to the composition field implied by the reconstructed interface.
Attempting to obtain an ideal approximation using a higher order element such as DG $Q_1$ or DG $P_1$ is more difficult, especially if the bounds on the composition $C$ are respected, since the result of a unconstrained least squares approximation for such an element is almost certain to violate these bounds.
Furthermore, the unconstrained least squares computation can be expected to be both more complex and more expensive.
Thus, any approximation using a non-constant DG element will require a heuristic approach.

In our implementation, in order to generate a DG $Q_1$ element approximation to the $C$ field that is implied by the reconstructed interface, we apply the following constraints.
\begin{enumerate}
    
    \vskip 06pt
        
    \item The gradient of the element is in the same direction as the normal of the 
        interface.
        
    \vskip 06pt
    
    \item The gradient is as large as possible while maintaining $0 \le C \le 1$ 
        everywhere.

    \vskip 06pt
        
    \item In order to maintain conservation of mass the volume fraction implied by the DG 
        $Q_1$ element approximation to the $C$ field must match the volume fraction $f_e$ in the VOF approximation to the $C$ field; i.e.,
        \begin{equation*}
           \int\limits_{\Omega_e} \; C (\mathbf{x}) \; d\mathbf{x} \, = \, f_e \, V_e
        \end{equation*}
          where $V_e$ is the volume of $\Omega_e$.
    \vskip 06pt
    
\end{enumerate}
On a square mesh, for a cell with the reconstructed interface
\begin{equation}
    \mathbf{n}_{\hat{e}} \cdot (\mathbf{x} \, - \, \mathbf{x}_{\hat{e}}^{c} )
       \, = \, d_{\hat{e}}
\end{equation}
the above constraints result in the approximation on the unit cell being
\begin{equation}
    C(\mathbf{x}) \, = \, f_{\hat{e}} \, - \,
      1 - | \, 2 \, f_{\hat{e}} \, - \, 0.5 | \; \;
      \frac{\phantom{\|}\mathbf{n}_{\hat{e}} \phantom{\|_1}}{\|\mathbf{n}_{\hat{e}}\|_1} \cdot (\mathbf{x} \, - \,  \mathbf{x}_{\hat{e}}^{c} )
\end{equation}

If we use an DG $Q_1$ element, the use of the above equation produces a bilinear approximation to the VOF method's reconstructed $C$ field, with little additional computational cost over the $P_0$ approximation.

\subsubsection{Coupling with the AMR Algorithm}
\label{SUBSUBSECTION:COUPLING WITH THE AMR Algorithm AMR Algorithm}

The deal.ii library~\cite{dealII85} upon which ASPECT is built manages the AMR algorithm through the p4est library~\cite{CB-LCW-OG:2011}. 
Deal.ii, and hence, ASPECT provides a mechanism for setting the refinement criteria; both when to refine a cell and when to coarsen a cell.
Since reconstructing and advecting the interface across different levels of refinement both increases algorithm complexity and decreases the accuracy with which the interface is resolved, in this work we ensure that the interface is always on the finest level of refinement.
This approach requires that the cells that contain the interface, including the case where the interface is on a cell boundary, and any cell that shares a vertex with any of those cells must also be at the finest level of refinement.

The criteria for refining a cell that we have adopted is a two step algorithm that requires one pass over the entire mesh and one pass over a subset of the entire mesh.
In the first step we check every cell in the entire mesh making a list of all cells that
contain a part of the interface.
More specifically, we regard all cells $\Omega_e$ that satisfy
$\epsilon_{vof} < f_e < 1 - \epsilon_{vof}$, where  $\epsilon_{vof}$ is a small parameter,
to contain a portion of the interface.
In addition, all cells $\Omega_e$ that have a neighboring cell $\Omega_e^{\prime}$ that shares a face with $\Omega_e$ and differ in volume fraction sufficiently (e.g., $| f_e - f_e^{\prime} | > \epsilon_{vof}$) are also added to this list.
In the computational results shown in~Section~\ref{SECTION:NUMERICAL RESULTS} we use the value $\epsilon_{vof} = 10^{-6}$.
In the second pass over a subset of the entire grid we make a list of all cells that share a vertex with any cell already in the list of cells that contain a portion of the interface and also flag each of these cells for refinement.
These flags are then passed to deal.ii, and thus on to p4est, which handles the details of the refinement of these cells and the coarsening of those cells that no longer need to be at the finest level of refinement. 

Given that the time step  $\Delta t$ is constrained by a CFL condition, the interface can move at most $\sigma$ cell widths where $\sigma \le 1$ is the CFL number.
This permits the reduction of the frequency with which we conduct the remeshing procedure to $N$ time steps where $N < \frac{W - 2}{2 \sigma}$ and $W$ is the minimum width of the maximally refined band of cells.
(See, for example, any of the AMR computations in the second (b) and fourth (d) frames in
Figures~\ref{FIG:B = 0.0}--\ref{FIG:B = 2.0} for explicit examples of $W$.)
For the refinement strategy described above, the safest assumption is that $W = 4$.
This takes into consideration the case where the interface is at the cell boundary.
A band of larger width $W > 4$ would both require a more complex algorithm to find the
necessary cells to flag and would increase the number of refined cells. 
Thus, there is a balance between cost associated with the frequency of running the algorithm to flag cells for refinement and cost of having a larger value of $W$.
This balance is problem dependent.


\section{Numerical Results}
\label{SECTION:NUMERICAL RESULTS}

In this section we present our numerical results.
First, in Section~\ref{SUBSECTION:INTERFACE TRACKING BENCHMARK PROBLEMS} we compute two test problems
with prescribed velocity fields to verify the accuracy of our implementation of
the VOF algorithm \cite{JEP-EGP:2004}.
Then, in Section~\ref{SUBSECTION:MANTLE CONVECTION BENCHMARK PROBLEMS}, we compute some mantle
convection benchmarks to verify the accuracy of the coupling to the mantle
convection code.
Finally, in Section~\ref{SUBSECTION:THERMOCHEMICAL CONVECTION IN A DENSITY STRATIFIED FLUID}, we apply the algorithm to a problem of interest in the field of geodynamics.

\subsection{Interface Tracking Benchmark Problems}
\label{SUBSECTION:INTERFACE TRACKING BENCHMARK PROBLEMS}

In this section, we compute two test problems with known exact solutions to ensure that our the implementation of the VOF algorithm converges at its design rate.
In particular, in Subsection~\ref{SUBSUBSECTION:ADVECTION OF A LINEAR INTERFACE} we compute one of the test problems from~\cite{JEP-EGP:2004}, and in Subsection~\ref{SUBSUBSECTION:ROTATION OF A CIRCULAR INTERFACE} we compute an modified version of a test problem from the same paper.

\subsubsection{Definition of the Error Measurement}
\label{SUBSUBSECTION:DEFINITION OF THE ERROR MEASUREMENT}

Since each volume fraction $f_e$ is constant on its grid cell $\Omega_e$, we use $P_0$ elements to store the value of the volume fraction $f_e$ on each $\Omega_e$.
Given a \textit{fixed} grid with cells $\Omega_e$ indexed by $e$ we define the error between the exact $f_e^{exact}$ and computed $f_e^{comp}$ volume fractions by
\begin{equation}
  \label{DEF:VOLUME FRACTION ERROR}
    \text{Error} \left(f^{exact} \, - \, f^{comp} \right) \, = \,         
      \sum\limits_{e} \, \left| \, f_e^{exact} \, - \, f_e^{comp} \, \right|
        V( \Omega_e )
\end{equation}
where $V( \Omega_e )$ denotes the volume of $\Omega_e$.
Note that this is the $L^1$ norm of the difference between $f_e^{exact}$ and$f_e^{comp}$ with weight $V( \Omega_e )$.

\subsubsection{Advection of a Linear Interface in a Constant Velocity Field}
\label{SUBSUBSECTION:ADVECTION OF A LINEAR INTERFACE}

\begin{figure}
    \centering
    \begin{tikzpicture}[scale=0.5]
        \draw[<->] (0,4.1) -- node[above] {1} (4,4.1);
        \draw (0,0) rectangle (4,4);
        \fill[color=blue, opacity=0.7] (4,0) -- (4,4) -- (0,4);
        \draw[<-] (1.75,1.75) -- (1.15,1.0) node[right] {\(\mathbf{u}\)};
    \end{tikzpicture}

    \caption{Initial conditions for the ``Advection of a Linear Interface in a Constant Velocity Field'' test problem}
    \label{FIG:ADVECTION OF A LINEAR INTERFACE INITIAL CONDITIONS}
\end{figure}
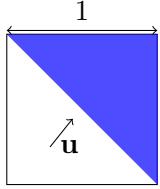

\begin{figure}
    \centering
    \includegraphics[width=0.49\linewidth]{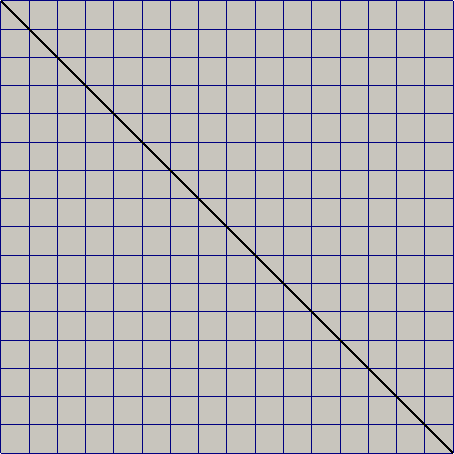}
    \includegraphics[width=0.49\linewidth]{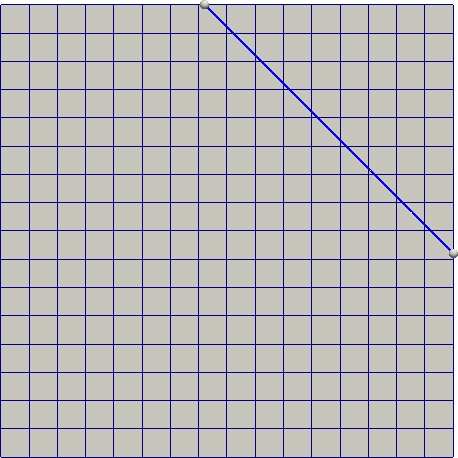}
    \caption{On the left is the initial condition as reconstructed by our VOF method.
        On the right is a comparison between the exact and computed interface at $t = 1$, with the exact interface in blue and the computed interface in black.
        It is apparent that the two interfaces are visually indistinguishable.}
    \label{FIG:LINEAR INTERFACE TRANSLATION}
\end{figure}

Our first VOF benchmark problem is the advection of a linear interface in a constant
velocity field as shown in
Figures~\ref{FIG:ADVECTION OF A LINEAR INTERFACE INITIAL CONDITIONS}
and~\ref{FIG:LINEAR INTERFACE TRANSLATION}.
The computational domain is $[0,1] \times [0,1]$ and the initial interface given by
$y = 1 - x$.
At each time step $t^k \to t^{k+1}$ the interface is advanced the velocity field
$\mathbf{u}=(-\frac{25}{100},-\frac{24}{100})$, and then compared with the exact solution, for which the interface is given by $y = \frac{51}{100} - x$.
In this computation we used a CFL number of $\sigma = \frac{1}{2}$, which resulted in, for example, of a total of $23$ time steps on the least refined grid of $h = 2^{-16}$ .

Since the ELVIRA interface reconstruction method reproduces lines exactly, we expect the error in the computations to be exact to machine precision
$\epsilon_{\text{mach}} \approx 10^{-16}$.
The results of computations for $h = ...$ are given in
Table~\ref{TABLE:ADVECTION OF A LINEAR INTERFACE}.
We note that in all cases the error is $O \left( \epsilon_{\text{mach}} \right)$.

\begin{table}
    \centering
    \begin{filecontents*}{translLin2.csv}
size,h,nproc,errorI,cRateI,errorF,cRateF,time(sec)
4,$2^{-4}$,4,1.226292e-16, nan,1.233822e-16, nan,1.310000e+00
5,$2^{-5}$,4,1.215846e-16,0.01,1.216746e-16,0.02,4.140000e+00
6,$2^{-6}$,4,2.959404e-16,-1.28,2.960829e-16,-1.28,2.680000e+01
7,$2^{-7}$,4,5.927208e-16,-1.00,5.927376e-16,-1.00,1.160000e+02
\end{filecontents*}
\pgfplotstabletypeset[
    search path={DATA_TABLES/VOF_Benchmarks},
    col sep=comma,
    columns={h, errorF},
    columns/h/.style={
        string type,
        column name={\(h\)},
        column type/.add={|}{|},
    },
    columns/errorF/.style={
        column name={Error},
        sci zerofill,
        dec sep align={c|},
        column type/.add={}{|},
        precision=5,
    },
    after row=\hline,
    every head row/.style={
        before row=\hline,
    },
    ]{translLin2.csv}
    \caption{Error when advection a linear interface using a constant velocity
        field not aligned to the mesh. Note that machine precision is
        $\epsilon_{\text{mach}}\approx 10^{-16}$ and the number of cells that the
    interface passes through is approximately $\frac{L}{h}$.}
    \label{TABLE:ADVECTION OF A LINEAR INTERFACE}
\end{table}

\subsubsection{Rotation of a Circular Interface}
\label{SUBSUBSECTION:ROTATION OF A CIRCULAR INTERFACE}

\begin{figure}
    \centering
    \begin{tikzpicture}[scale=0.5]
        \draw[<->] (0,4.2) -- node[above] {1} (4,4.2);
        \draw (0,0) rectangle (4,4);
        \fill[color=blue, opacity=0.7] (2.5,2) circle (0.5);
        \fill[color=red] (2,2) circle (0.05);
        \draw[->] (2,3) arc (90:270:1) node[below] {\(\mathbf{u}\)};
    \end{tikzpicture}
    \caption{Circular interface rotation benchmark problem, the red dot is the center of rotation}
    \label{FIG:Circular_Interface_Rotation}
\end{figure}
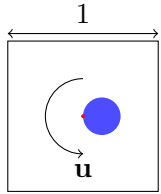

\begin{figure}
    \includegraphics[width=0.49\linewidth]{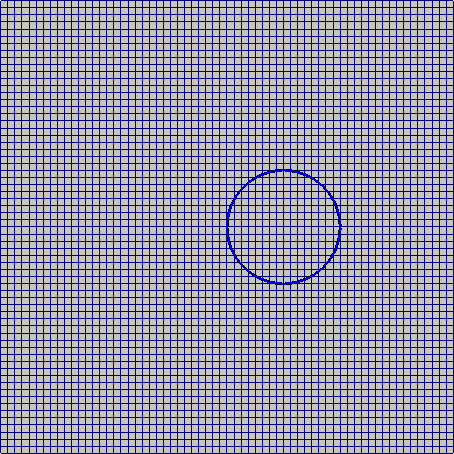}
    \includegraphics[width=0.49\linewidth]{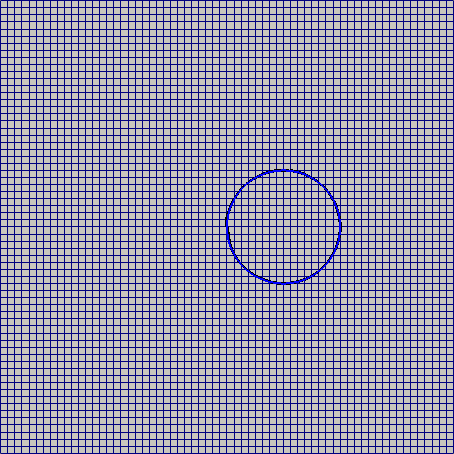}
    \caption{Comparison}
    \label{FIG:Circular_Interface_Rotation_comparison}
\end{figure}

\begin{table}
    \centering
    \begin{filecontents*}{rotCirc2.csv}
size,h,nproc,errorI,cRateI,errorF,cRateF,time(sec)
4,$2^{-4}$,8,6.067056e-03, nan,6.038973e-03, nan,4.400000e+00
5,$2^{-5}$,8,1.759717e-03,1.79,1.745156e-03,1.79,7.010000e+00
6,$2^{-6}$,8,3.968449e-04,2.15,3.927448e-04,2.15,1.910000e+01
7,$2^{-7}$,8,1.066424e-04,1.90,1.056049e-04,1.89,6.510000e+01
8,$2^{-8}$,8,2.659000e-05,2.00,2.634644e-05,2.00,3.680000e+02
9,$2^{-9}$,8,6.580688e-06,2.01,6.489521e-06,2.02,9.770000e+02
\end{filecontents*}
\pgfplotstabletypeset[
    search path={DATA_TABLES/VOF_Benchmarks},
    col sep=comma,
    columns={h, errorF, cRateF},
    columns/h/.style={
        string type,
        column name={\(h\)},
        column type/.add={|}{|},
    },
    columns/errorF/.style={
        column name={Error},
        sci zerofill,
        dec sep align={c|},
        column type/.add={}{|},
        precision=5,
    },
    columns/cRateF/.style={
        column name={Rate},
        clear infinite,
        fixed,
        fixed zerofill,
        dec sep align={c|},
        column type/.add={}{|},
        precision=2,
    },
    after row=\hline,
    every head row/.style={
        before row=\hline,
    },
    ]{rotCirc2.csv}
    \caption{Rotation of a circular interface offset from the center of rotation}
    \label{TABLE:Offset Rotation of a Circle}
\end{table}

The second benchmark problem is the advection of a circle containing composition~1 in a 
rotating velocity field as shown in
figure~\ref{FIG:Circular_Interface_Rotation}.
In this problem the angular velocity is $\pi$ radians per unit time with an end time of
$t = 2.0$. 
Note that the center of rotation is \textit{not} at the center of the circle, but rather it lies on the boundary of the circle.
Since our reconstruction and advection algorithms, are designed to be second-order accurate, we expect the approximate interface to be a second-order accurate approximation to the true circle.
The $L^1$ error in the volume fractions $f_e$ for this benchmark for six computations 
with increasing grid resolution $h = 2^{-4}, 2^{-5}, \ldots ,  2^{-9}$ is shown in
Table~\ref{TABLE:Offset Rotation of a Circle}.
In each of these computations we used a CFL number of $\sigma = \frac{1}{2}$.

\subsection{Mantle Convection Benchmark Problems}
\label{SUBSECTION:MANTLE CONVECTION BENCHMARK PROBLEMS}

In this section we compute two `benchmark' problems that are well-known and frequently used in the computational mantle convection community to demonstrate that our VOF interface tracking algorithm can reproduce previously published computational results of the same problem.
In our view the first problem, commonly known as the ``van Keken problem'' or the ``van Keken isoviscous Rayleigh-Taylor problem'' is not a reasonable `benchmark', since the problem is mathematically ill-posed. 
In other words, it is unstable~\cite{SC:1961} and perturbations due to different numerical methods can yield vastly differing results.
In fact, in~\cite{EGP-DLT-YH-HL-JM-LHK:2017} we demonstrated that it suffices to change only the algorithm with which the composition variable $C$ is advected in order to obtain clearly different results at the same output time.
For example, see Figure~11 of~\cite{EGP-DLT-YH-HL-JM-LHK:2017} or compare Figures 5(c)--(d) of~\cite{HS-ME:2010} to our results here or in~\cite{EGP-DLT-YH-HL-JM-LHK:2017} or to the results in~\cite{PEVK-SDK-HS-URC-DN-MPD:1997}.

\subsubsection{The van Keken Isoviscous Rayleigh-Taylor Problem}
\label{SUBSUBSECTION:THE VAN KEKEN ISOVISCOUS RAYLEIGH-TAYLOR PROBLEM}

In this section, we present our computation of the van Keken isoviscous Rayleigh-Taylor
problem~\cite{PEVK-SDK-HS-URC-DN-MPD:1997}.
In spite of the fact that the problem is unstable and hence ill-posed, it has become a 
standard 'benchmark' in the computational geodynamics community.
In this problem a less dense (buoyant) fluid lies beneath a denser fluid, with a perturbed
interface between the two layers.
The problem is computed in a $[D, \, 1]$ computational domain where $D = 0.9142$ is the 
width of the domain.
The initial discontinuity between the two compositional / density layers is given by
\begin{align}
  \label{EQ:van Keken Initial Condition}
    C(x,y,t=0)=
      \left\{
        \begin{array}{ll}
          0,\quad & \text{if} \quad \, 0 \, \le \, y \, < \, 0.2 \,
                           + \, 0.02 \, \cos \left(\pi \, x \, / \,D \right) \, ,\\
          1,\quad & \text{otherwise} \, .
        \end{array}
      \right.
\end{align}
This initial condition has a (discontinuous) interface along the curve
\begin{align}
  \label{EQ:The Initial van Keken Interface}
    y \, = \, 0.2 \, + \, 0.02 \, \cos \left(\frac{\pi x}{D} \right).
\end{align}

\begin{figure}
  \begin{center}
    \includegraphics[width=1.0\textwidth]{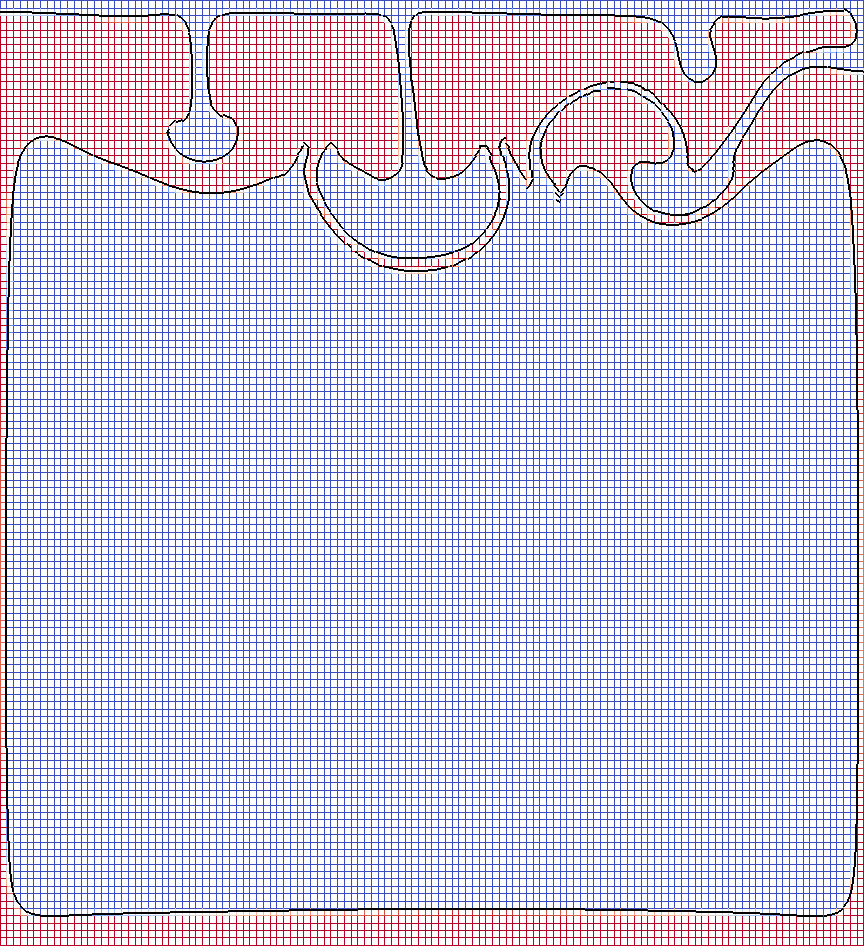}
    \caption{Computed solution of the van Keken isoviscous Rayleigh-Taylor problem
          at time $t = 2000$ on a uniform grid of $128 \times 128$ cells.
          Compare with the computational results in~\protect\cite{MK-TH-WB:2012}, \protect\cite{HS-ME:2010}, and~\protect\cite{PEVK-SDK-HS-URC-DN-MPD:1997}.
          }
    \label{FIG:van Keken Isoviscous Rayleigh-Taylor Problem}

  \end{center}
\end{figure}

\subsubsection{The Gerya-Yuen Sinking Box Benchmark}
\label{SUBSUBSECTION:THE GERYA-YUEN SINKING BOX BENCHMARK}

Following the original authors, we pose the Gerya-Yuen `sinking box' problem~\cite{TVG-DAY:2003} in dimensional form.
The problem is defined on a $500~\mathrm{km} \times 500~\mathrm{km}$ 
two-dimensional Cartesian computational domain.
A small horizontally centered $100~\mathrm{km} \times 100~\mathrm{km}$ square
is placed with its top edge $50~\mathrm{km}$ below the top of the domain so that the
the initial location and dimension of the box is defined by the composition field $C(\mathbf{x},t)$ as follows:
\begin{align}
  \label{cd:paras2}
    C(\mathbf{x},0)=
      \left \{
        \begin{matrix}
          1, \quad & \text{ if } \, (x,y )\in [200\textrm{ km},300\textrm{ km}]\times[350\textrm{ km},450\textrm{ km}] \, , \\
         0,\quad &\text{otherwise} \, .
        \end{matrix} \right .
\end{align}
The block's density is $\rho_1 = 3300 \textrm{ kg/m}^3$, while the background
density is $\rho_1=3300 \textrm{ kg/m}^3$.
We approximate the solution of the incompressible Stokes equations
(i.e.,equations~\eqref{EQ:Dimensional Conservation of Mass}--\eqref{EQ:Dimensional Stokes Equation for v} without the term $\rho_0 \alpha (T - T_0) g$ in equation~\eqref{EQ:Dimensional Stokes Equation for v}) with these initial conditions and holding the following parameters fixed:
\begin{align}
\label{cd:paras1}
  \begin{array}{llll}
    {\bf g}=(0,9.8) \textrm{ m/s}^2, & \textrm{acceleration due to gravity} \\
          L=500 \textrm{ km} & \textrm{domain height and width}             \\
        \mu=10^{21} \textrm{ Pa}\cdot \textrm{s} & \textrm{viscosity}       \\
     \rho_0=3200\textrm{ kg/m}^3, & \textrm{background density}             \\
     \rho_1=3300 \textrm{ kg/m}^3, & \textrm{small box density}             \\
  \end{array}
\end{align}




\begin{figure}
  \begin{center}
    \begin{subfigure}[t]{0.45\textwidth}
    	\centering
      \includegraphics[width=1.0\textwidth]{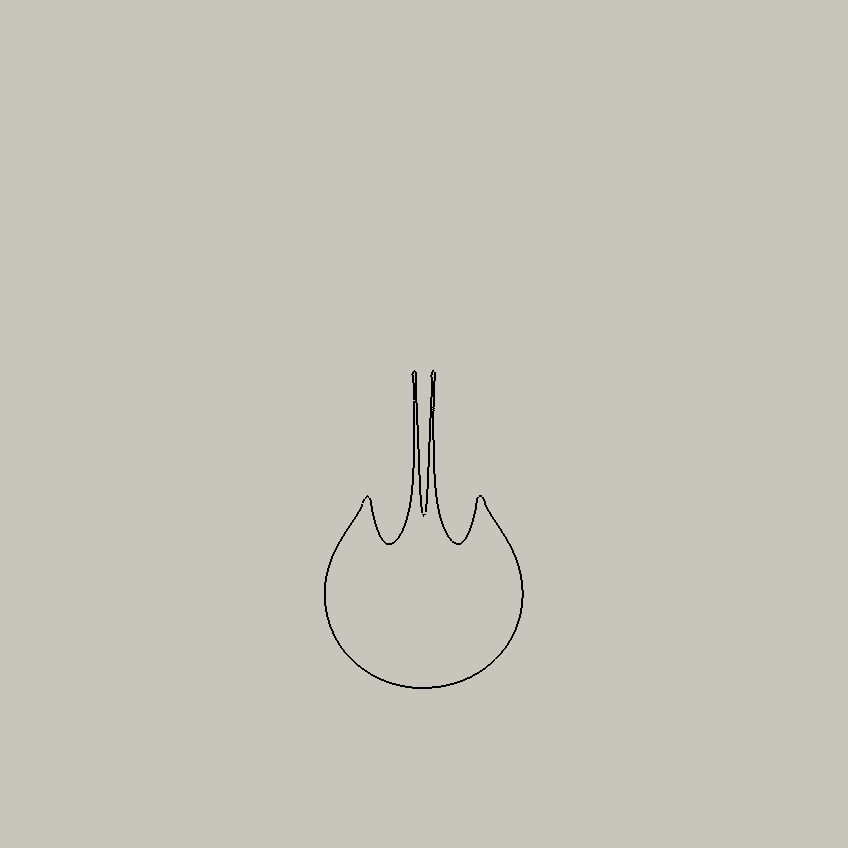}
      \caption{\small The interface against a tan background.}
      \label{FIG:SINKING BOX BENCHMARK WITH NO GRID}
    \end{subfigure}
    \quad
    \begin{subfigure}[t]{0.45\textwidth}  
    	\centering   
      \includegraphics[width=1.00\textwidth]{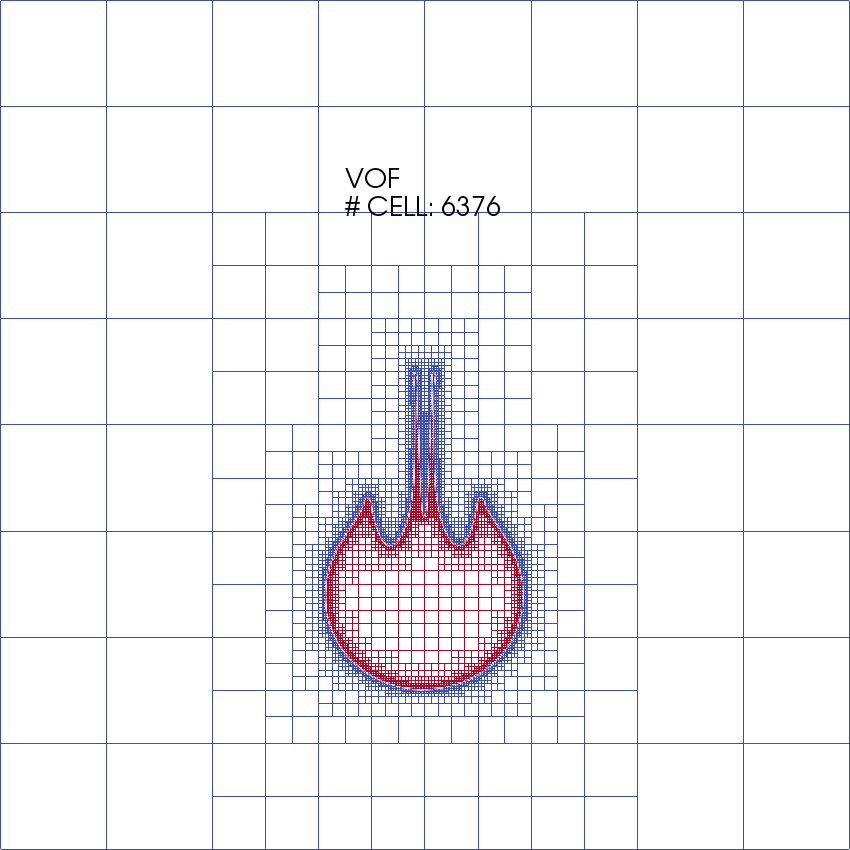}
      \caption{The interface together with the underlying AMR grid.}
      \label{FIG:SINKING BOX BENCHMARK WITH AMR}
    \end{subfigure}
    \caption{Fluid interface for the Gerya-Yuen~\protect\cite{TVG-DAY:2003} `sinking box' problem at time $t = 9.81$ Myr computed with AMR as shown on the right.}    
    \label{FIG:SINKING BOX BENCHMARK}
  \end{center}
\end{figure}

\subsection{Computations of Thermochemical Convection in a Density Stratified Fluid}
\label{SUBSECTION:THERMOCHEMICAL CONVECTION IN A DENSITY STRATIFIED FLUID}

We now present the results of our computations of the model problem for
thermochemical convection with density stratification, the equations for which were presented in
Section~\ref{SECTION: THERMOCHEMICAL CONVECTION WITH DENSITY STRATIFICATION}.
In these computations the Rayleigh number is fixed at $\mathrm{Ra} = 10^5$ and we vary only the buoyancy ratio as follows: $\mathrm{B} = 0.0, \, 0.1, \, 0.2, \, \ldots , \, 1.0$ and
$\mathrm{B} = 2.0$.
The domain for all of the computational results shown below is a two-dimensional rectangular region that we denote by $\Omega= [0,3]\times[0,1]$ as shown in
Figure~\ref{Fig:The Nondimensional Computational Domain}.

The initial conditions for the temperature $T$ are,
\begin{align}
    \label{EQ: OUT-OF-PHASE TEMPERATURE INITIAL CONDITIONS}
    T(\mathbf{x}, 0) =
    \left\{
    \begin{array}{ll}
        (1 - 5 \, y) \, + \, A \, \sin(10 \, \pi \, y ) \, (1 - \cos(\frac{2}{3} \, k \, \pi \, x ))
        & \text{if} \quad 0 \le y \le \frac{1}{10}    \, , \\
        (5 - 5 \, y) \,  + \, A \, \sin(10 \, \pi \, y ) \,
        (1 - \cos(\frac{2}{3} \, k \, \pi \, x \, + \, \pi )) \,
        &    \text{if} \quad \frac{9}{10} \le y \le 1 \, , \\
        0.5
        & \phantom{if} \quad \text{otherwise}         \; , \\
    \end{array}
    \right .
\end{align}
where the period of the perturbation $k \, = \, 1.5$ and the amplitude of the perturbation
$A \, = \, 0.05$.
Note that $A \, = \, 0.05$ ensures that $0 \le \, T(x, y;0)\, \le \, 1$ throughout  the
entire computational domain.
The initial conditions for the composition are,
\begin{align}
    \label{EQ: NUMERICAL COMPOSITION INITIAL CONDITION}
    C(x,y;t=0) &=
    \left\{
    \begin{array}{ll}
        1 & \text{if} \quad 0 \le y < \frac{1}{2} \, ,\\
        0 & \text{if} \quad \frac{1}{2} \le y \le 1
    \end{array}
    \right .
\end{align}
and the boundary conditions for the velocity and temperature are as specified
in~\eqref{EQ:No-FLow BCs}--\eqref{EQ:Temperature BCs at x = 3}.

All of the results shown below were computed twice: once on a fixed, uniform grid with $192 \times 64$ square cells each with side $h= 64^{-1}$ and then on the same underlying grid but with the addition of two levels of an adaptively refined mesh, on and in, a neighborhood of the interface.  
Each level of refinement increases the grid resolution by a factor of two; i.e.,
$h \to \frac{h}{4}$ with two levels of refinement.

\begin{figure}
    \begin{center}
        \begin{subfigure}[b]{1.0\textwidth}
            \includegraphics[width=1.0\textwidth]{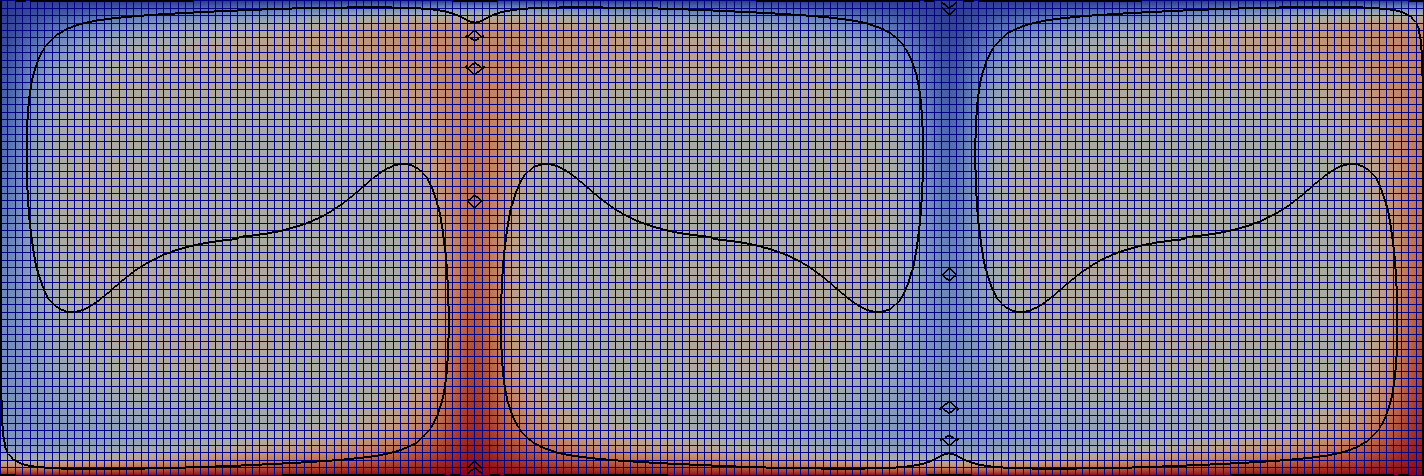}
            \caption{$\mathrm{B} = 0.0$ at $t' = 1.97 \cdot 10^{-2}$ ON A UNIFORM 
                GRID of $196 \times 64$ cells.}
            \label{FIG:B = 0.0 at t' = 0.0197 ON A UNIFORM GRID}
        \end{subfigure}

        \begin{subfigure}[b]{1.0\textwidth}
            \includegraphics[width=1.0\textwidth]{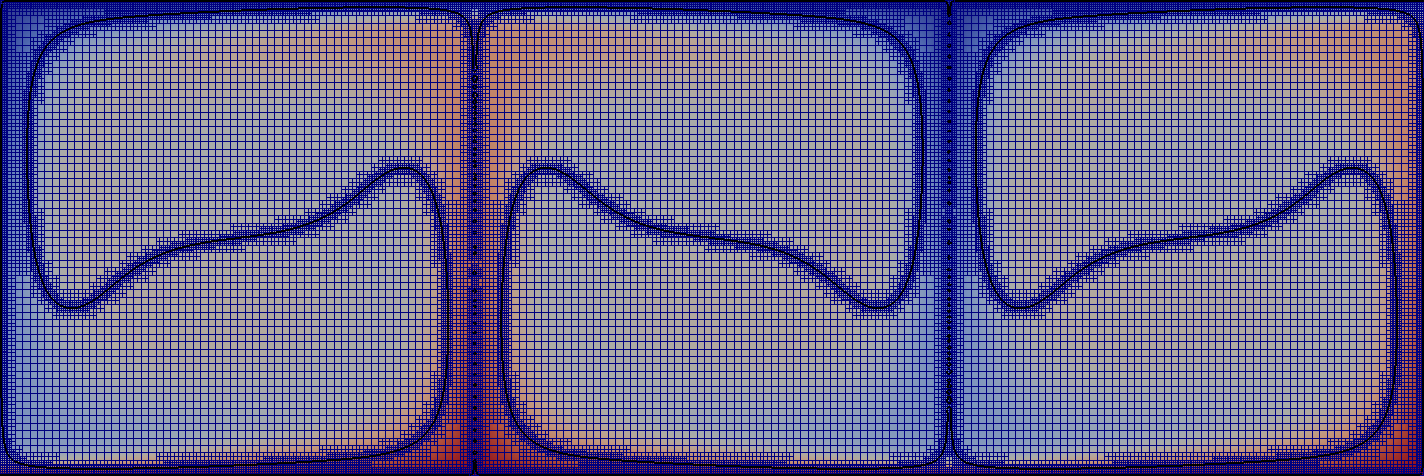}
            \caption{$\mathrm{B} =  0.0$ at $t' = 1.97\cdot10^{-2}$ with two levels of AMR}
            \label{FIG:B = 0.0 at t' = 0.0197 with AMR}
        \end{subfigure}

        \begin{subfigure}[b]{1.0\textwidth}
            \includegraphics[width=1.0\textwidth]{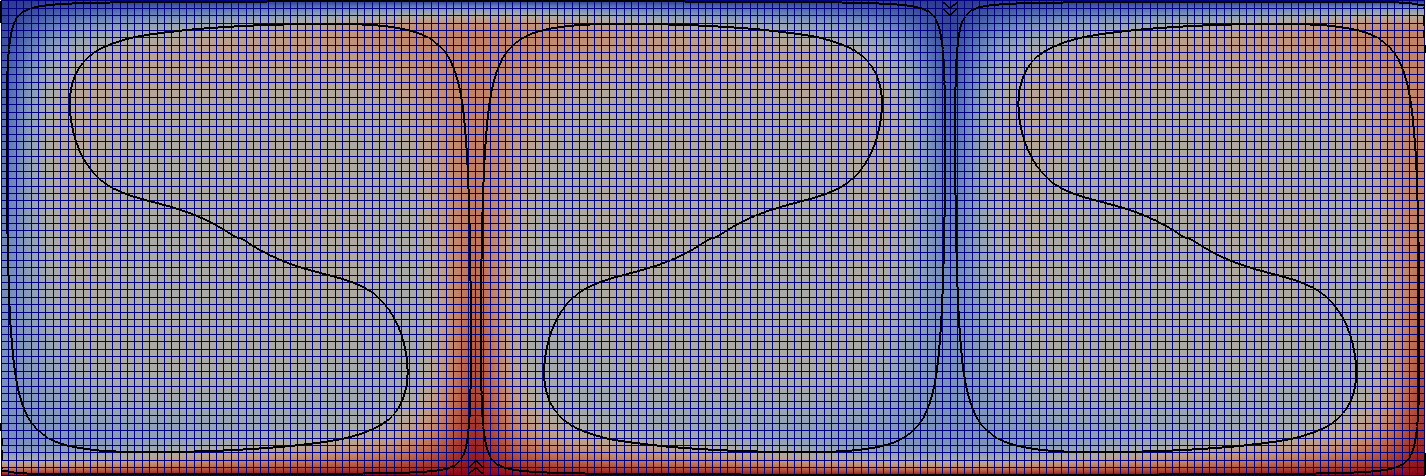}
            \caption{$B = 0.0$ at $t' = 2.36 \cdot 10^{-2}$ ON A UNIFORM GRID of 
                $196 \times 64$ cells.}
            \label{FIG:B = 0.0 at t' = 0.0236  ON A UNIFORM GRID}
        \end{subfigure}

        \begin{subfigure}[b]{1.0\textwidth}
            \includegraphics[width=1.0\textwidth]{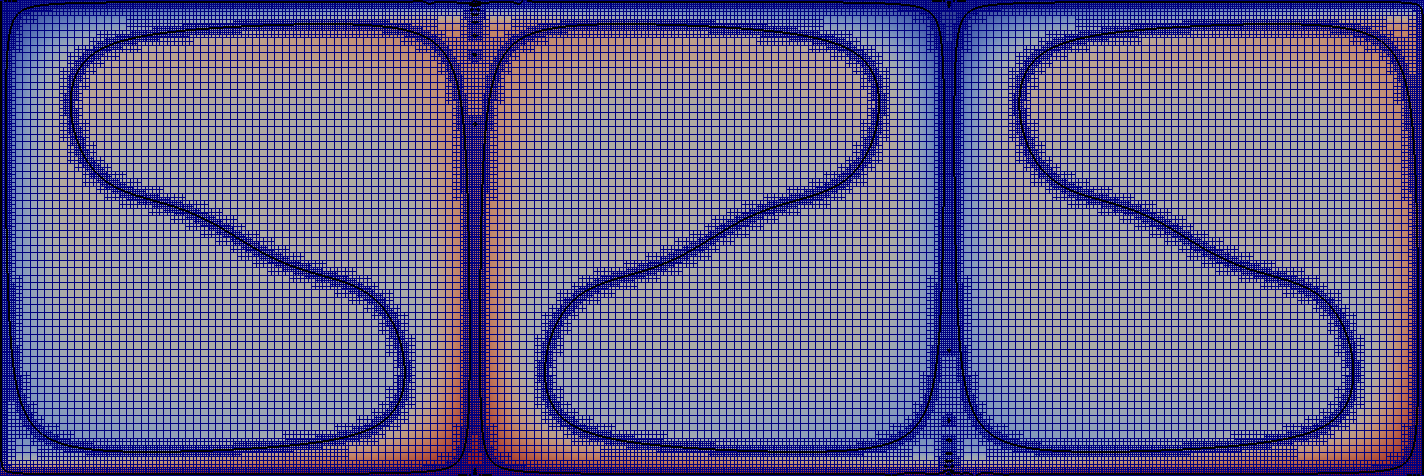}
            \caption{$\mathrm{B} = 0.0$ at $t'=2.36\cdot10^{-2}$ with two levels of AMR}
            \label{FIG:B = 0.0 at t' = 0.0236 WITH AMR}
        \end{subfigure}
        \caption{Computations with $\mathrm{B} = 0.0$ and $\mathrm{Ra} = 10^5$ on an 
            underlying uniform grid of $196 \times 64$ square cells at
            $t' = 1.97 \cdot 10^{-2}$ and $t' = 2.36 \cdot 10^{-2}$. 
        	The background color is the temperature, which varies from $T = 0.0$ (dark blue) to $T = 1.0$ (dark red).
                 }
        \label{FIG:B = 0.0}
    \end{center}
\end{figure} 

\begin{figure}
    \begin{center}
        \begin{subfigure}[b]{1.0\textwidth}
            \includegraphics[width=1.0\textwidth]{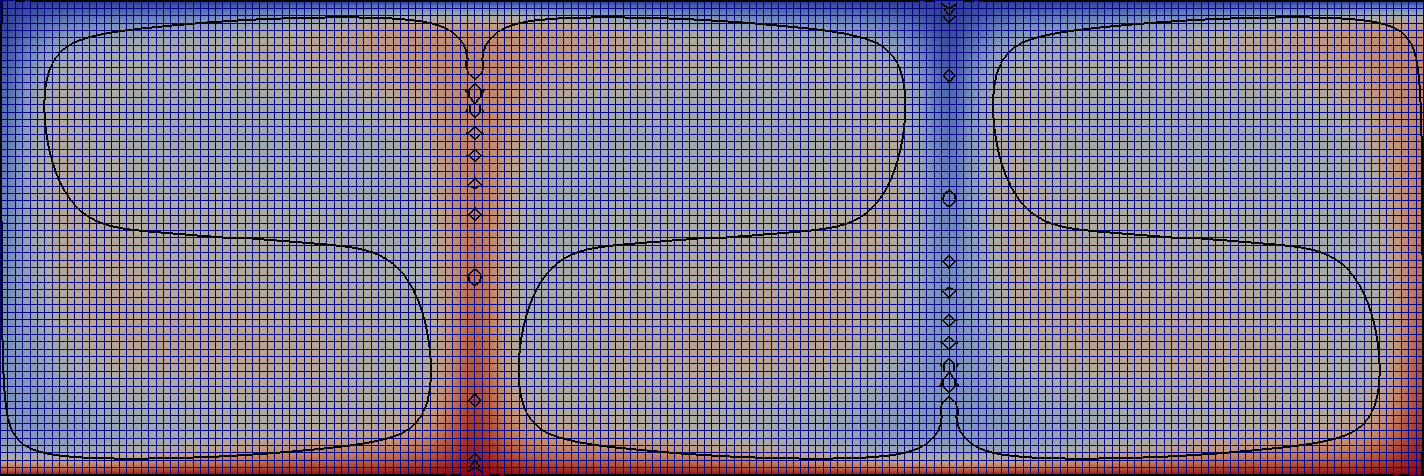}
            \caption{$\mathrm{B} = 0.1$ at $t' = 1.97\cdot10^{-2}$ ON A UNIFORM GRID of $196 \times 64$ cells.}
            \label{FIG:B = 0.1 at t' = 0.0197  ON A UNIFORM GRID}
        \end{subfigure}

        \begin{subfigure}[b]{1.0\textwidth}
            \includegraphics[width=1.0\textwidth]{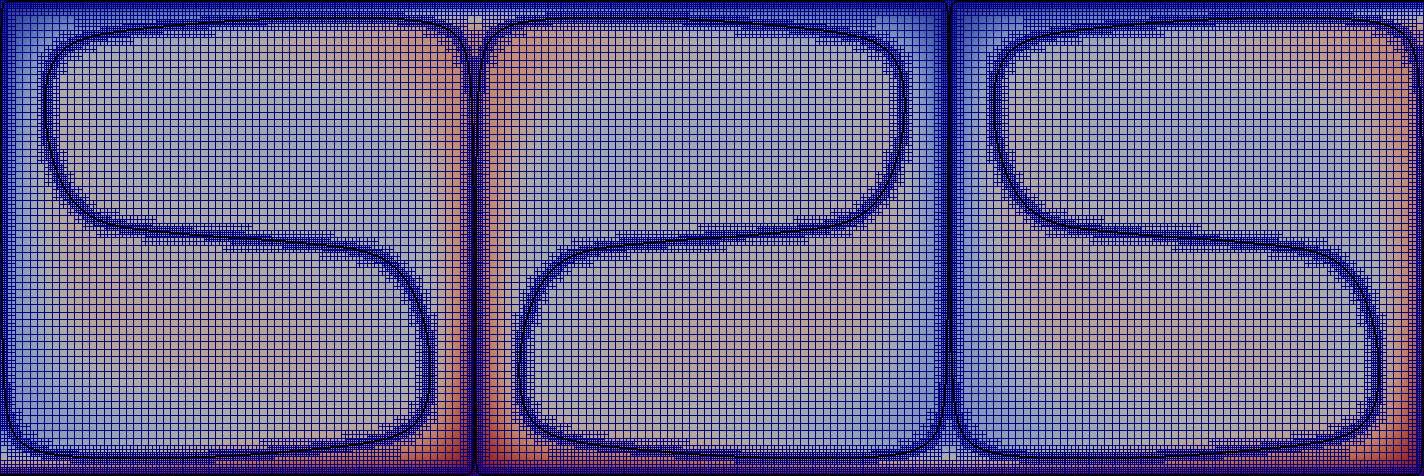}
            \caption{$\mathrm{B} =  0.1$ at $t' = 1.97\cdot10^{-2}$ with two levels of AMR}
            \label{FIG:B = 0.1 at t' = 0.0197 with AMR}
        \end{subfigure}

        \begin{subfigure}[b]{1.0\textwidth}
            \includegraphics[width=1.0\textwidth]{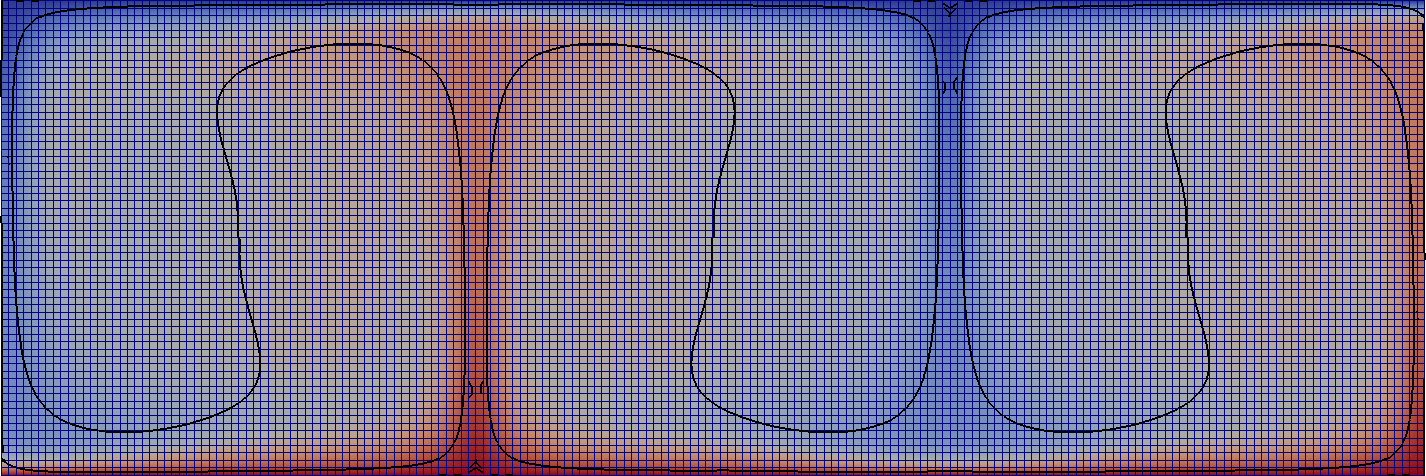}
            \caption{$B = 0.1$ at $t'=2.36\cdot10^{-2}$ ON A UNIFORM GRID of $196 \times 64$ cells.}
            \label{FIG:B = 0.1 at t' = 0.0236 ON A UNIFORM GRID}
        \end{subfigure}

        \begin{subfigure}[b]{1.0\textwidth}
            \includegraphics[width=1.0\textwidth]{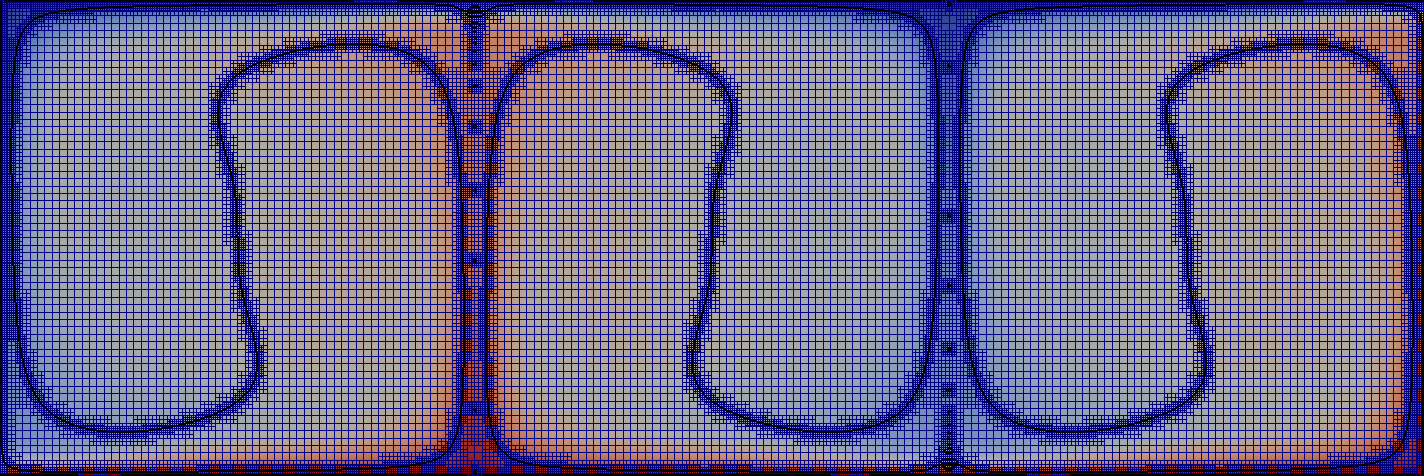}
            \caption{$\mathrm{B} = 0.1$ at $t'=2.36\cdot10^{-2}$ with two levels of AMR}
            \label{FIG:B = 0.1 at t' = 0.0236 WITH AMR}
        \end{subfigure}
         \caption{Computations with $\mathrm{B} = 0.1$ and $\mathrm{Ra} = 10^5$ on 
             an underlying uniform grid of $196 \times 64$ square cells at
             $t' = 1.97 \cdot 10^{-2}$ and $t' = 2.36 \cdot 10^{-2}$. 
         	 The background color is the temperature, which varies from $T = 0.0$ (blue) to $T = 1.0$ (dark red).
                 }
        \label{FIG:B = 0.1}
    \end{center}
\end{figure} 

\begin{figure}
    \begin{center}
        \begin{subfigure}[b]{1.0\textwidth}
            \includegraphics[width=1.0\textwidth]{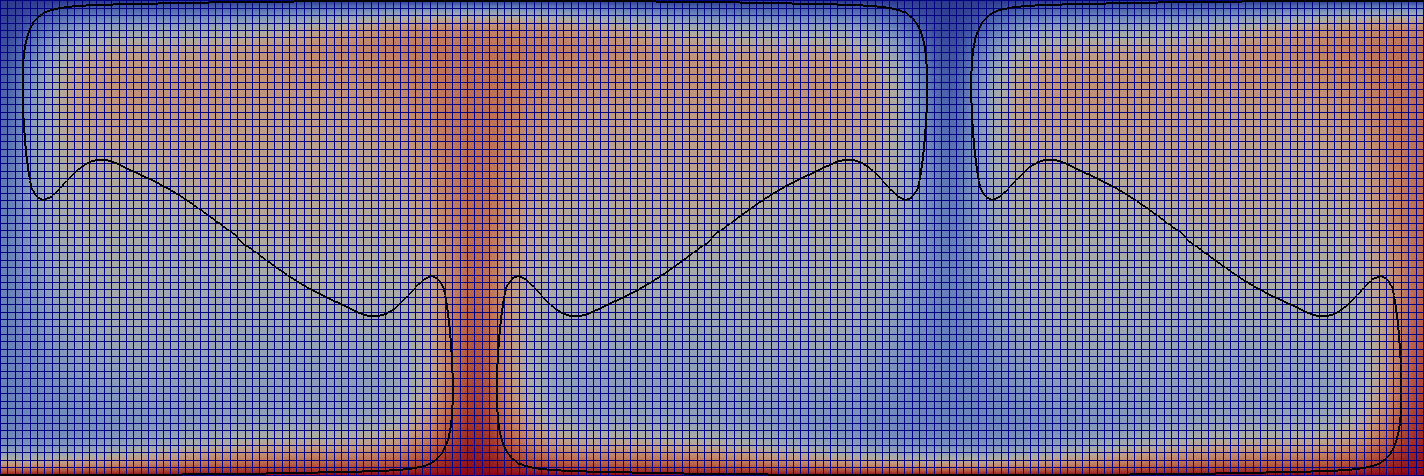}
            \caption{$\mathrm{B} = 0.2$ at $t' = 1.97\cdot10^{-2}$ ON A UNIFORM GRID of $196 \times 64$ cells.}
            \label{FIG:B = 0.2 at t' = 0.0197 on a uniform grid}
        \end{subfigure}

        \begin{subfigure}[b]{1.0\textwidth}
            \includegraphics[width=1.0\textwidth]{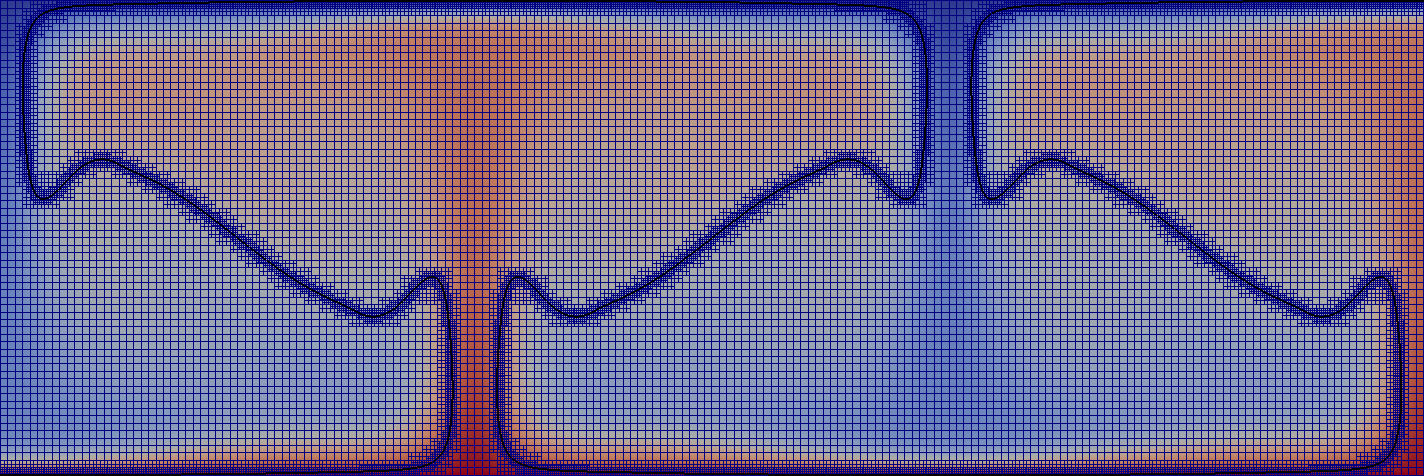}
            \caption{$\mathrm{B} =  0.2$ at $t' = 1.97\cdot10^{-2}$ with two levels of AMR}
            \label{FIG:B = 0.2 at t' = 0.0197 with AMR}
        \end{subfigure}

        \begin{subfigure}[b]{1.0\textwidth}
            \includegraphics[width=1.0\textwidth]{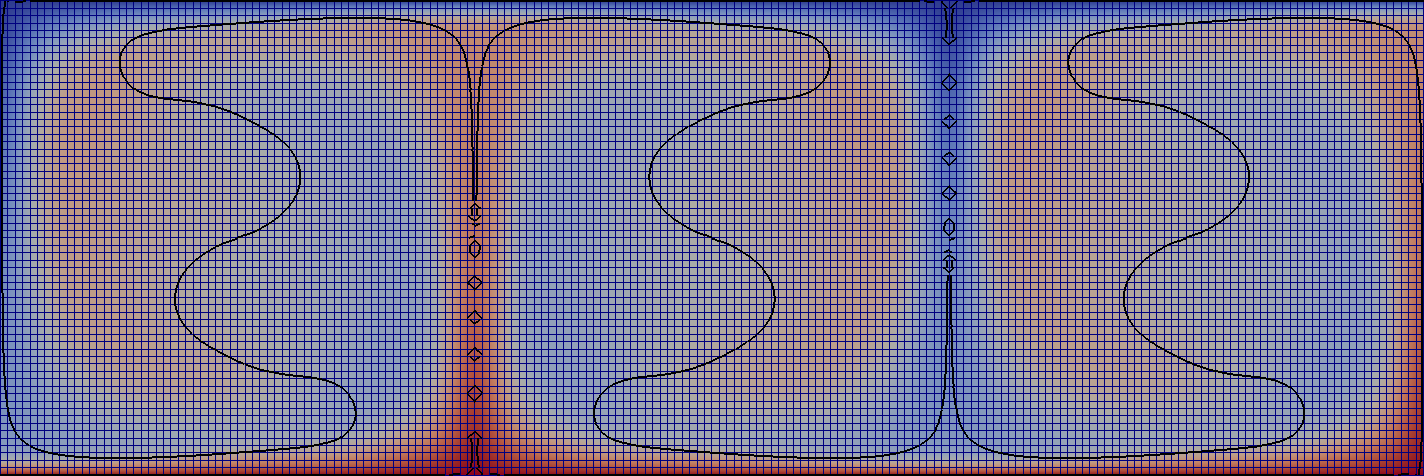}
            \caption{$B = 0.2$ at $t'=2.36\cdot10^{-2}$ ON A UNIFORM GRID of $196 \times 64$ cells.}
            \label{FIG:B = 0.2 at t' = 0.0236 ON A UNIFORM GRID}
        \end{subfigure}

        \begin{subfigure}[b]{1.0\textwidth}
            \includegraphics[width=1.0\textwidth]{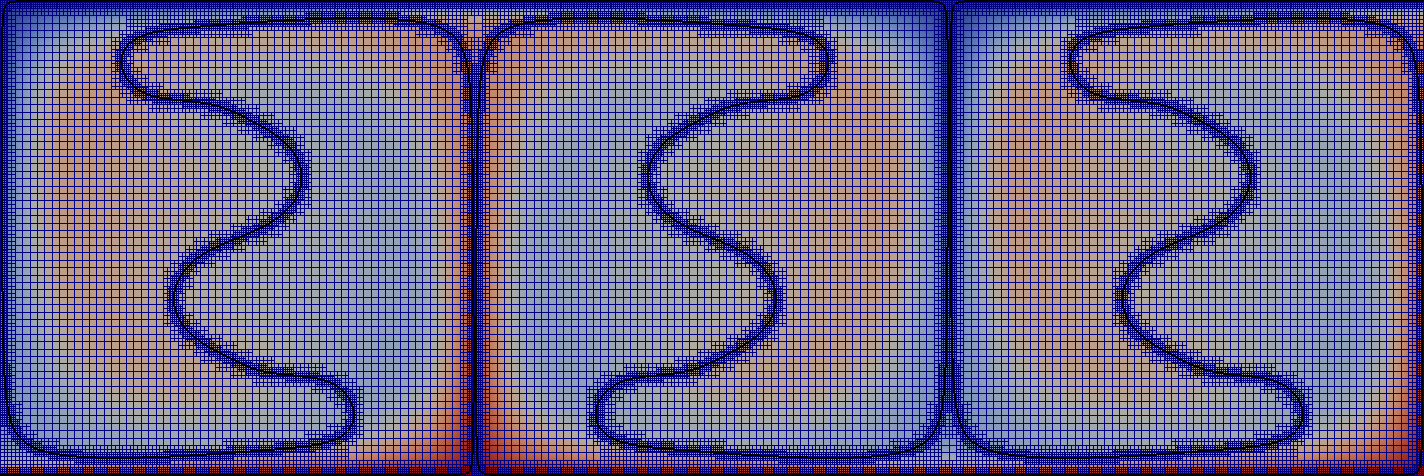}
            \caption{$\mathrm{B} = 0.2$ at $t'=2.36\cdot10^{-2}$ with two levels of AMR}
            \label{FIG:B = 0.2 at t' = 0.0236 WITH AMR}
        \end{subfigure}
         \caption{Computations with $\mathrm{B} = 0.2$ and $\mathrm{Ra} = 10^5$ on 
             an underlying uniform grid of $196 \times 64$ square cells at
             $t' = 1.97 \cdot 10^{-2}$ and $t' = 2.36 \cdot 10^{-2}$. 
         	 The background color is the temperature, which varies from $T = 0.0$ (blue) to $T = 1.0$ (dark red).
                 }
        \label{FIG:B = 0.2}
    \end{center}
\end{figure} 

\begin{figure}
    \begin{center}
        \begin{subfigure}[b]{1.0\textwidth}
            \includegraphics[width=1.0\textwidth]{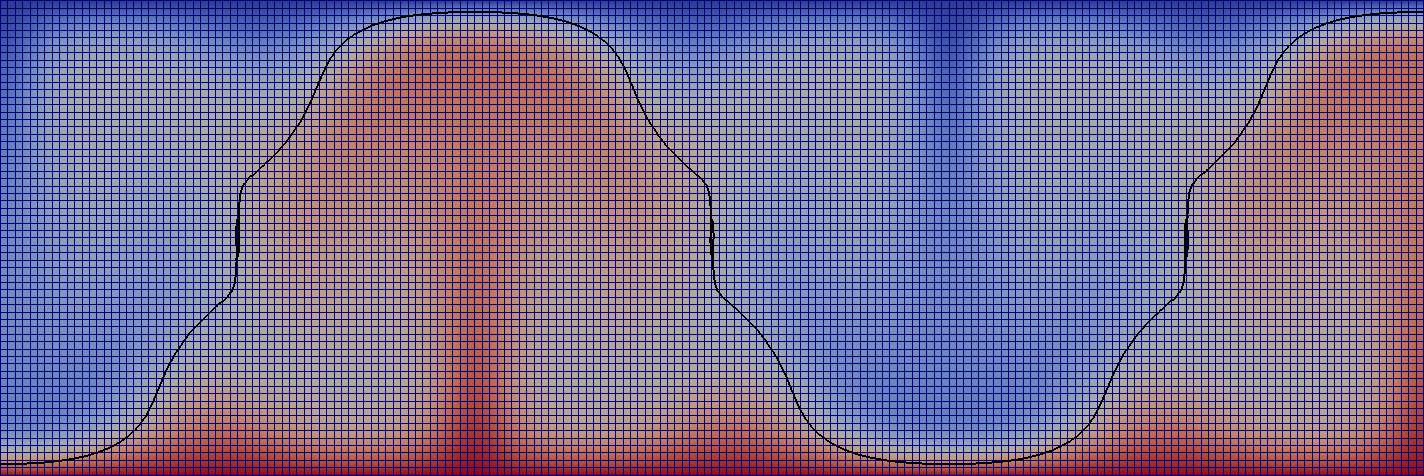}
            \caption{$\mathrm{B} = 0.3$ at $t' = 1.97\cdot10^{-2}$ ON A UNIFORM GRID of $196 \times 64$ cells.}
            \label{FIG:B = 0.3 at t' = 1.97 ON A UNIFORM GRID}
        \end{subfigure}

        \begin{subfigure}[b]{1.0\textwidth}
            \includegraphics[width=1.0\textwidth]{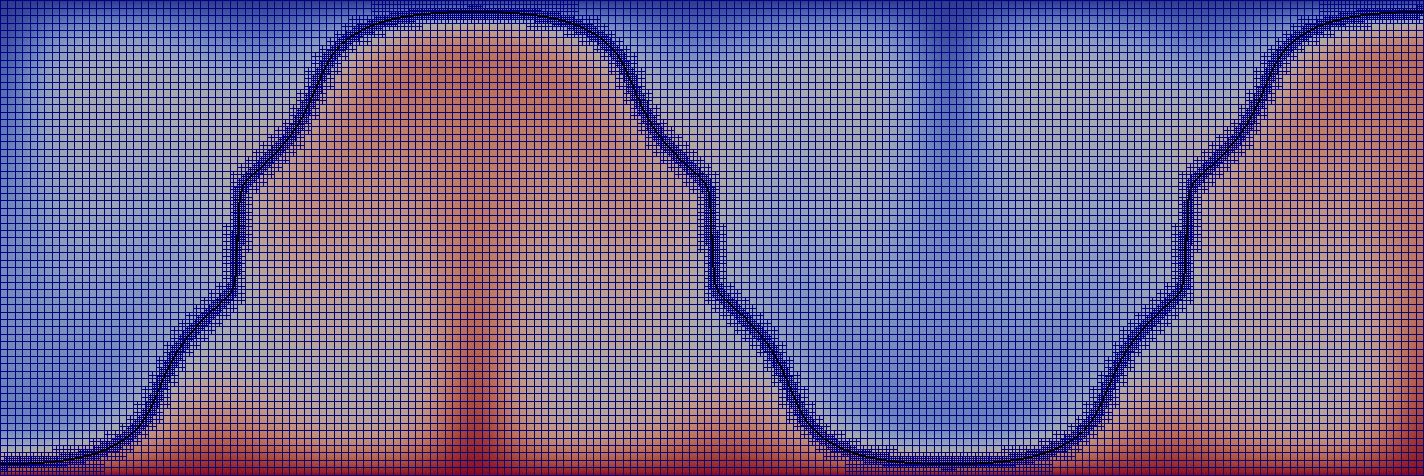}
            \caption{$\mathrm{B} =  0.3$ at $t' = 1.97\cdot10^{-2}$ with two levels of AMR}
            \label{FIG:B = 0.3 at t' = 1.97 WITH AMR}
        \end{subfigure}

        \begin{subfigure}[b]{1.0\textwidth}
            \includegraphics[width=1.0\textwidth]{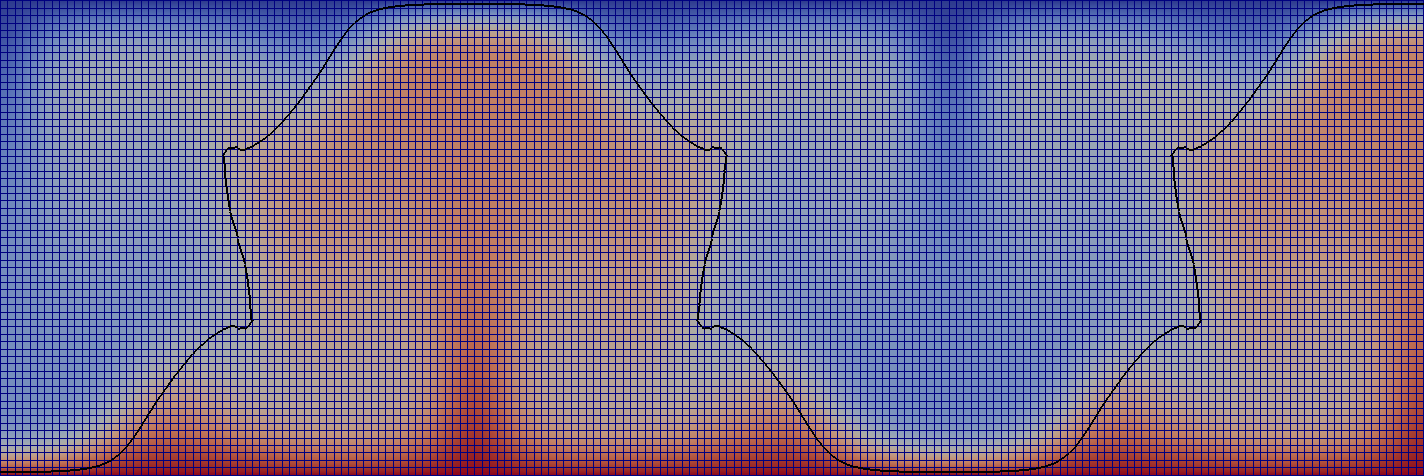}
            \caption{$B = 0.3$ at $t'=2.36\cdot10^{-2}$ ON A UNIFORM GRID of $196 \times 64$ cells.}
            \label{FIG:B = 0.3 at t' = 2.36 ON A UNIFORM GRID}
        \end{subfigure}

        \begin{subfigure}[b]{1.0\textwidth}
            \includegraphics[width=1.0\textwidth]{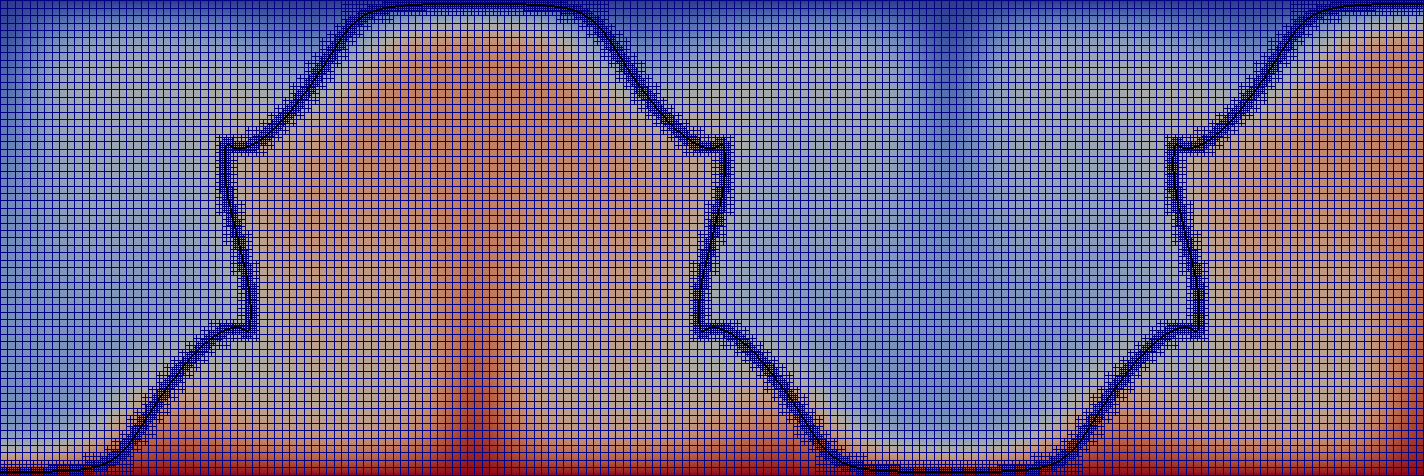}
            \caption{$\mathrm{B} = 0.3$ at $t'=2.36\cdot10^{-2}$ with two levels of AMR}
            \label{FIG:B = 0.3 at t' = 2.36 WITH AMR}
        \end{subfigure}
         \caption{Computations with $\mathrm{B} = 0.3$ and $\mathrm{Ra} = 10^5$ on 
             an underlying uniform grid of $196 \times 64$ square cells at
             $t' = 1.97 \cdot 10^{-2}$ and $t' = 2.36 \cdot 10^{-2}$. 
         	 The background color is the temperature, which varies from $T = 0.0$ (dark blue) to $T = 1.0$ (dark red).
                 }
        \label{FIG:B = 0.3}
    \end{center}
\end{figure} 

\begin{figure}
    \begin{center}
        \begin{subfigure}[b]{1.0\textwidth}
            \includegraphics[width=1.0\textwidth]{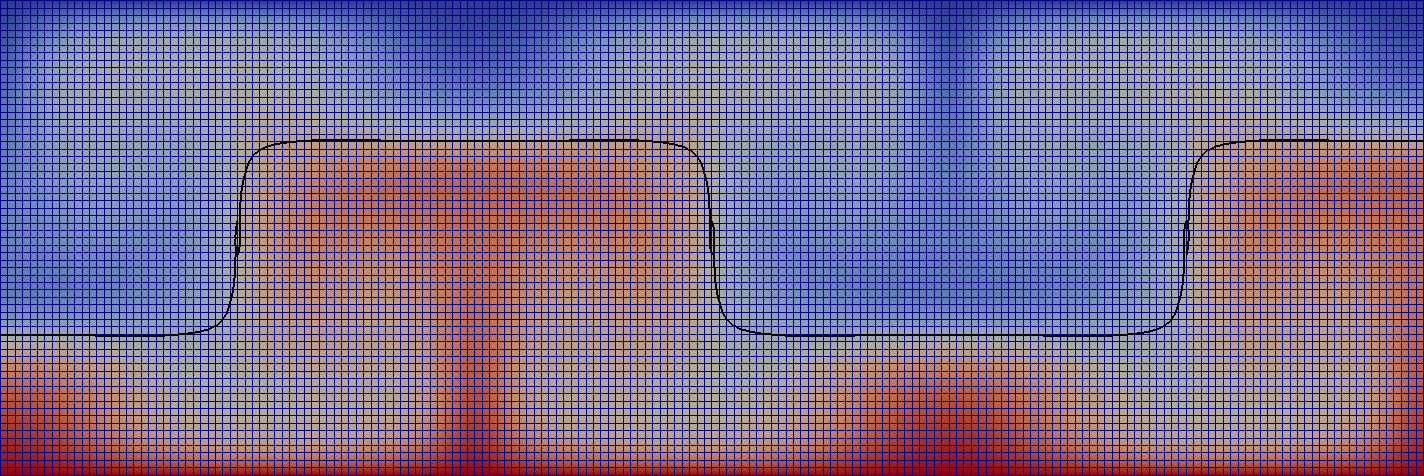}
            \caption{$\mathrm{B} = 0.4$ at $t' = 1.97\cdot10^{-2}$ ON A UNIFORM GRID of $196 \times 64$ cells.}
            \label{FIG:B = 0.4 at t' = 1.97 ON A UNIFORM GRID}
        \end{subfigure}

        \begin{subfigure}[b]{1.0\textwidth}
            \includegraphics[width=1.0\textwidth]{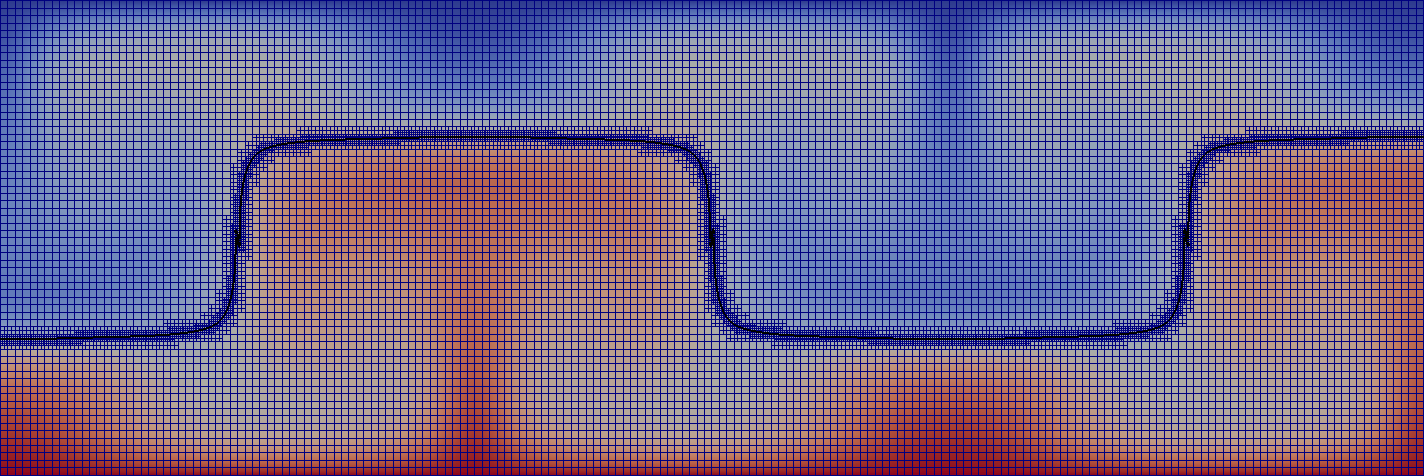}
            \caption{$\mathrm{B} =  0.4$ at $t' = 1.97\cdot10^{-2}$ with two levels of AMR}
            \label{FIG:B = 0.4 at t' = 1.97 WITH AMR}
        \end{subfigure}

        \begin{subfigure}[b]{1.0\textwidth}
            \includegraphics[width=1.0\textwidth]{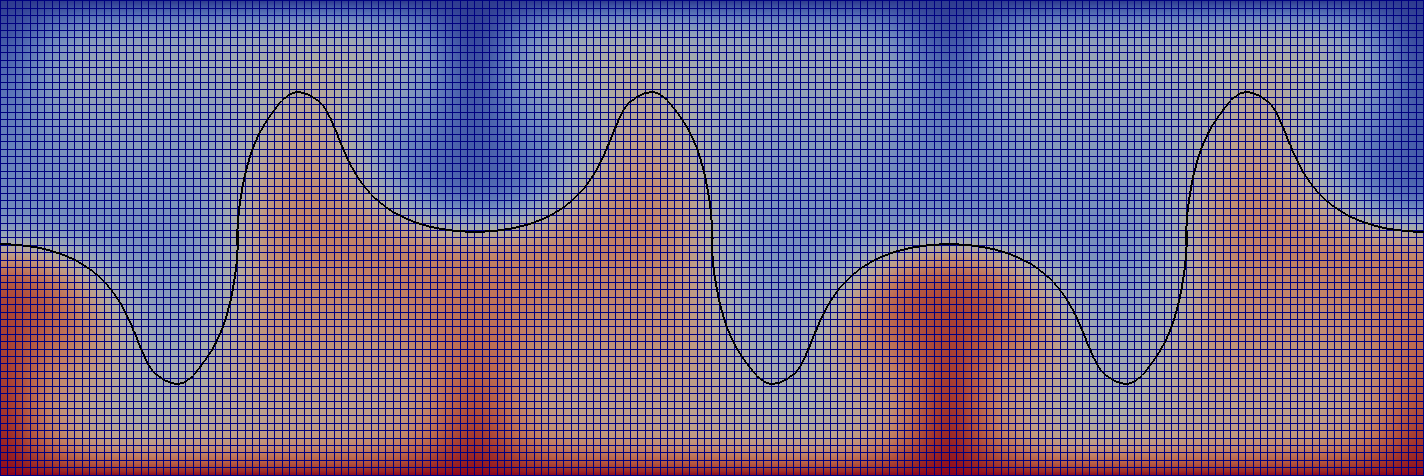}
            \caption{$B = 0.4$ at $t'=2.36\cdot10^{-2}$ ON A UNIFORM GRID of $196 \times 64$ cells.}
            \label{FIG:B = 0.4 at t' = 2.36 ON A UNIFORM GRID}
        \end{subfigure}

        \begin{subfigure}[b]{1.0\textwidth}
            \includegraphics[width=1.0\textwidth]{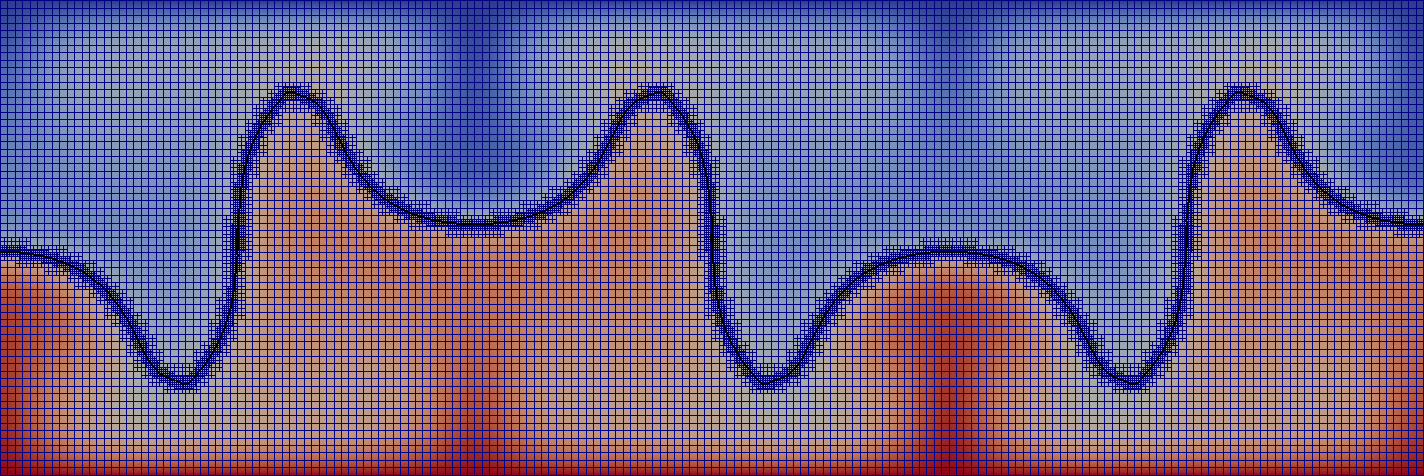}
            \caption{$\mathrm{B} = 0.4$ at $t'=2.36\cdot10^{-2}$ with two levels of AMR}
            \label{FIG:B = 0.4 at t' = 2.36 WITH AMR}
        \end{subfigure}
         \caption{Computations with $\mathrm{B} = 0.4$ and $\mathrm{Ra} = 10^5$ on 
             an underlying uniform grid of $196 \times 64$ square cells  at
             $t' = 1.97 \cdot 10^{-2}$ and $t' = 2.36 \cdot 10^{-2}$. 
         	 The background color is the temperature, which varies from $T = 0.0$ (dark blue) to $T = 1.0$ (dark red).
                  }
        \label{FIG:B = 0.4}
    \end{center}
\end{figure} 

\begin{figure}
    \begin{center}
        \begin{subfigure}[b]{1.0\textwidth}
            \includegraphics[width=1.0\textwidth]{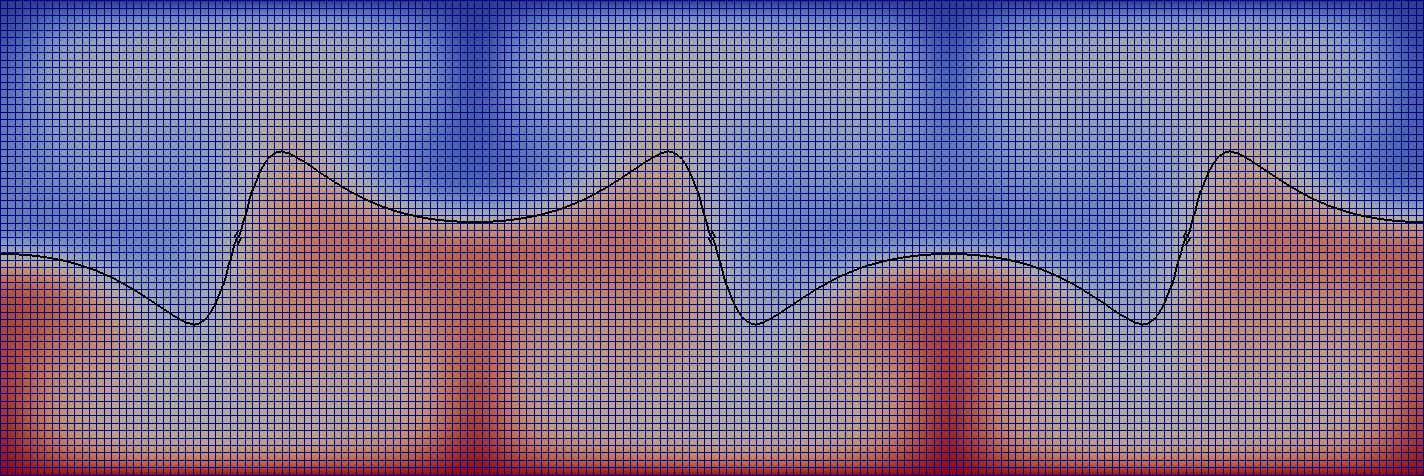}
            \caption{$\mathrm{B} = 0.5$ at $t' = 1.97\cdot10^{-2}$ ON A UNIFORM GRID of $196 \times 64$ cells.}
            \label{FIG:B = 0.5 at t' = 0.0197 on a uniform grid}
        \end{subfigure}

        \begin{subfigure}[b]{1.0\textwidth}
            \includegraphics[width=1.0\textwidth]{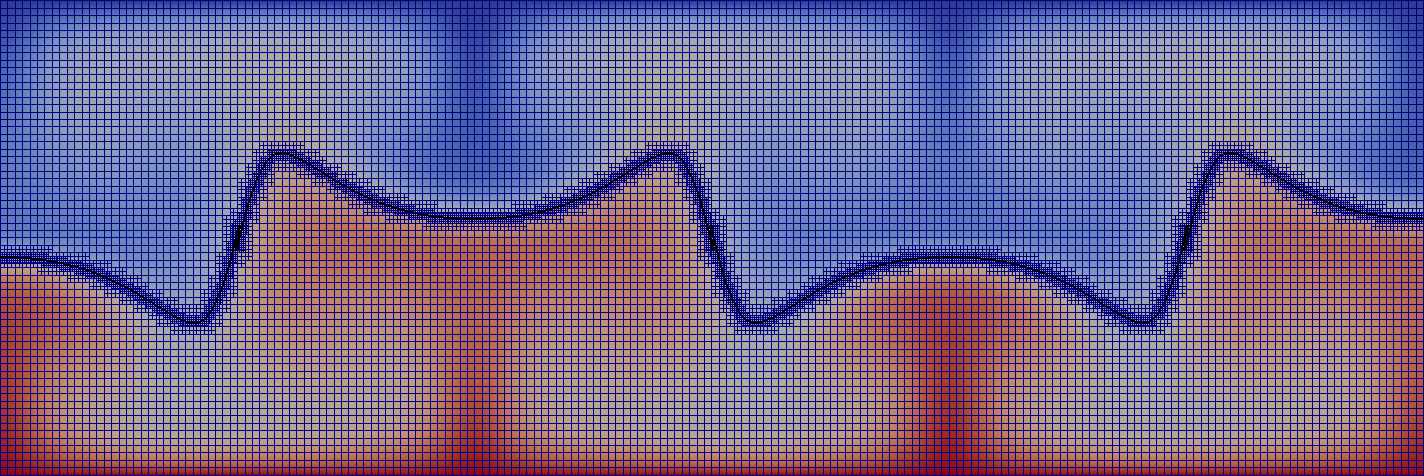}
            \caption{$\mathrm{B} =  0.5$ at $t' = 1.97\cdot10^{-2}$ with two levels of AMR}
            \label{FIG:B = 0.5 at t' = 0.0197 with AMR}
        \end{subfigure}

        \begin{subfigure}[b]{1.0\textwidth}
            \includegraphics[width=1.0\textwidth]{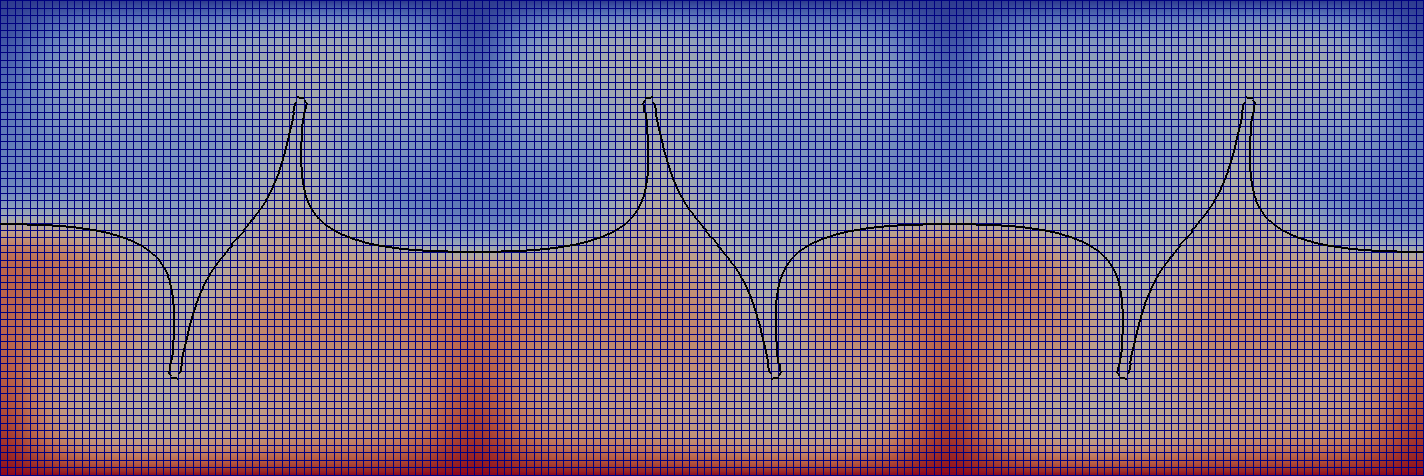}
            \caption{$B = 0.5$ at $t'=2.36\cdot10^{-2}$ ON A UNIFORM GRID of $196 \times 64$ cells.}
            \label{FIG:B = 0.5 at t' = 0.0236 ON A UNIFORM GRID}
        \end{subfigure}

        \begin{subfigure}[b]{1.0\textwidth}
            \includegraphics[width=1.0\textwidth]{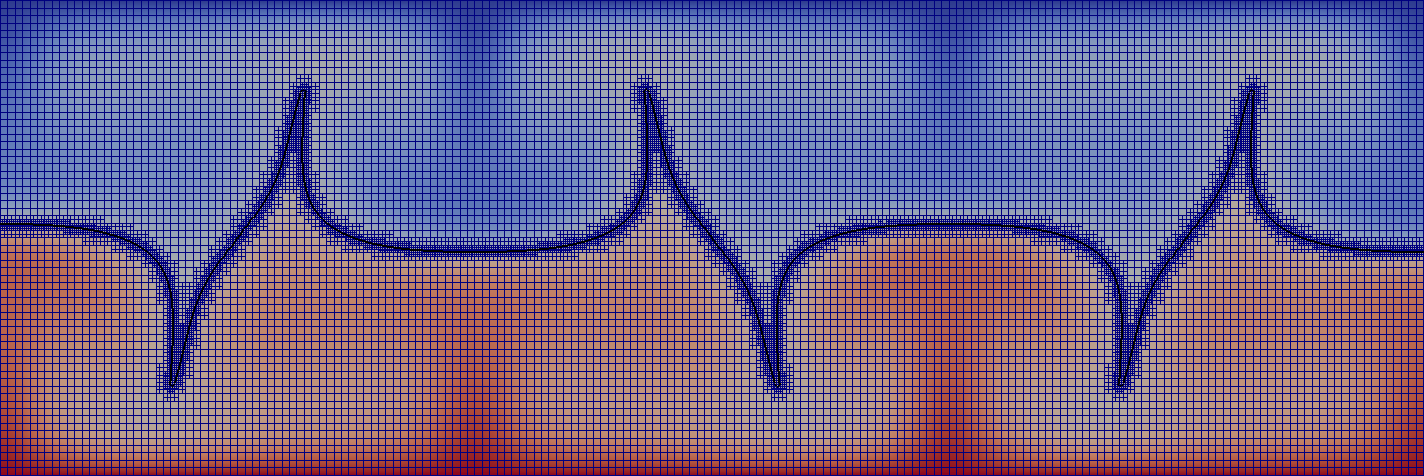}
            \caption{$\mathrm{B} = 0.5$ at $t'=2.36\cdot10^{-2}$ with two levels of AMR}
            \label{FIG:B = 0.5 at t' = 0.0236 WITH AMR}
        \end{subfigure}
         \caption{Computations with $\mathrm{B} = 0.5$ and $\mathrm{Ra} = 10^5$ on 
             an underlying uniform grid of $196 \times 64$ square cells at
             $t' = 1.97 \cdot 10^{-2}$ and $t' = 2.36 \cdot 10^{-2}$. 
         	 The background color is the temperature, which varies from $T = 0.0$ (dark blue) to $T = 1.0$ (dark red).
                 }
        \label{FIG:B = 0.5}
    \end{center}
\end{figure} 

\begin{figure}
    \begin{center}
        \begin{subfigure}[b]{1.0\textwidth}
            \includegraphics[width=1.0\textwidth]{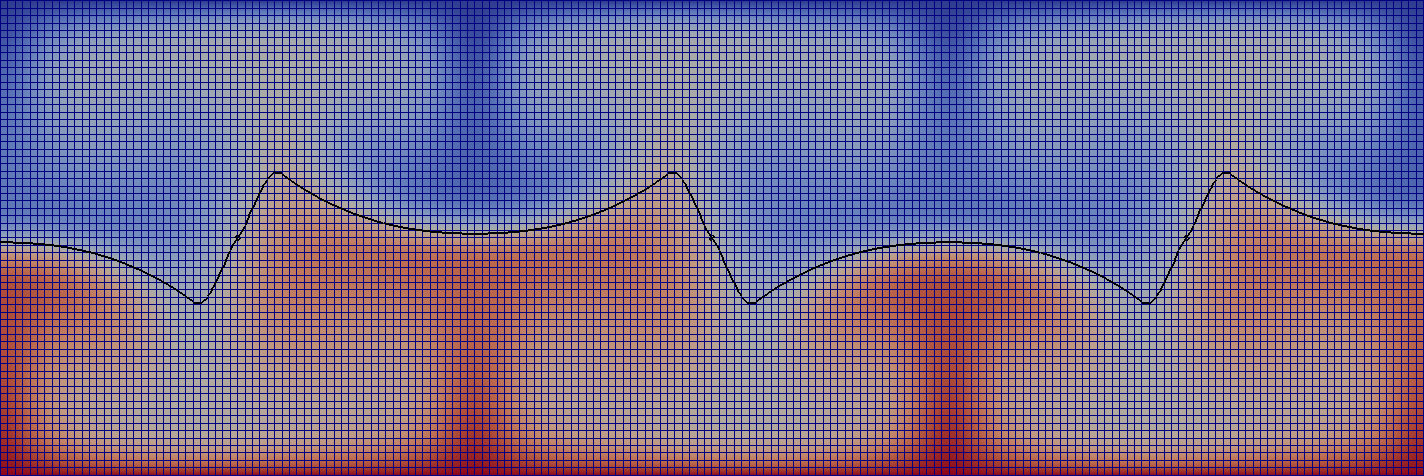}
            \caption{$\mathrm{B} = 0.6$ at $t' = 1.97\cdot10^{-2}$ ON A UNIFORM GRID of $196 \times 64$ cells.}
            \label{FIG:B = 0.6 at t' = 0.0197 on a uniform grid}
        \end{subfigure}

        \begin{subfigure}[b]{1.0\textwidth}
            \includegraphics[width=1.0\textwidth]{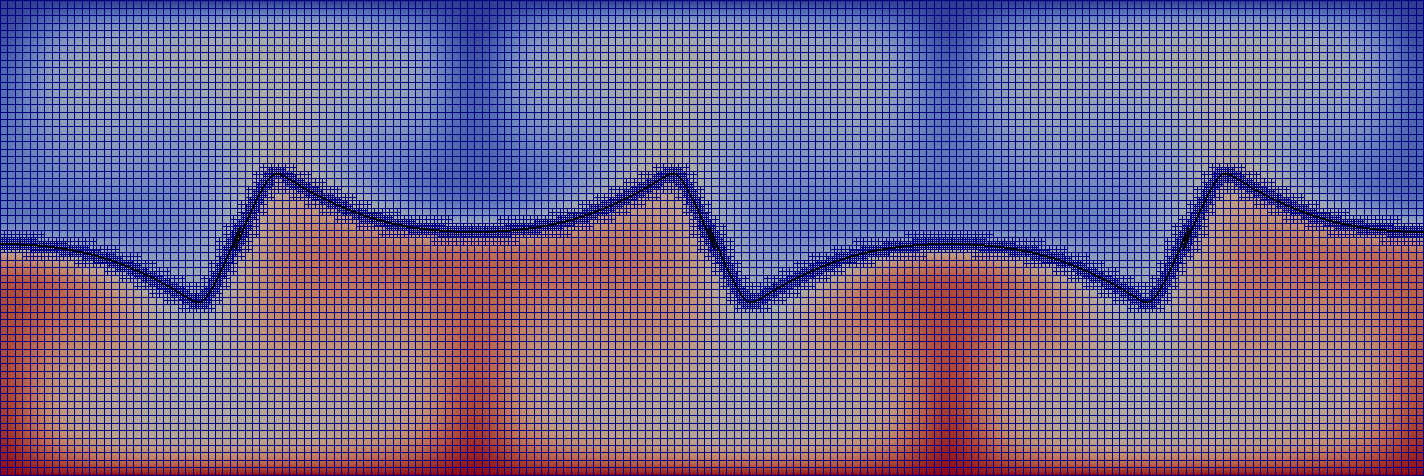}
            \caption{$\mathrm{B} =  0.6$ at $t' = 1.97\cdot10^{-2}$ with two levels of AMR}
            \label{FIG:B = 0.6 at t' = 0.0197 with AMR}
        \end{subfigure}

        \begin{subfigure}[b]{1.0\textwidth}
            \includegraphics[width=1.0\textwidth]{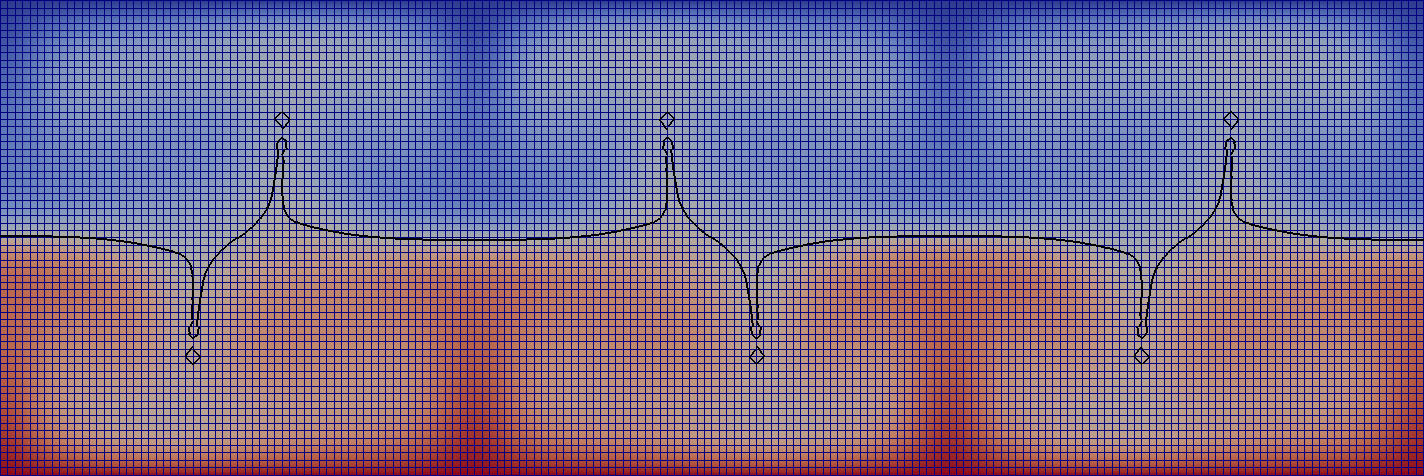}
            \caption{$B = 0.6$ at $t'=2.36\cdot10^{-2}$ ON A UNIFORM GRID of $196 \times 64$ cells.}
            \label{FIG:B = 0.6 at t' = 0.0236 ON A UNIFORM GRID}
        \end{subfigure}

        \begin{subfigure}[b]{1.0\textwidth}
            \includegraphics[width=1.0\textwidth]{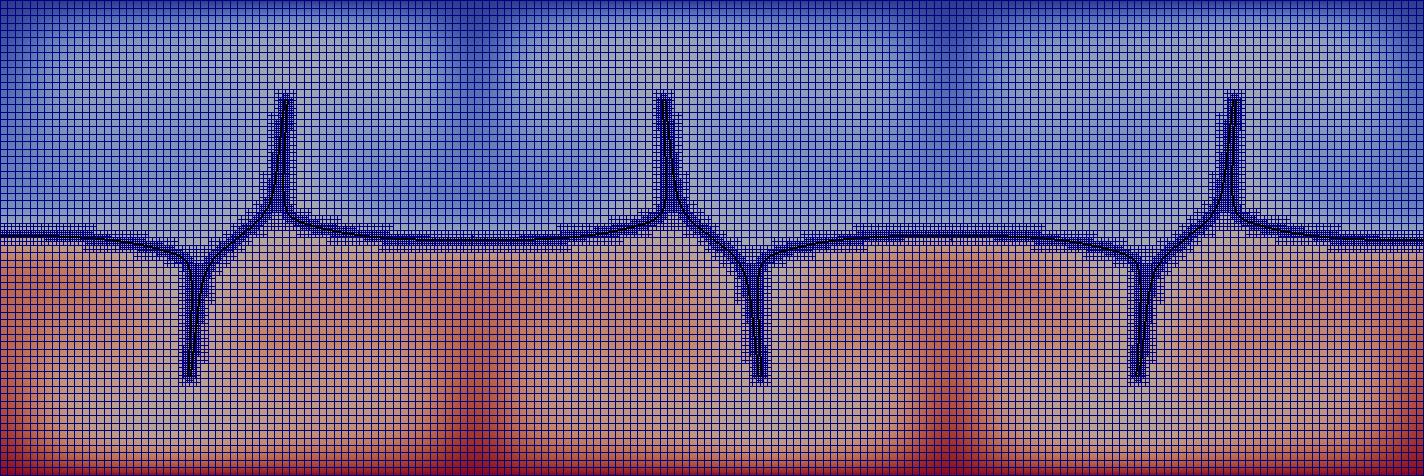}
            \caption{$\mathrm{B} = 0.6$ at $t'=2.36\cdot10^{-2}$ with two levels of AMR}
            \label{FIG:B = 0.6 at t' = 0.0236 WITH AMR}
        \end{subfigure}
         \caption{Computations with $\mathrm{B} = 0.6$ and $\mathrm{Ra} = 10^5$ on 
             an underlying uniform grid of $196 \times 64$ square cells  at
             $t' = 1.97 \cdot 10^{-2}$ and $t' = 2.36 \cdot 10^{-2}$. 
         	 The background color is the temperature, which varies from $T = 0.0$ (dark blue) to $T = 1.0$ (dark red).
                 }
        \label{FIG:B = 0.6}
    \end{center}
\end{figure} 

\begin{figure}
    \begin{center}
        \begin{subfigure}[b]{1.0\textwidth}
            \includegraphics[width=1.0\textwidth]{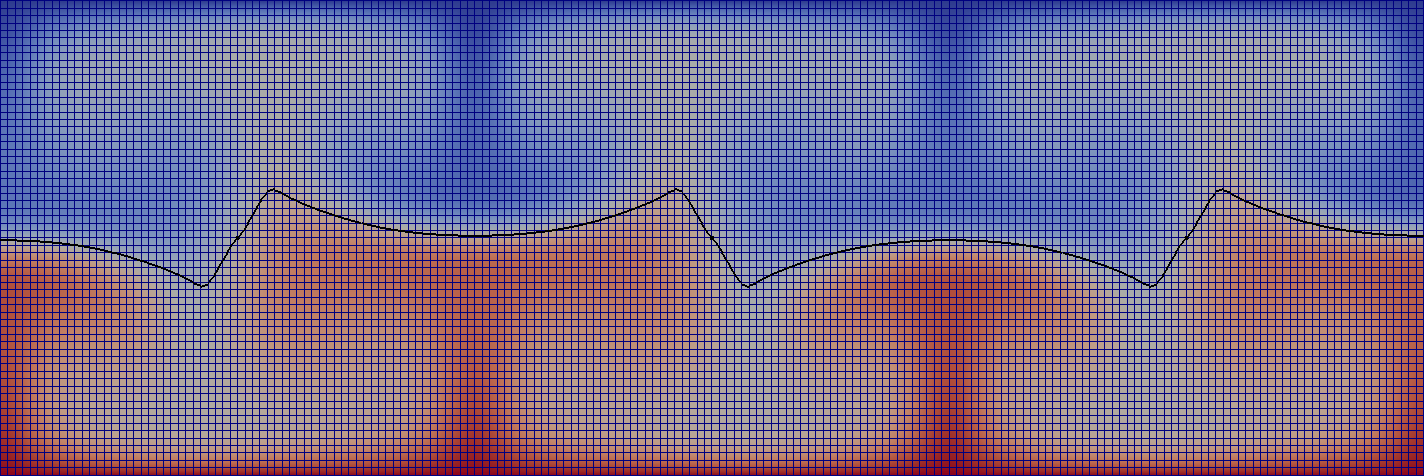}
            \caption{$\mathrm{B} = 0.7$ at $t' = 1.97\cdot10^{-2}$ ON A UNIFORM GRID of $196 \times 64$ cells.}
            \label{FIG:B = 0.7 at t' = 0.0197 on a uniform grid}
        \end{subfigure}

        \begin{subfigure}[b]{1.0\textwidth}
            \includegraphics[width=1.0\textwidth]{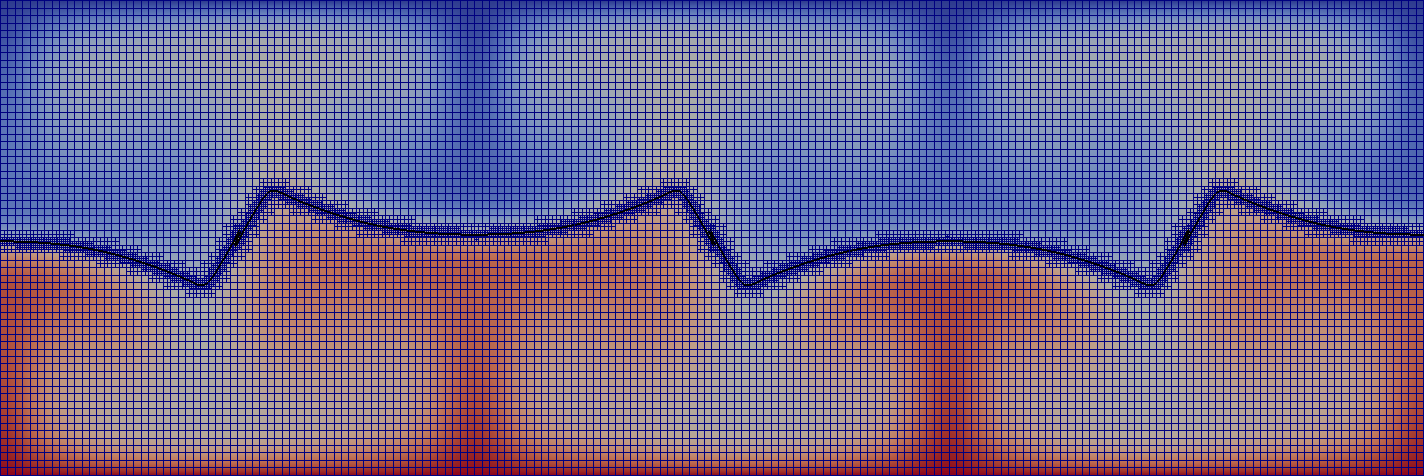}
            \caption{$\mathrm{B} =  0.7$ at $t' = 1.97\cdot10^{-2}$ with two levels of AMR}
            \label{FIG:B = 0.7 at t' = 0.0197 with AMR}
        \end{subfigure}

        \begin{subfigure}[b]{1.0\textwidth}
            \includegraphics[width=1.0\textwidth]{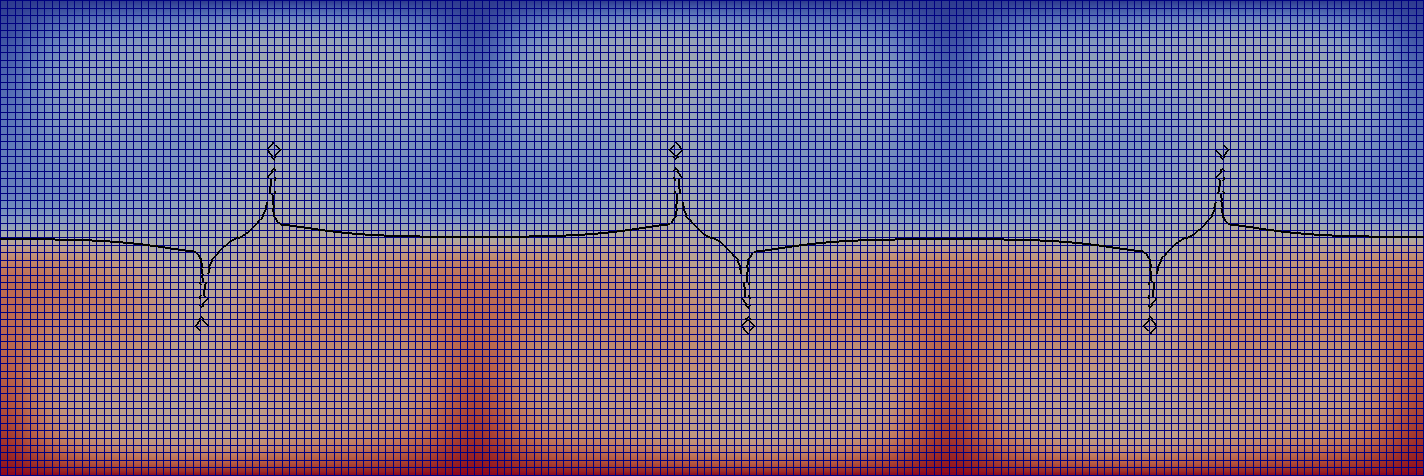}
            \caption{$B = 0.7$ at $t'=2.36\cdot10^{-2}$ ON A UNIFORM GRID of $196 \times 64$ cells.}
            \label{FIG:B = 0.7 at t' = 0.0236 ON A UNIFORM GRID}
        \end{subfigure}

        \begin{subfigure}[b]{1.0\textwidth}
            \includegraphics[width=1.0\textwidth]{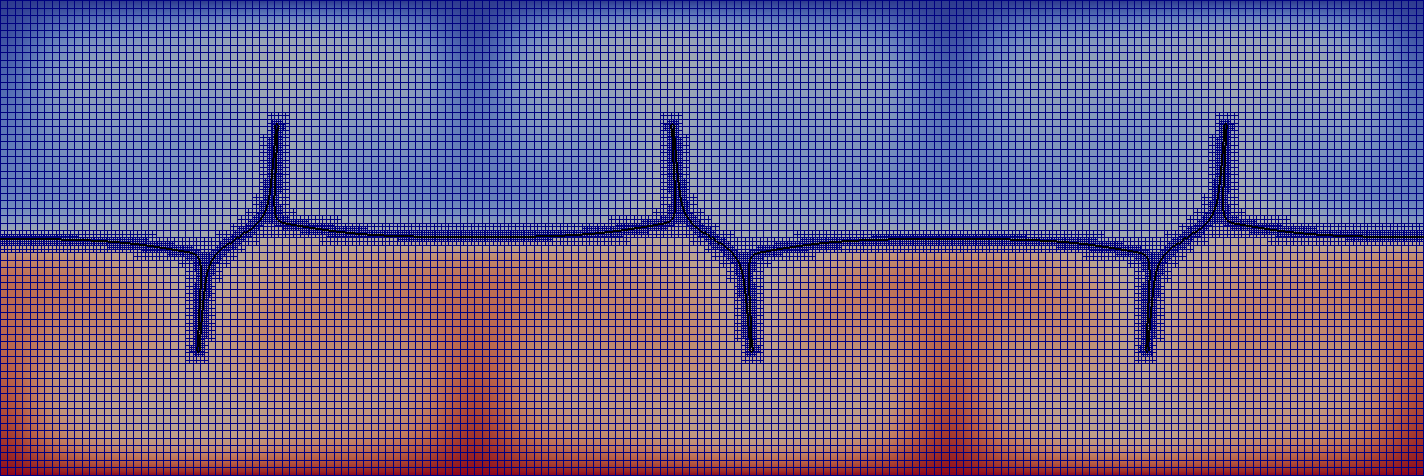}
            \caption{$\mathrm{B} = 0.7$ at $t'=2.36\cdot10^{-2}$ with two levels of AMR}
            \label{FIG:B = 0.7 at t' = 0.0236 WITH AMR}
        \end{subfigure}
         \caption{Computations with $\mathrm{B} = 0.7$ and $\mathrm{Ra} = 10^5$ on 
             an underlying uniform grid of $196 \times 64$ square cells at
             $t' = 1.97 \cdot 10^{-2}$ and $t' = 2.36 \cdot 10^{-2}$. 
         	 The background color is the temperature, which varies from $T = 0.0$ (dark blue) to $T = 1.0$ (dark red).
                 }
        \label{FIG:B = 0.7}
    \end{center}
\end{figure} 

\begin{figure}
    \begin{center}
        \begin{subfigure}[b]{1.0\textwidth}
            \includegraphics[width=1.0\textwidth]{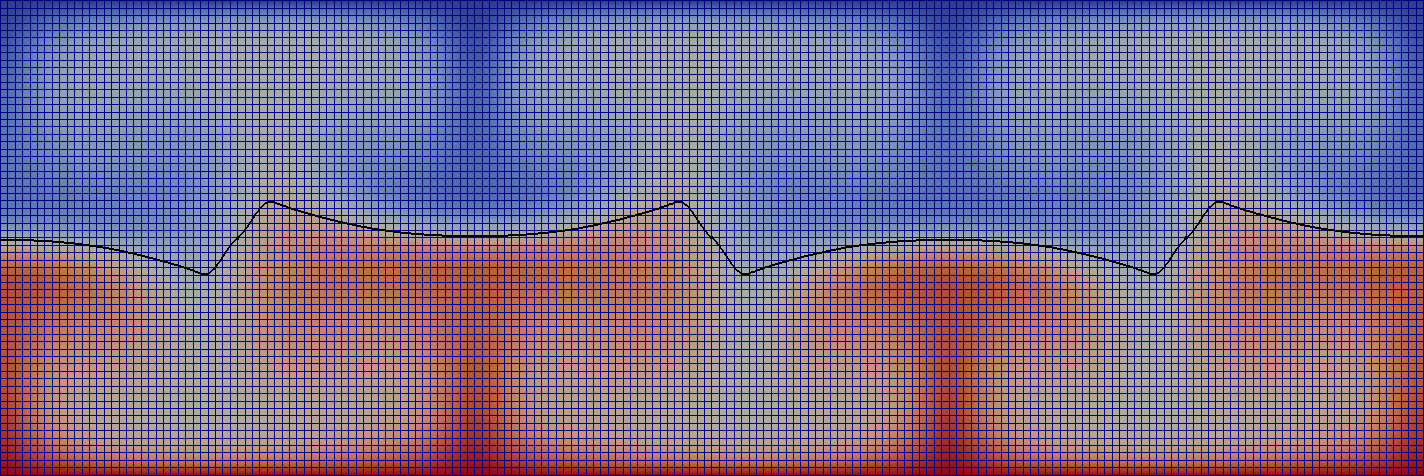}
            \caption{$\mathrm{B} = 0.8$ at $t' = 1.97\cdot10^{-2}$ ON A UNIFORM GRID of $196 \times 64$ cells.}
            \label{FIG:B = 0.8 at t' = 0.0197 on a uniform grid}
        \end{subfigure}

        \begin{subfigure}[b]{1.0\textwidth}
            \includegraphics[width=1.0\textwidth]{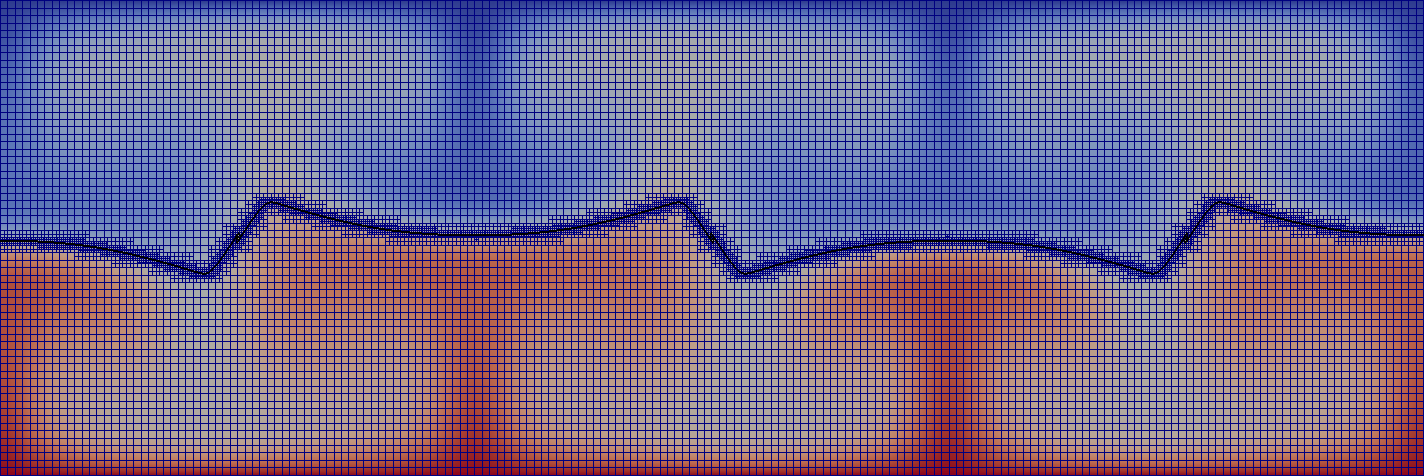}
            \caption{$\mathrm{B} =  0.8$ at $t' = 1.97\cdot10^{-2}$ with two levels of AMR}
            \label{FIG:B = 0.8 at t' = 0.0197 with AMR}
        \end{subfigure}

        \begin{subfigure}[b]{1.0\textwidth}
            \includegraphics[width=1.0\textwidth]{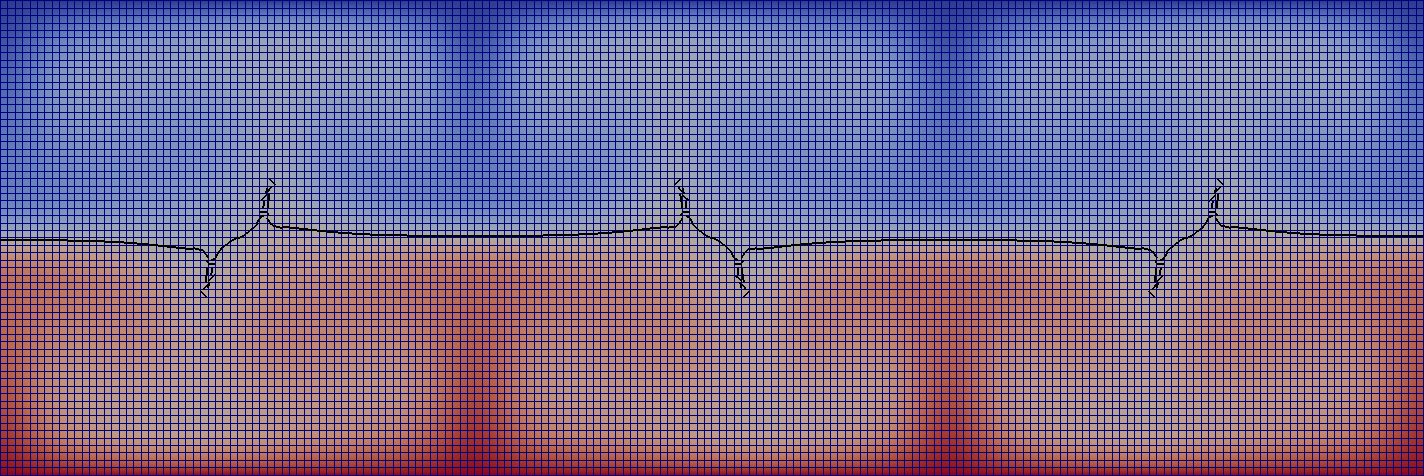}
            \caption{$B = 0.8$ at $t'=2.36\cdot10^{-2}$ ON A UNIFORM GRID of $196 \times 64$ cells.}
            \label{FIG:B = 0.8 at t' = 0.0236 ON A UNIFORM GRID}
        \end{subfigure}

        \begin{subfigure}[b]{1.0\textwidth}
            \includegraphics[width=1.0\textwidth]{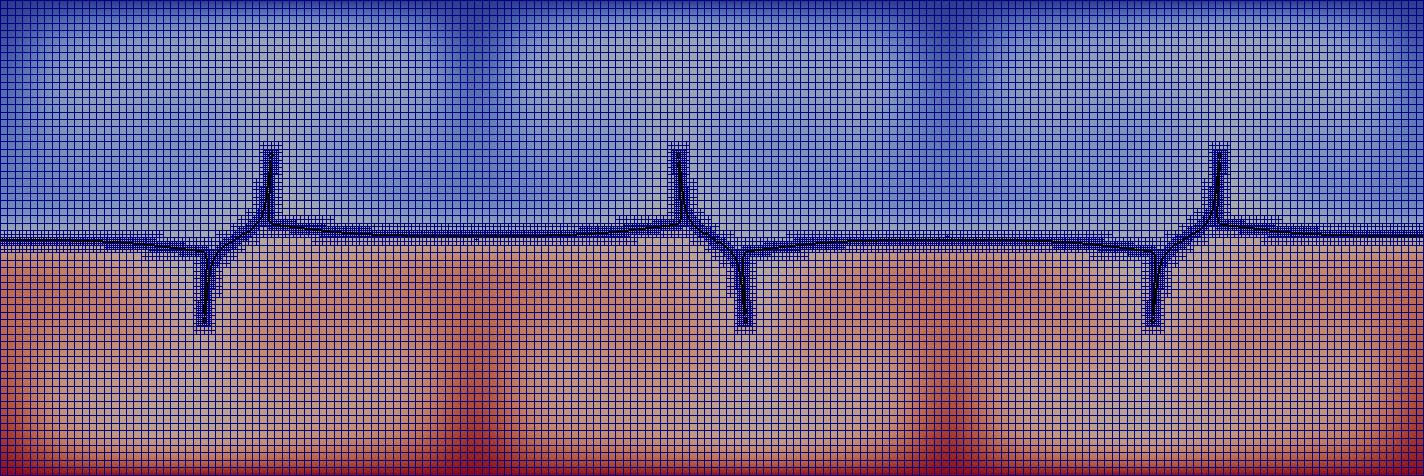}
            \caption{$\mathrm{B} = 0.8$ at $t'=2.36\cdot10^{-2}$ with two levels of AMR}
            \label{FIG:B = 0.8 at t' = 0.0236 WITH AMR}
        \end{subfigure}
         \caption{Computations with $\mathrm{B} = 0.8$ and $\mathrm{Ra} = 10^5$ on 
             an underlying uniform grid of $196 \times 64$ square cells at
             $t' = 1.97 \cdot 10^{-2}$ and $t' = 2.36 \cdot 10^{-2}$. 
         	 The background color is the temperature, which varies from $T = 0.0$ (dark blue) to $T = 1.0$ (dark red).
                 }
        \label{FIG:B STUDY ALL B = 0 8}
    \end{center}
\end{figure} 

\begin{figure}
    \begin{center}
        \begin{subfigure}[b]{1.0\textwidth}
            \includegraphics[width=1.0\textwidth]{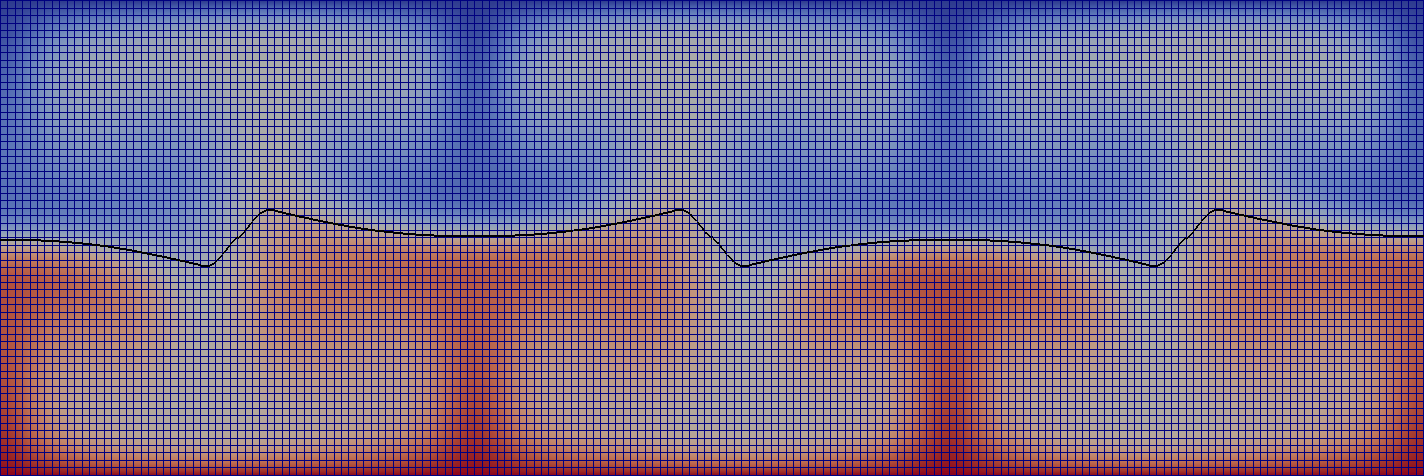}
            \caption{$\mathrm{B} = 0.9$ at $t' = 1.97\cdot10^{-2}$ ON A UNIFORM GRID of $196 \times 64$ cells.}
            \label{FIG:B = 0.9 at t' = 0.0197 on a uniform grid}
        \end{subfigure}

        \begin{subfigure}[b]{1.0\textwidth}
            \includegraphics[width=1.0\textwidth]{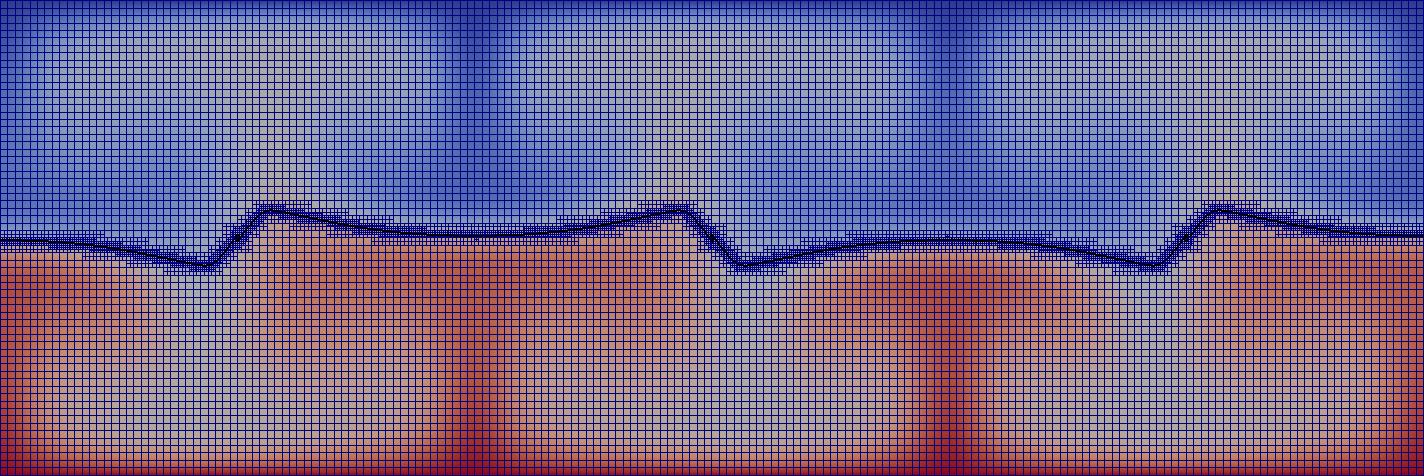}
            \caption{$\mathrm{B} =  0.9$ at $t' = 1.97\cdot10^{-2}$ with two levels of AMR}
            \label{FIG:B = 0.9 at t' = 0.0197 with AMR}
        \end{subfigure}

        \begin{subfigure}[b]{1.0\textwidth}
            \includegraphics[width=1.0\textwidth]{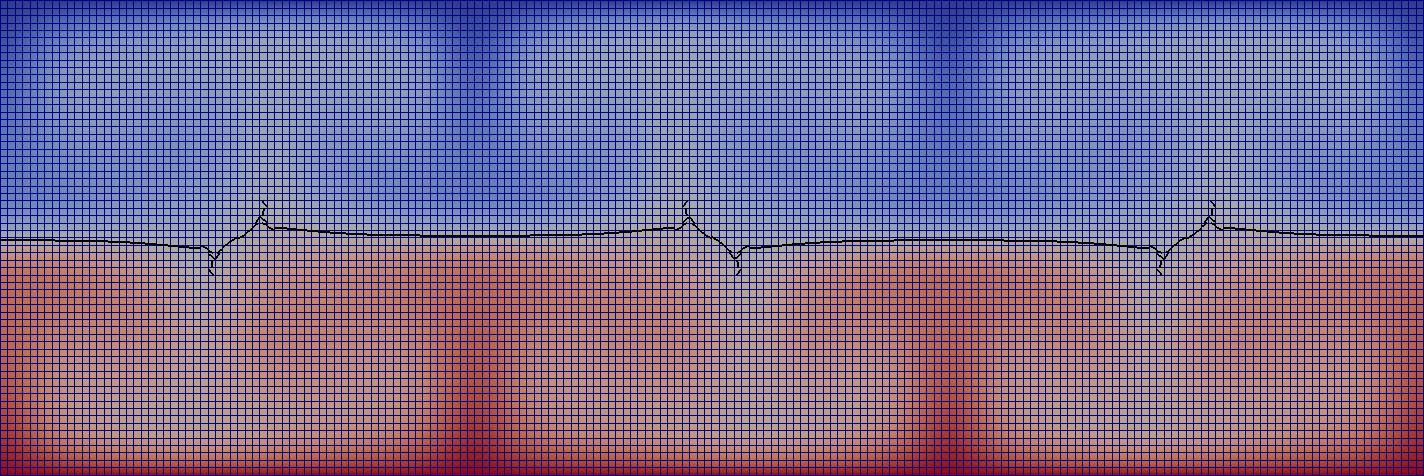}
            \caption{$B = 0.9$ at $t'=2.36\cdot10^{-2}$ ON A UNIFORM GRID of $196 \times 64$ cells.}
            \label{FIG:B = 0.9 at t' = 0.0236 ON A UNIFORM GRID}
        \end{subfigure}

        \begin{subfigure}[b]{1.0\textwidth}
            \includegraphics[width=1.0\textwidth]{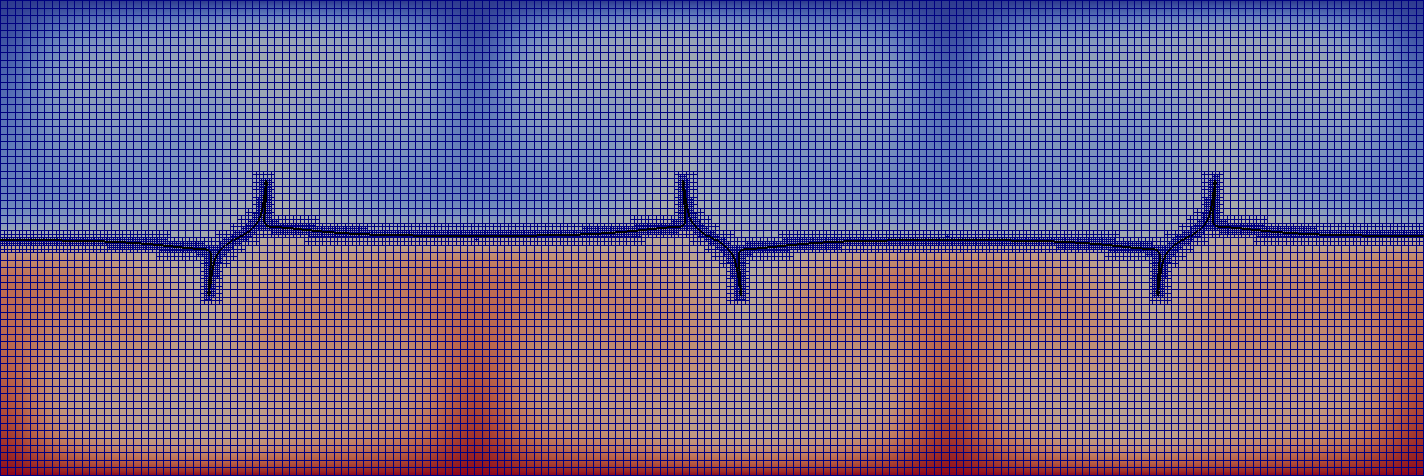}
            \caption{$\mathrm{B} = 0.9$ at $t'=2.36\cdot10^{-2}$ with two levels of AMR}
            \label{FIG:B = 0.9 at t' = 0.0236 WITH AMR}
        \end{subfigure}
         \caption{Computations with $\mathrm{B} = 0.9$ and $\mathrm{Ra} = 10^5$ on 
             an underlying uniform grid of $196 \times 64$ square cells at
             $t' = 1.97 \cdot 10^{-2}$ and $t' = 2.36 \cdot 10^{-2}$. 
         	 The background color is the temperature, which varies from $T = 0.0$ (dark blue) to $T = 1.0$ (dark red).
                 }
        \label{FIG:B STUDY ALL B = 0 9}
    \end{center}
\end{figure} 

\begin{figure}
    \begin{center}
        \begin{subfigure}[b]{1.0\textwidth}
            \includegraphics[width=1.0\textwidth]{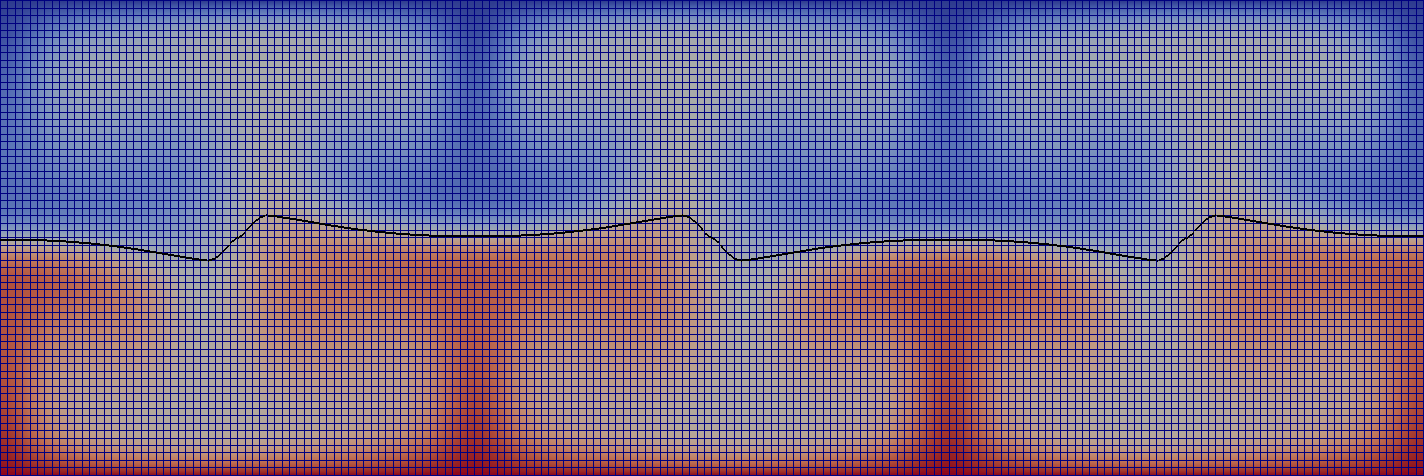}
            \caption{$\mathrm{B} = 1.0$ at $t' = 1.97\cdot10^{-2}$ ON A UNIFORM GRID of $196 \times 64$ cells.}
            \label{FIG:B = 1.0 at t' = 0.0197 on a uniform grid}
        \end{subfigure}

        \begin{subfigure}[b]{1.0\textwidth}
            \includegraphics[width=1.0\textwidth]{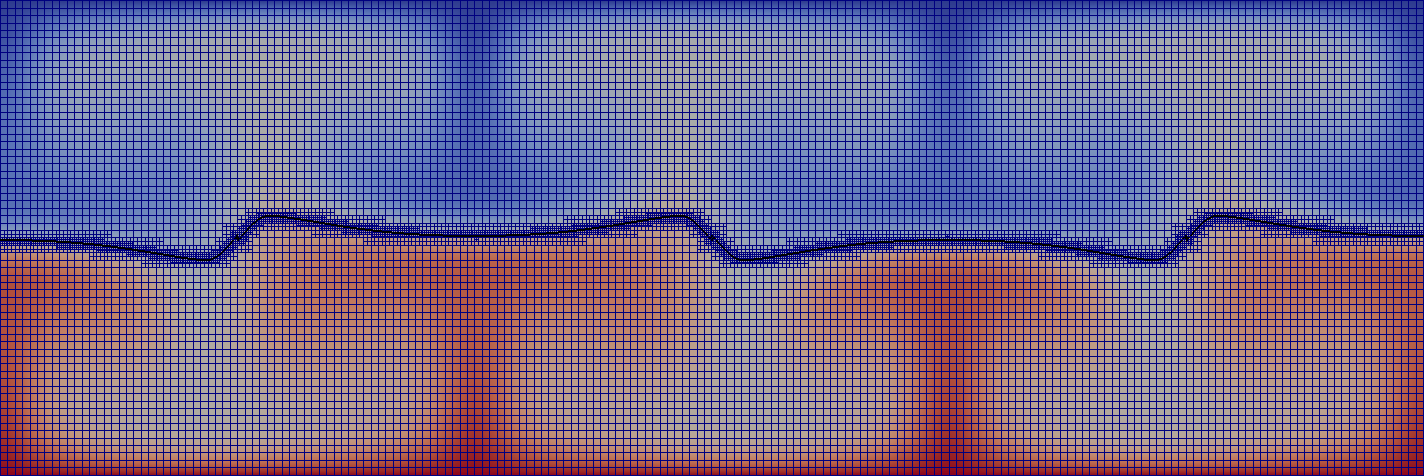}
            \caption{$\mathrm{B} =  1.0$ at $t' = 1.97\cdot10^{-2}$ with two levels of AMR}
            \label{FIG:B = 1.0 at t' = 0.0197 with AMR}
        \end{subfigure}

        \begin{subfigure}[b]{1.0\textwidth}
            \includegraphics[width=1.0\textwidth]{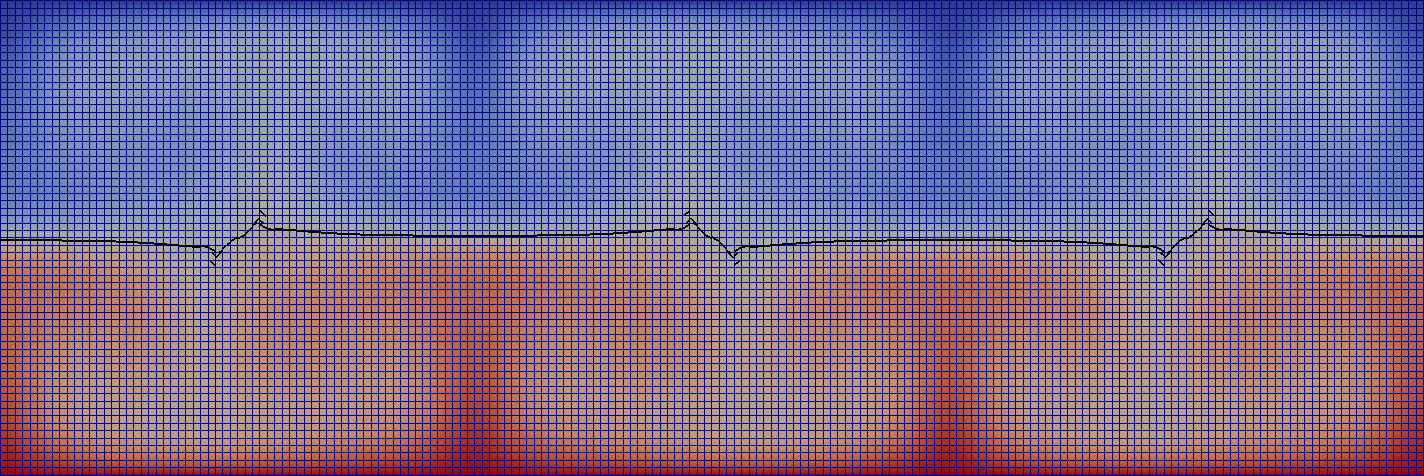}
            \caption{$B = 1.0$ at $t'=2.36\cdot10^{-2}$ ON A UNIFORM GRID of $196 \times 64$ cells.}
            \label{FIG:B = 1.0 at t= 3.0 ON A UNIFORM GRID}
        \end{subfigure}

        \begin{subfigure}[b]{1.0\textwidth}
            \includegraphics[width=1.0\textwidth]{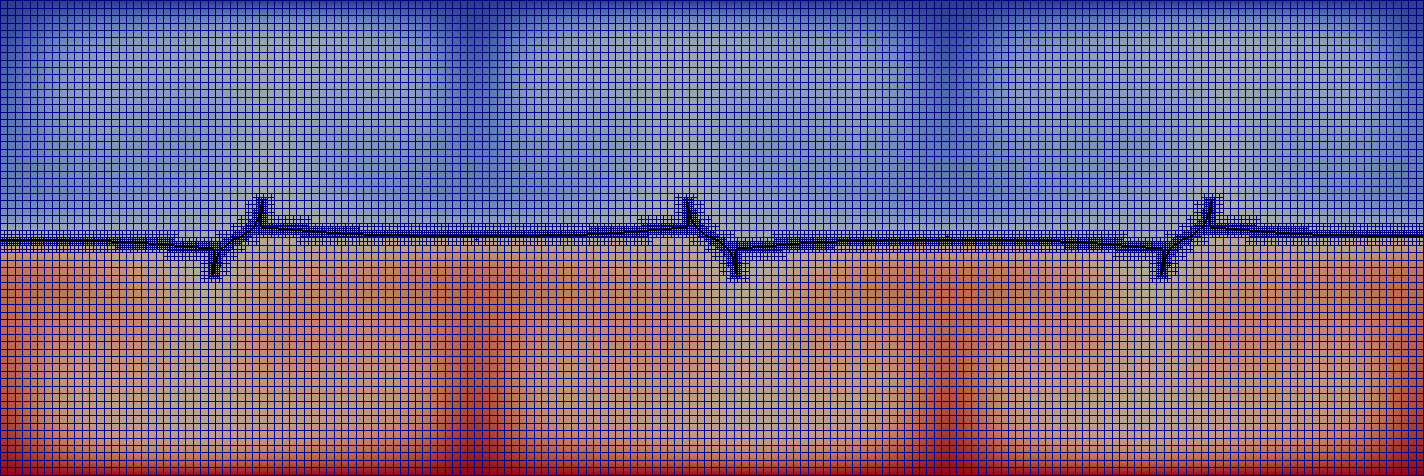}
            \caption{$\mathrm{B} = 1.0$ at $t'=2.36\cdot10^{-2}$ with two levels of AMR}
            \label{FIG:B = 1.0 at t' = 0.0236 WITH AMR}
        \end{subfigure}
         \caption{Computations with $\mathrm{B} = 1.0$ and $\mathrm{Ra} = 10^5$ on 
             an underlying uniform grid of $196 \times 64$ square cells at
             $t' = 1.97 \cdot 10^{-2}$ and $t' = 2.36 \cdot 10^{-2}$. 
         	The background color is the temperature, which varies from $T = 0.0$ (dark blue) to $T = 1.0$ (dark red).
                 }
        \label{FIG:B STUDY ALL B = 1.0}
    \end{center}
\end{figure} 

\begin{figure}
    \begin{center}
        \begin{subfigure}[b]{1.0\textwidth}
            \includegraphics[width=1.0\textwidth]{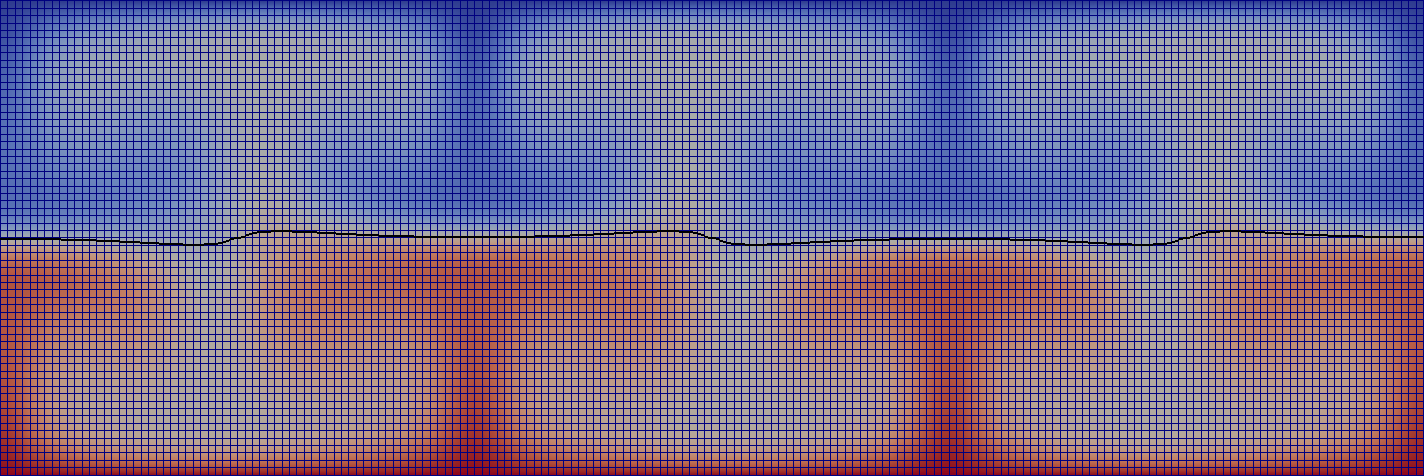}
            \caption{$\mathrm{B} = 2.0$ at $t' = 1.97\cdot10^{-2}$ ON A UNIFORM GRID of $196 \times 64$ cells.}
            \label{FIG:B = 2.0 at t' = 0.0197 on a uniform grid}
        \end{subfigure}

        \begin{subfigure}[b]{1.0\textwidth}
            \includegraphics[width=1.0\textwidth]{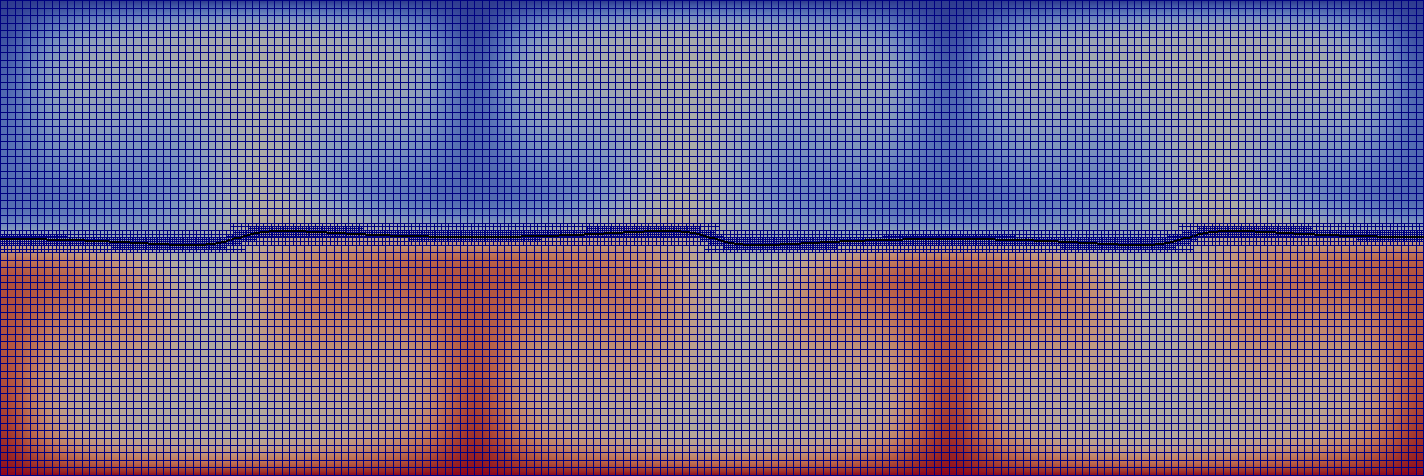}
            \caption{$\mathrm{B} =  2.0$ at $t' = 1.97\cdot10^{-2}$ with two levels of AMR}
            \label{FIG:B = 2.0 at t' = 0.0197 with AMR}
        \end{subfigure}

        \begin{subfigure}[b]{1.0\textwidth}
            \includegraphics[width=1.0\textwidth]{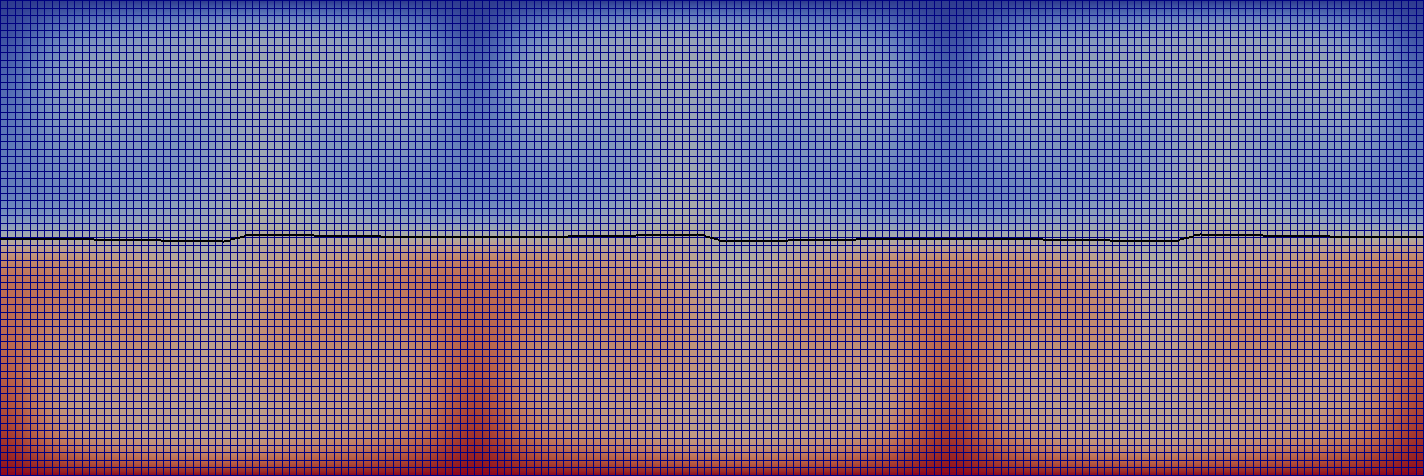}
            \caption{$\mathrm{B} = 2.0$ at $t' = 2.36 \cdot 10^{-2}$ ON A UNIFORM GRID of
            	$196 \times 64$ square cells.}
            \label{FIG:B = 2.0 at t'=2.36 ON A UNIFORM GRID}
        \end{subfigure}

        \begin{subfigure}[b]{1.0\textwidth}
            \includegraphics[width=1.0\textwidth]{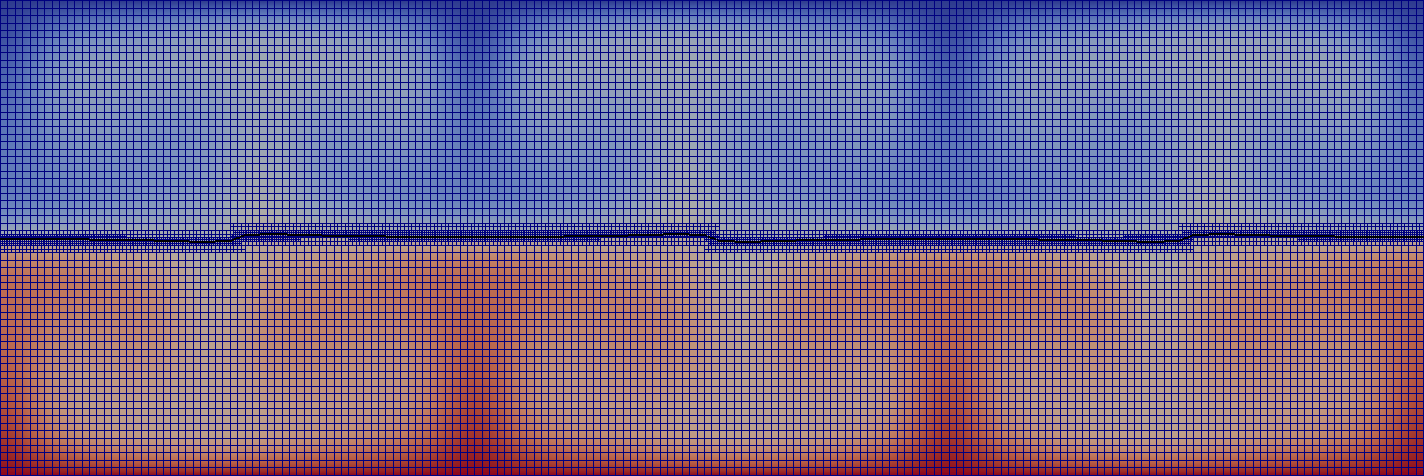}
            \caption{$\mathrm{B} = 2.0$ at $t' = 2.36 \cdot 10^{-2}$ on a uniform $196 \times 64$ grid and two levels of AMR}
            \label{FIG:B = 2.0 at t' = 0.0236 WITH AMR}
        \end{subfigure}
         \caption{Computations with $\mathrm{B} = 2.0$ and $\mathrm{Ra} = 10^5$ on an 
         	underlying uniform grid of $196 \times 64$ square cells at
         	$t' = 1.97 \cdot 10^{-2}$ and $t' = 2.36 \cdot 10^{-2}$. 
         	The background color is the temperature, which varies from $T = 0.0$ (dark blue) to $T = 1.0$ (dark red).
                 }
        \label{FIG:B = 2.0}
    \end{center}
\end{figure}


\section{Discussion}
\label{Section:Discussion}

In Section~\ref{SUBSECTION:INTERFACE TRACKING BENCHMARK PROBLEMS} we demonstrated that our implementation of the VOF method in ASPECT is second-order accurate on smooth flows in the norm ddefined in~\eqref{DEF:VOLUME FRACTION ERROR}. 
In Section~\ref{SUBSECTION:MANTLE CONVECTION BENCHMARK PROBLEMS} we demonstrated that the method correctly reproduces two benchmarks from the computational mantle convection literature.
We now present a detailed discussion of the results of our computations of thermochemical convection in density stratified flow shown in
Section~\ref{SUBSECTION:THERMOCHEMICAL CONVECTION IN A DENSITY STRATIFIED FLUID}.
This model problem is designed to study the basic physics underlying the formation of thermal plumes that form at LLSVPs, entrain some of the material in the LLSVP, and bring it to the Earth's surface.
It is also a two dimensional analog of the experimental results of Davaille~\cite{AD:1999} and Le Bars and Davaille~\cite{MLB-AD:2004,MLB-AD:2005}.

\subsection{Computations of Thermochemical Convection in a Density Stratified Fluid}

Examining the results in
Section~\ref{SUBSECTION:THERMOCHEMICAL CONVECTION IN A DENSITY STRATIFIED FLUID} of our computations of thermochemical convection in a density stratified fluid  for values of the nondimensional buoyancy parameter
$\mathrm{B} = 0.0, \, 0.1, \, 0.2, \ldots, 1.0$ and $\mathrm{B} = 2.0$ at Rayleigh number $\mathrm{Ra} = 10^5$, we note a fundamental change in the dynamics and structure of the flow field as $\mathrm{B}$ increases from $\mathrm{B} = 0.0$ to $\mathrm{B} = 2.0$.
First, considering only the extreme values $\mathrm{B}=0.0$ and $\mathrm{B}=2.0$, we observe the following difference in the qualitative behavior of the interface. 
For $\mathrm{B}=0.0$ (Figure~\ref{FIG:B = 0.0}), which is the classic Rayleigh-B\'{e}nard problem in which there is no difference in the densities of the two fluids (i.e., $\Delta \rho = 0$), the height of the convection cells is equal to the height of the domain $\Omega$ and we observe the steady cellular convection structure with three $1 \times 1$ counter rotating cells as predicted by the analysis in Section~6.21 of~\cite{DLT-GS:2014}.
That the flow is steady, (i.e., independent of time) in Figure~\ref{FIG:B = 0.0} is apparent after comparing the temperature fields at  $t' = 1.97\cdot10^{-2}$ and
$t' = 2.36 \cdot 10^{-2}$.


Note that for $\mathrm{B} = 0.0$ each of the three  $1 \times 1$ convection cells overturn at the same fixed rate.
On the other hand, for $\mathrm{B} = 2.0$ the magnitude of $\Delta \rho$ prevents the denser fluid from reaching the top of the domain and producing overturns, and hence convection cells, on the scale of the height of the domain. 
Rather, the structure of the flow shown in Figure~\ref{FIG:B = 2.0} consists of six (roughly) square counter rotating $\frac{1}{2} \times \frac{1}{2}$ cells below
$y = 0.5$ and a similar structure above $y \, = \, 0.5$. 
Thus, for $\mathrm{B} = 2.0$ we observe a permanently stratified convection structure.
Furthermore, from $\mathrm{B} = 0.7$ in Figure~\ref{FIG:B = 0.7} and, \textit{perhaps}, from $\mathrm{B} = 0.4$ in Figure~\ref{FIG:B = 0.4} or  $\mathrm{B} = 0.5$ in Figure~\ref{FIG:B = 0.5}, on; i.e., as
$\mathrm{B} \, \to \, 2.0$ from below with $\mathrm{B} > \mathrm{B}_{c}$ where
$0.3 < \mathrm{B}_{c} \le 0.7$, it appears that \textit{at the times shown} the flow is tending \textit{continuously} toward the stratified convection pattern shown in
Figure~\ref{FIG:B = 2.0}.

The features at either end of the interval $\mathrm{B} = [0.0 \, , 2.0]$ are consistent with the diagrams - obtained from experiments - on the left and right of Figure~1 in~\cite{MLB-AD:2005}, although in the diagram on the right the authors have only drawn three cells above and three cells below the centerline and, in both drawings, the cells appear to be more rectangular than square in shape. 
We assume that these diagrams are simply rough sketches of the dynamics of what the authors of~\cite{MLB-AD:2005} refer to as ``Whole Layer'' (left) and ``Stratified'' (right) convection.
Perhaps, also, these diagrams are for different values of the other two nondimensional parameters the authors varied in the work described in the sequence of papers~\cite{MLB-AD:2004,MLB-AD:2005} and \cite{AD:1999}; namely, the ratio $a$ of the height of the lower layer to the height of the entire domain and the ratio $\gamma$ of viscosity of the lower layer to that of the upper layer.
In the work we present in 
Section~\ref{SUBSECTION:THERMOCHEMICAL CONVECTION IN A DENSITY STRATIFIED FLUID} we did not vary these other two parameters; they were held fixed at $a = 0.5$ and
$\gamma = 1.0$.
In short, we conclude that our computational results correctly correspond \textit{qualitatively} to what the authors of~\cite{MLB-AD:2005} observe in their experiments when the nondimensional parameters $a$ and $\gamma$ are held fixed at
$a = 0.5$ and $\gamma = 1.0$.
Finally, note the similarity of the two counter rotating convection cells on the right in 
Figures~\ref{FIG:B = 0.0 at t' = 0.0236 ON A UNIFORM GRID}--\ref{FIG:B = 0.0 at t' = 0.0236 WITH AMR}
and~\ref{FIG:B = 0.1 at t' = 0.0236 ON A UNIFORM GRID}--\ref{FIG:B = 0.1 at t' = 0.0236 WITH AMR} to the structure of the flow in Figure~4(a) of~\cite{MLB-AD:2004}.

It is possible to obtain additional insight into the structure and dynamics of the flow for various values of $\mathrm{B}$ from the results shown in
 Figures~\ref{FIG:B = 0.0}--\ref{FIG:B = 2.0}.
As $\mathrm{B}$ increases from $0.0$ to $0.1$, $0.2$, and $0.3$ in
Figures~\ref{FIG:B = 0.0}--\ref{FIG:B = 0.3} we observe
that the rate of overturn decreases, until for $\mathrm{B}=0.3$ the denser material has just reached the top of the domain at $t' = 2.36 \cdot 10^{-2}$
(Figures~\ref{FIG:B = 0.3 at t' = 2.36 ON A UNIFORM GRID} and 
\ref{FIG:B = 0.3 at t' = 2.36 WITH AMR}), whereas for smaller values of $\mathrm{B}$ the overturn has passed beyond the top of the domain by $t' = 2.36 \cdot 10^{-2}$. 
For $\mathrm{B} = 0.4$ we can see from 
Figures~\ref{FIG:B = 0.4 at t' = 2.36 ON A UNIFORM GRID}
and~\ref{FIG:B = 0.4 at t' = 2.36 WITH AMR}, that the fluid does not reach a full overturn by $t' = 2.36 \cdot 10^{-2}$ suggesting that there may be a transition between the qualitative dynamics of the flow at some $\mathrm{B}_{c}$ in the range
$0.3 \le \mathrm{B} \le 0.4$.
In~\cite{MLB-AD:2005} the authors find $\mathrm{B}_{c} = 0.302$ when the viscosity ratio is $\gamma = 6.7$.

For $0.5 \le \mathrm{B} \le 1.0$ in
Figures~\ref{FIG:B = 0.5}--\ref{FIG:B STUDY ALL B = 1.0} the general interface structures are similar, although with smaller volumes for the ``pinched'' regions that are  produced during the transition from ``Whole Layer'' convection to ``Stratified'' flow.
As shown in Figure~\ref{FIG:B = 2.0}, for $\mathrm{B} = 2.0$, the stratification is sufficiently strong that the pinched structures do not form, although a standing wave does form as a slight perturbation from the initial location of the interface at $y = \frac{1}{2}$ with boundaries at $x \approxeq 0.5, \, 1.5, \, 2.5$.

\subsubsection{A qualitative comparison to the experiments Davaille and Le Bars}

\begin{figure}[t!]
	\begin{center}
		\includegraphics[width=\textwidth]{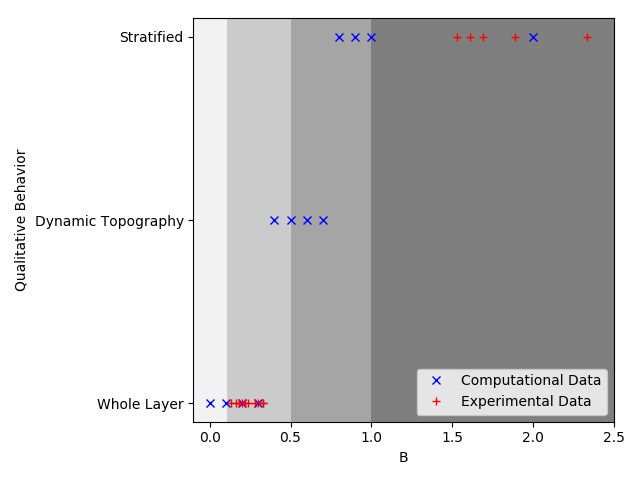}
		\caption{A \textit{qualitative} comparison of the computations presented in this 
			paper to the experimental results of Davaille~\protect\cite{AD:1999} and Le Bars \& Davaille~\protect\cite{MLB-AD:2004,MLB-AD:2005}.
			The grayscale regions correspond to boundaries of the qualitative regions shown in  in Figure~3 of~\protect\cite{AD:1999} and Figure~2 of~\protect\cite{MLB-AD:2005} for $a=0.5$.
			The experimental data is from Table~3~of~\protect\cite{AD:1999} and Table~3~of~\protect\cite{MLB-AD:2004} with $a=0.5$ as is the case for all of the computations in this article.
            The terms ``Stratified'', ``Dynamic Topography'', and ``Whole Layer'' used to describe the qualitative state of the flow are the same as those used by the authors of~\protect\cite{MLB-AD:2004,MLB-AD:2005} and~\protect\cite{AD:1999}.
		}
		\label{FIG:COMPUTATION VERSUS EXPERIMENT}
	\end{center}
\end{figure}

In this section we briefly make some additional \textit{qualitative} comparisons of our computational results to the experimental results of Davaille~\cite{AD:1999} and Le Bars \& Davaille~\cite{MLB-AD:2004,MLB-AD:2005}. 
Before doing so however, it is first necessary to make several caveats concerning this comparison.
First, as we mentioned above, in the experiments the authors varied two additional nondimensional parameters; namely, (1) the ratio $a$ of the height of the lower layer to the height of the entire domain and (2) the ratio $\gamma$ of the viscosity of the fluid that initially occupies the lower layer to the viscosity that initially occupies the upper layer.
In our computations, shown in
 Section~\ref{SUBSECTION:THERMOCHEMICAL CONVECTION IN A DENSITY STRATIFIED FLUID}, we kept these parameters fixed at $a = 0.5$ and $\gamma = 1.0$.
Second, in the experiments the two fluids are \textit{miscible}, whereas in our computations the two fluids are \textit{immiscible}. 
In both cases there is no surface tension at the boundary between the two fluids.

The general transition between one type of structure and another (e.g., ``Whole Layer'' convection to ``Stratified Convection'') is similar to that found in the experiments shown in~\cite{MLB-AD:2005}, although the precise location of the transition may differ.
A rough comparison is show in Figure~\ref{FIG:COMPUTATION VERSUS EXPERIMENT}. The different grayscale backgrounds in Figure~\ref{FIG:COMPUTATION VERSUS EXPERIMENT} correspond to the grayscale regions in Figure~3 of~\cite{AD:1999} and Figure~2 of~\cite{MLB-AD:2005} for $a=0.5$.
In the results presented in this paper we do not continue the computation for a
sufficiently long times to confirm that in the $0.3 \le \mathrm{B} \le 0.5$ regime the flow oscillates before beginning an overturn.
However, the observed behavior does produce structures that match those described in~\cite{MLB-AD:2005} for the length of time for which we do have computational results.
This difference may be in part due to the fact that in~\cite{MLB-AD:2005} the two fluids also vary in viscosity ratio $\gamma$, and Rayleigh number $\mathrm{Ra}$.

\subsubsection{Numerical artifacts that occur when the interface is underresolved}
Since the VOF method maintains a sharp interface between the two compositional fields, it is able to capture features that are approximately on the order of the grid scale $h$. 
However, in cases where the structures formed by the interface become sufficiently small, for example, a thin column of fluid of width $2 \, h$, the interface reconstruction algorithm might produce numerical artifacts that are ``characteristic'' of the combination of the particular reconstruction algorithm and advection algorithm one chooses to use in the VOF method.\footnote{It is important to recognize that this is not a failing of the VOF method in general or of the specific interface reconstruction and advection algorithms we have chosen for our work here, since whenever a computation is underresolved, \textit{all} numerical methods will exhibit some sort of numerical artifact or artifacts that are ``characteristic'' of that particular method.}
Here we briefly examine of the nature of one particular numerical artifact that appears frequently in
Section~\ref{SUBSECTION:THERMOCHEMICAL CONVECTION IN A DENSITY STRATIFIED FLUID}.

The most common numerical artifact in the computational results shown in Section~\ref{SUBSECTION:THERMOCHEMICAL CONVECTION IN A DENSITY STRATIFIED FLUID} is the tendency for the reconstructed interface to form `droplets' that are diamond shaped and generally occupy a square of four cells, each edge having two square cells of side $h$.
For example, droplets such as these appear in
Figure~\ref{FIG:B = 0.1 at t' = 0.0197 ON A UNIFORM GRID}.
In the computations shown in
Section~\ref{SUBSECTION:THERMOCHEMICAL CONVECTION IN A DENSITY STRATIFIED FLUID} these droplets typically resolve into a thin vertical column of fluid of approximately
$2 h-4 h$ in width with a length that is nearly the entire height of the computational domain.
For example Figures~\ref{FIG:B = 0.1 at t' = 0.0197 ON A UNIFORM GRID} and 
Figure~\ref{FIG:B = 0.1 at t' = 0.0197 with AMR}, in which  the more refined computation
in Figure~\ref{FIG:B = 0.1 at t' = 0.0197 with AMR} appears to be sufficiently well-resolved to draw the conclusion that a thin column of fluid is forming in the locations where in Figure~\ref{FIG:B = 0.1 at t' = 0.0197 ON A UNIFORM GRID} there are only a few droplets and no real indication of what the flow ``should'' look like.
Or the droplets may resolve into a thin finger that is shorter than the height of the computational domain such as in
Figures~\ref{FIG:B = 0.6 at t' = 0.0236 ON A UNIFORM GRID}
~\ref{FIG:B = 0.6 at t' = 0.0236 WITH AMR}.
~\ref{FIG:B = 0.7 at t' = 0.0236 ON A UNIFORM GRID}
and~\ref{FIG:B = 0.7 at t' = 0.0236 WITH AMR}.

We note that if a feature of the interface is underresolved, it can help the user determine if additional refinement is required. 
In some instances, perhaps after making a second, more refined computation, it will be clear that additional refinement is necessary, sometimes even more refined than the second computation was.
For example see Figures~\ref{FIG:B = 0.1 at t' = 0.0236 ON A UNIFORM GRID} 
and~\ref{FIG:B = 0.1 at t' = 0.0236 WITH AMR}, neither of which appear sufficiently well resolved to accept the computation in
Figure~\ref{FIG:B = 0.1 at t' = 0.0236 WITH AMR} as well resolved enough to determine the true nature of the flow.
On the other hand, there are instances when the numerical artifact is sufficiently small so as not to affect the dynamics of the interface that are of interest and additional resolution might not be required.
For example, depending on the user and the underlying scientific application, this might be the case for Figures~\ref{FIG:B = 0.8 at t' = 0.0236 ON A UNIFORM GRID},
\ref{FIG:B = 0.8 at t' = 0.0236 WITH AMR},
\ref{FIG:B = 0.9 at t' = 0.0236 ON A UNIFORM GRID},
and~\ref{FIG:B = 0.9 at t' = 0.0236 WITH AMR}, even though under magnification the fingers in the refined computations shown in
Figures~\ref{FIG:B = 0.8 at t' = 0.0236 WITH AMR}
and~\ref{FIG:B = 0.9 at t' = 0.0236 WITH AMR} do not yet appear fully resolved.
In other words, depending on the application, these computations may or may not be well resolved enough for the user to arrive at conclusions
\textit{appropriate for their application} concerning the flow at this point in time.

In conclusion, we emphasize that the required degree of resolution for a given computation will depend on the purpose of the computation and the user's need for fine detail as opposed to general qualitative information concerning the flow.

\section{Conclusions}
\label{Section:Conclusions}

We have implemented a Volume-of-Fluid (VOF) interface tracking method in the open source finite element code ASPECT, which is designed to model convection and other processes in the Earth's mantle.
Our VOF method works efficiently and effectively in ASPECT's parallel environment and with it's adaptive mesh refinement (AMR) algorithm.
We show that the VOF method reproduces linear interfaces in a constant flow to machine precision and is second-order accurate when we use it to compute a standard, smooth, interface tracking benchmark problem.
We also demonstrate that the method shows excellent agreement with two benchmark problems from the computational mantle convection literature.
In particular, in the second of these benchmarks we use AMR to allow us to compute at a much higher effective resolution at lower computational cost than would otherwise be possible. 


Finally, we use the new interface tracking methodology to study a problem involving thermochemical convection in density stratified flow. 
This model problem is relevant to the study of structures at the core mantle boundary known as  Large Low Shear Velocity Provinces (LLSVPs).
Recent studies utilizing seismic imaging have revealed large regions with anomalous seismic properties in the lower mantle.
There are two dome-like regions beneath Africa and the Pacific with low shear-wave velocities that extend some 1000 km above the core-mantle boundary and have horizontal dimensions of several thousand kilometers~\cite{SC-BR:2012,SWF-BR:2015}.
Most interpretations propose that the heterogeneities are compositional in nature, differing from the surrounding mantle, an interpretation that would be consistent with chemical geodynamic models.
Based on geological and geochemical studies it has been argued that LLSVPs have persisted for
billions of years~\cite{KB-BS-THT-MAS:2008}.

The model problem is designed to study the basic physics underlying the formation of thermal plumes that bring some of this material to the Earth's surface.
In our computations of we use AMR to obtain an effective grid resolution of $768 \times 256$ square cells overlaying the fluid interface on an underlying grid of $192 \times 64$ square cells.
This increase in resolution confirms that for a certain range of the nondimensional buoyancy  parameter $\mathrm{B}$ at Rayleigh number $\mathrm{Ra} = 10^5$ our computations of the interface have converged well enough to interpret with confidence the large scale dynamics of the two regions of differing densities.

In conclusion, the results of the work presented here demonstrate that our VOF interface tracking method should perform well on a number of problems of interest to the computational mantle convection community.


\section*{Acknowledgements}
This work was supported by the National Science Foundation's (NSF) SI2-SSE Program under 
Award number 1440811.
The development of ASPECT was supported by the Computational Infrastructure for
Geodynamics (CIG) under NSF Award numbers 0949446 and 1550901. 
The computations were made under the auspices of CIG on the U.C.~Davis Division of 
Mathematical and Physical Sciences distributed computing cluster Peloton.


\bibliographystyle{elsarticle/elsarticle-num}

\bibliography{COMPUTERS_AND_FLUIDS_PAPER}

\begin{thebibliography}{76}
\expandafter\ifx\csname natexlab\endcsname\relax\def\natexlab#1{#1}\fi
\expandafter\ifx\csname url\endcsname\relax
  \def\url#1{\texttt{#1}}\fi
\expandafter\ifx\csname urlprefix\endcsname\relax\def\urlprefix{URL }\fi

\bibitem[{Anbarlooei and Mazaheri(2011)}]{HRA-KM:2011}
Anbarlooei, H.~R., Mazaheri, K., 2011. {Moment of fluid} interface
  reconstruction method in axisymmetric coordinates. International Journal for
  Numerical Methods in Biomedical Engineering 27~(10), 1640--1651.

\bibitem[{Arndt et~al.(2017{\natexlab{a}})Arndt, Bangerth, Davydov, Heister,
  Heltai, Kronbichler, Maier, Pelteret, Turcksin, and
  Wells}]{DA-WB-DD-TH-LH-MK-MM-JPP-BT-DW:2017}
Arndt, D., Bangerth, W., Davydov, D., Heister, T., Heltai, L., Kronbichler, M.,
  Maier, M., Pelteret, J.-P., Turcksin, B., Wells, D., 2017{\natexlab{a}}. The
  \texttt{deal.II} library, version 8.5. Journal of Numerical Mathematics
  25~(3), 137--146.

\bibitem[{Arndt et~al.(2017{\natexlab{b}})Arndt, Bangerth, Davydov, Heister,
  Heltai, Kronbichler, Maier, Pelteret, Turcksin, and Wells}]{dealII85}
Arndt, D., Bangerth, W., Davydov, D., Heister, T., Heltai, L., Kronbichler, M.,
  Maier, M., Pelteret, J.-P., Turcksin, B., Wells, D., 2017{\natexlab{b}}. The
  \texttt{deal.II} library, version 8.5. Journal of Numerical Mathematics.

\bibitem[{Bangerth et~al.(2017{\natexlab{a}})Bangerth, Dannberg, Gassmoeller,
  Heister, et~al.}]{ASPECT-v1_5_0:2017}
Bangerth, W., Dannberg, J., Gassmoeller, R., Heister, T., et~al., March
  2017{\natexlab{a}}. {ASPECT}.

\bibitem[{Bangerth et~al.(2017{\natexlab{b}})Bangerth, Dannberg,
  Gassm{\"o}ller, Heister, et~al.}]{ASPECT-MANUAL:2017}
Bangerth, W., Dannberg, J., Gassm{\"o}ller, R., Heister, T., et~al.,
  2017{\natexlab{b}}. {ASPECT: A}dvanced Solver for Problems in Earth's
  ConvecTion, User Manual. CIG.
\newline\urlprefix\url{http://www.math.clemson.edu/~heister/manual.pdf}

\bibitem[{Bangerth et~al.(2007)Bangerth, Hartmann, and
  Kanschat}]{WB-RH-GK:2007}
Bangerth, W., Hartmann, R., Kanschat, G., 2007. {deal.II} -- a general purpose
  object oriented finite element library. ACM Trans. Math. Softw. 33~(4),
  24/1--24/27.

\bibitem[{Bars and Davaille(2004)}]{MLB-AD:2004}
Bars, M.~L., Davaille, A., 2004. Large interface deformation in two-layer
  thermal convection of miscible viscous fluids. J. Fluid Mech. 499, 75–110.

\bibitem[{Bars and Davaille(2005)}]{MLB-AD:2005}
Bars, M.~L., Davaille, A., 2005. Thermochemical convection in two superimposed
  miscible viscous fluids. In: Gutkowski, W., Kowalewski, T. (Eds.), {Mechanics
  of the 21st Century, Proceedings of the 21st International Congress of
  Theoretical and Applied Mechanics}. Springer Verlag, pp. FM7--12126.

\bibitem[{Burke et~al.(2008)Burke, Steinberger, Torsvik, and
  Smethurst}]{KB-BS-THT-MAS:2008}
Burke, K., Steinberger, B., Torsvik, T.~H., Smethurst, M.~A., 2008. Plume
  generation zones at the margins of large low shear velocity provinces on the
  core--mantle boundary. Earth and Planetary Science Letters 265~(1), 49--60.

\bibitem[{Burstedde et~al.(2011{\natexlab{a}})Burstedde, Wilcox, and
  Ghattas}]{p4est}
Burstedde, C., Wilcox, L.~C., Ghattas, O., 2011{\natexlab{a}}.
  {\texttt{p4est}}: Scalable algorithms for parallel adaptive mesh refinement
  on forests of octrees. SIAM Journal on Scientific Computing 33~(3),
  1103--1133.

\bibitem[{Burstedde et~al.(2011{\natexlab{b}})Burstedde, Wilcox, and
  Ghattas}]{CB-LCW-OG:2011}
Burstedde, C., Wilcox, L.~C., Ghattas, O., 2011{\natexlab{b}}.
  {\texttt{p4est}}: Scalable algorithms for parallel adaptive mesh refinement
  on forests of octrees. {SIAM Journal on Scientific Computing} 33~(3),
  1103--1133.

\bibitem[{Chandrasekhar(1961)}]{SC:1961}
Chandrasekhar, S., 1961. Hydrodynamic and Hydromagnetic Stability. Dover, New
  York.

\bibitem[{Chorin(1985)}]{AJC:1985}
Chorin, A.~J., 1985. Curvature and solidification. J. Comput. Phys. 57,
  472--490.

\bibitem[{Chorin and Marsden(1993)}]{AJC-JEM:1993}
Chorin, A.~J., Marsden, J.~E., 1993. A Mathematical Introduction to Fluid
  Mechanics, 4th Edition. No.~4 in Texts in Applied Mathematics.
  Springer-Verlag, New York, qA901.C53 1992.

\bibitem[{Colella(1990)}]{PC:1990}
Colella, P., 1990. Multidimensional upwind methods for hyperbolic conservation
  laws. J. Comput. Phys. 87, 171--200.

\bibitem[{Cottaar and Romanowicz(2012)}]{SC-BR:2012}
Cottaar, S., Romanowicz, B., 2012. An unusually large {ULVZ} at the base of the
  mantle near {H}awaii. Earth and Planetary Science Letters 355, 213--222.

\bibitem[{Dannberg et~al.(2017)Dannberg, Eilon, Faul, Gassmöller, Moulik, and
  Myhill}]{JD-ZE-UF-RG-PM-RM:2017}
Dannberg, J., Eilon, Z., Faul, U., Gassmöller, R., Moulik, P., Myhill, R.,
  2017. The importance of grain size to mantle dynamics and seismological
  observations. Geochemistry, Geophysics, Geosystems 18~(8), 3034--3061.

\bibitem[{Dannberg and Heister(2016)}]{JD-TH:2016}
Dannberg, J., Heister, T., 2016. Compressible magma/mantle dynamics: 3d,
  adaptive simulations in {ASPECT}. Geophysical Journal International 207~(3),
  1343--1366.

\bibitem[{Davaille(1999)}]{AD:1999}
Davaille, A., 1999. Two-layer thermal convection in miscible viscous fluids. J.
  Fluid Mech. 379, 223–253.

\bibitem[{Donea and Huerta(2005)}]{JD-AH:2005}
Donea, J., Huerta, A., 2005. Steady Transport Problems. John Wiley and Sons.

\bibitem[{French and Romanowicz(2015)}]{SWF-BR:2015}
French, S.~W., Romanowicz, B., 2015. Broad plumes rooted at the base of the
  {E}arth's mantle beneath major hotspots. Nature 525~(7567), 95--99.

\bibitem[{Gassmoeller(2016)}]{RG:2016}
Gassmoeller, R., 2016. Open source support for massively parallel, generic
  finite element methods. {Poster presented at the 2016 NSF SI2 PI Workshop}.
\newline\urlprefix\url{http://maxim.ucsd.edu/suave/index.html?file=si2n.cxml}

\bibitem[{Gerya and Yuen(2003)}]{TVG-DAY:2003}
Gerya, T.~V., Yuen, D.~A., 2003. Characteristics-based marker-in-cell method
  with conservative finite-differences schemes for modeling geological flows
  with strongly variable transport properties. Physics of the Earth and
  Planetary Interiors 140~(4), 293--318.

\bibitem[{Guermond et~al.(2011)Guermond, Pasquetti, and Popov}]{JLG-RP-BP:2011}
Guermond, J.-L., Pasquetti, R., Popov, B., 2011. Entropy viscosity method for
  nonlinear conservation laws. J. Comput. Phys. 230~(11), 4248 -- 4267,
  {Special issue High Order Methods for {CFD} Problems}.

\bibitem[{Hager and O{\textquoteright}Connell(1981)}]{BHH-RJO-1981}
Hager, B., O{\textquoteright}Connell, R., 1981. A simple global model of plate
  dynamics and mantle convection. Journal of Geophysical Research: Solid Earth
  86~(B6), 4843--4867.

\bibitem[{Hager and Clayton(1989)}]{BHH-RWC:1989}
Hager, B.~H., Clayton, R.~W., 1989. Constraints on the structure of mantle
  convection using seismic observations, flow models, and the geoid. In: Mantle
  Convection, Plate Tectonics and Global Dynamics. Gordon and Breach Science
  Publishers, pp. 657--763.

\bibitem[{He et~al.(2017)He, Puckett, and Billen}]{YH-EGP-MIB:2017}
He, Y., Puckett, E.~G., Billen, M.~I., 2017. A discontinuous {G}alerkin method
  with a bound preserving limiter for the advection of non-diffusive fields in
  solid {E}arth geodynamics. PEPI 263, 23--37.

\bibitem[{Heister et~al.(2017)Heister, Dannberg, Gassm{\"o}ller, and
  Bangerth}]{TH-JD-RG-WB:2017}
Heister, T., Dannberg, J., Gassm{\"o}ller, R., Bangerth, W., 2017. High
  accuracy mantle convection simulation through modern numerical methods. {II}:
  {R}ealistic models and problems. Geophysical Journal International 210~(2),
  833--851.

\bibitem[{Helmsen et~al.(1997)Helmsen, Colella, and Puckett}]{JJH-PC-EGP:1997}
Helmsen, J.~J., Colella, P., Puckett, E.~G., 1997. Non-convex profile evolution
  in two dimensions using volume of fluids. Technical Report LBNL-40693,
  Lawrence Berkeley National Laboratory.

\bibitem[{Helmsen et~al.(1996)Helmsen, Colella, Puckett, and
  Dorr}]{JJH-PC-EGP-MD:1996}
Helmsen, J.~J., Colella, P., Puckett, E.~G., Dorr, M., January 1996. Two new
  methods for simulating photolithography development in three dimensions. In:
  Proceeedings of the 10th SPIE Optical/Laser Microlithography Conference. Vol.
  2726. SPIE, San Jose, CA, pp. 253--261.

\bibitem[{Henderson et~al.(1991)Henderson, Colella, and
  Puckett}]{LFH-PC-EGP:1991}
Henderson, L.~F., Colella, P., Puckett, E.~G., March 1991. On the refraction of
  shock waves at a slowfast gas interface. J. Fluid Mech. 224, 1--27.

\bibitem[{Hill and Shashkov(2013)}]{RNH-MS:2013}
Hill, R.~N., Shashkov, M., 2013. The symmetric moment-of-fluid interface
  reconstruction algorithm. J. Comput. Phys. 249, 180 -- 184.

\bibitem[{Hirt and Nichols(1981)}]{CWH-BDN:1981}
Hirt, C.~W., Nichols, B.~D., 1981. {Volume of Fluid (VOF)} method for the
  dynamics of free boundaries. J. Comput. Phys. 39, 201--225.

\bibitem[{Huber and Helmig(1999)}]{RH-RH:1999}
Huber, R., Helmig, R., 1999. Multiphase flow in heterogeneous porous media: A
  classical finite element method versus an implicit pressure--explicit
  saturation-based mixed finite element--finite volume approach. International
  Journal for Numerical Methods in Fluids 29~(8), 899--920.

\bibitem[{Jemison et~al.(2015)Jemison, Sussman, and Shashkov}]{MJ-MMS-MS:2015}
Jemison, M., Sussman, M., Shashkov, M., 2015. Filament capturing with the
  multimaterial moment-of-fluid method. J. Comput. Phys. 285, 149--172.

\bibitem[{King et~al.(1990)King, Raefsky, and Hager}]{SDK-AR-BHH:1990}
King, S., Raefsky, A., Hager, B., 1990. {ConMan: V}ectorizing a finite element
  code for incompressible two-dimensional convection in the
  {Earth{\textquoteright}s} mantle. Physics of the Earth and Planetary
  Interiors 59~(3), 195--207.

\bibitem[{Kothe et~al.(1999)Kothe, Puckett, and Williams}]{DBK-EGP-MWW-3:1999}
Kothe, D.~B., Puckett, E.~G., Williams, M.~W., 1999. Robust finite volume
  modeling of 3-d free surface flows on unstructured meshes. In: Proceedings of
  the 14th AIAA Computational Fluid Dynamics Conference. American Institute of
  Aeronautics and Astronautics, Norfolk, VA, pp. 1--6.

\bibitem[{Kronbichler et~al.(2012)Kronbichler, Heister, and
  Bangerth}]{MK-TH-WB:2012}
Kronbichler, M., Heister, T., Bangerth, W., 2012. High accuracy mantle
  convection simulation through modern numerical methods. Geophysical Journal
  International 191~(1), 12--29.

\bibitem[{LeVeque(1996)}]{RLV:1996}
LeVeque, R.~J., April 1996. High-resolution conservative algorithms for
  advection in incompressible flow. SIAM J. Numer. Anal. 33~(2), 627--665.

\bibitem[{Manga(1996)}]{MM:1996}
Manga, M., 1996. Mixing of heterogeneities in the mantle: Effect of viscosity
  differences. Geophysical Research Letters 23~(4), 403--406.

\bibitem[{Manga and Stone(1994)}]{MM-HAS:1993}
Manga, M., Stone, H., 1994. Interactions between bubbles in magmas and lavas:
  effects of bubble deformation. Journal of Volcanology and Geothermal Research
  63~(3), 267--279.

\bibitem[{Manga et~al.(1993)Manga, Stone, and O'Connell}]{MM-HAS-RJO:1993}
Manga, M., Stone, H.~A., O'Connell, R.~J., 1993. The interaction of plume heads
  with compositional discontinuities in the earth's mantle. Journal of
  Geophysical Research: Solid Earth 98~(B11), 19979--19990.

\bibitem[{McNamara and Zhong(2004)}]{AKM-SZ:2004}
McNamara, A., Zhong, S., 2004. Thermochemical structures within a spherical
  mantle: Superplumes or piles?: Thermochemical structures. Journal of
  Geophysical Research: Solid Earth 109~(B7), n/a--n/a.

\bibitem[{Miller and Puckett(1994)}]{GHM-EGP:1994}
Miller, G.~H., Puckett, E.~G., February 1994. Edge effects in
  molybdenum-encapsulated molten silicate shock wave targets. J. Appl. Phys.
  75~(3), 1426--1434.

\bibitem[{Miller and Puckett(1996)}]{GHM-EGP:1996}
Miller, G.~H., Puckett, E.~G., August 1996. A high-order {Godunov} method for
  multiple condensed phases. J. Comput. Phys. 128~(1), 134--164.

\bibitem[{Moresi and Gurnis(1996)}]{LM-MG:1996}
Moresi, L., Gurnis, M., 1996. Constraints on the lateral strength of slabs from
  three-dimensional dynamic flow models. Earth and Planetary Science Letters
  138~(1), 15--28.

\bibitem[{Nichols et~al.(1980)Nichols, Hirt, and Hotchkiss}]{BDN-CWH-RSH:1980}
Nichols, B.~D., Hirt, C.~W., Hotchkiss, R.~S., August 1980. {SOLA-VOF:} a
  solution algorithm for transient fluid flow with multiple free boundaries.
  Technical Report LA-8355, Los Alamos National Laboratory.

\bibitem[{Noh and Woodward(1976)}]{WFN-PRW:1976}
Noh, W.~F., Woodward, P.~R., June 28--July 3 1976. {SLIC (Simple Line Interface
  Calculation)}. In: van~de Vooren, A.~I., Zandbergen, P.~J. (Eds.),
  Proceedings of the Fifth International Conference on Numerical Methods in
  Fluid Dynamics. Vol.~59 of Lecture Notes in Physics. Springer-Verlag, Twente
  University, Enschede, pp. 330--340.

\bibitem[{Pilliod(1992)}]{JEP:1992}
Pilliod, J.~E., September 1992. An analysis of piecewise linear interface
  reconstruction algorithms for volume-of-fluid methods. {MS} {T}hesis,
  {G}raduate {G}roup in {A}pplied {M}athematics, {U}niversity of {C}alifornia,
  {D}avis.

\bibitem[{Pilliod and Puckett(1998)}]{JEP-EGP:1998}
Pilliod, J.~E., Puckett, E.~G., 1998. An unsplit, second-order accurate
  {Godunov} method for tracking deflagrations and detonations. In: Houwing, A.
  F.~P., Paull, A., Boyce, R.~R., Danehy, P.~M., Hannemann, H., Kurtz, J.~J.,
  McIntyre, T.~J., McMahon, S.~J., Mee, D.~J., Sandeman, R.~J., Tanno, H.
  (Eds.), Proceedings of the 21st International Symposium on Shock Waves.
  Vol.~II. Panther Publishing, Fyshwick, Australia, pp. 1053--1058.

\bibitem[{Pilliod and Puckett(2004)}]{JEP-EGP:2004}
Pilliod, J.~E., Puckett, E.~G., September 2004. Second-order accurate
  volume-of-fluid algorithms for tracking material interfaces. J. Comput. Phys.
  199~(2), 465--502.

\bibitem[{Puckett(1991)}]{EGP:1991}
Puckett, E.~G., 1991. A volume-of-fluid interface tracking algorithm with
  applications to computing shock wave refraction. In: Proceedings of the
  Fourth International Symposium on Computational Fluid Dynamics. pp. 933--938.

\bibitem[{Puckett(2010{\natexlab{a}})}]{EGP:2010a}
Puckett, E.~G., February 2010{\natexlab{a}}. On the second-order accuracy of
  volume-of-fluid interface reconstruction algorithms: Convergence in the max
  norm. CAMCoS 5~(1), 99--148.

\bibitem[{Puckett(2010{\natexlab{b}})}]{EGP:2010b}
Puckett, E.~G., October 2010{\natexlab{b}}. A volume-of-fluid interface
  reconstruction algorithm that is second-order accurate in the max norm.
  CAMCoS 5~(2), 199--220.

\bibitem[{Puckett(2014)}]{EGP:2014}
Puckett, E.~G., January 2014. On the second-order accuracy of volume-of-fluid
  interface reconstruction algorithms {II}: An improved constraint on the cell
  size. CAMCoS 8~(1), 123--158.

\bibitem[{Puckett et~al.(1997)Puckett, Almgren, Bell, Marcus, and
  Rider}]{EGP-ASA-JBB-DLM-WJR:1997}
Puckett, E.~G., Almgren, A.~S., Bell, J.~B., Marcus, D.~L., Rider, W.~J.,
  January 1997. A high-order projection method for tracking fluid interfaces in
  variable density incompressible flows. J. Comput. Phys. 130~(2), 269--282.

\bibitem[{Puckett and Miller(1996)}]{EGP-GHM:1996}
Puckett, E.~G., Miller, G.~H., 1996. The numerical computation of jetting
  impacts. In: Sturtevant, B., Shepherd, J.~E., Hornung, H. (Eds.), Proceedings
  of the 20th International Symposium on Shock Waves. Vol.~II. World
  Scientific, New Jersey, pp. 1467--1472.

\bibitem[{Puckett et~al.(2017)Puckett, Turcotte, He, Lokavarapu, Robey, and
  Kellogg}]{EGP-DLT-YH-HL-JM-LHK:2017}
Puckett, E.~G., Turcotte, D.~L., He, Y., Lokavarapu, H., Robey, J.~M., Kellogg,
  L.~H., 2017. New numerical approaches for modeling thermochemical convection
  in a compositionally stratified fluid. Physics of the Earth and Planetary
  Interiors~(online November 11, 2017).

\bibitem[{Samuel and Evonuk(2010)}]{HS-ME:2010}
Samuel, H., Evonuk, M., 2010. Modeling advection in geophysical flows with
  particle level sets. Geochemistry, Geophysics, Geosystems 11~(8).

\bibitem[{Scardovelli and Zaleski(2000)}]{RS-SZ:2000}
Scardovelli, R., Zaleski, S., 2000. Analytical relations connecting linear
  interfaces and volume fractions in rectangular grids. Journal of
  Computational Physics 164~(1), 228 -- 237.

\bibitem[{Schubert et~al.(2001)Schubert, Turcotte, and Olson}]{GS-DLT-PO:2001}
Schubert, G., Turcotte, D.~L., Olson, P., 2001. Mantle convection in the
  {E}arth and planets. Cambridge University Press.

\bibitem[{Sethian(1999)}]{JAS:1999}
Sethian, J.~A., June 1999. Level Set Methods and Fast Marching Methods:
  Evolving Interfaces in Computational Geometry, Fluid Mechanics, Computer
  Vision, and Materials Sciences, 2nd Edition. Vol.~3 of Cambridge Monographs
  on Applied and Computational Mathematics. Cambridge University Press,
  Cambridge, U.K. ; New York.

\bibitem[{Sheldon et~al.(1959)Sheldon, Cardwell~Jr, et~al.}]{JWS-WTC:1959}
Sheldon, J., Cardwell~Jr, W., et~al., 1959. One-dimensional, incompressible,
  noncapillary, two-phase fluid flow in a porous medium. Petroleum
  Transactions, AIME 216, 290--296.

\bibitem[{Steinberger(2000)}]{BS:2000}
Steinberger, B., 2000. Plumes in a convecting mantle: Models and observations
  for individual hotspots. Journal of Geophysical Research: Solid Earth
  105~(B5), 11127--11152.

\bibitem[{Strang(1968)}]{WGS:1968}
Strang, W.~G., 1968. On the construction and comparison of difference schemes.
  SIAM J. Numer. Anal. 5~(3), 506--517.

\bibitem[{Strang(2016)}]{WGS:2016}
Strang, W.~G., 2016. Introduction to Linear Algebra, 5th Edition. Cambridge
  Wellsely Press.

\bibitem[{Sussman and Puckett(2000)}]{MMS-EGP:2000}
Sussman, M.~S., Puckett, E.~G., 2000. A coupled level set and volume of fluid
  method for computing {3D} and axisymmetric incompressible two-phase flows. J.
  Comput. Phys. 162, 301--337.

\bibitem[{Tan et~al.(2006)Tan, Choi, Thoutireddy, Gurnis, and
  Aivazis}]{ET-EC-PT-MG-MA:2006}
Tan, E., Choi, E., Thoutireddy, P., Gurnis, M., Aivazis, M., 2006.
  Geoframework: Coupling multiple models of mantle convection within a
  computational framework: Geoframework-mantle convection models. Geochemistry,
  Geophysics, Geosystems 7~(6), n/a--n/a.

\bibitem[{Torrey et~al.(1985)Torrey, Cloutman, Mjolsness, and
  Hirt}]{MDT-LDC-RCM-CWH:1985}
Torrey, M.~D., Cloutman, L.~D., Mjolsness, R.~C., Hirt, C.~W., December 1985.
  {NASA-VOF2D: A} computer program for incompressible flows with free surfaces.
  Technical Report LA-10612-MS, Los Alamos National Laboratory.

\bibitem[{Torrey et~al.(1987)Torrey, Mjolsness, and Stein}]{MDT-RCM-LRS:1987}
Torrey, M.~D., Mjolsness, R.~C., Stein, L.~R., July 1987. {NASA-VOF3D: A}
  three-dimensonal computer program for incompressible flows with free
  surfaces. Technical Report LA-11009-MS, Los Alamos National Laboratory.

\bibitem[{Tryggvason et~al.(2011)Tryggvason, Scardovelli, and
  Zaleski}]{GT-RS-SZ:2007}
Tryggvason, G., Scardovelli, R., Zaleski, S., 2011. Direct Numerical
  Simulations of Gas-Liquid Multiphase Flows, 1st Edition. Cambridge Monographs
  on Applied \& Computational Mathematics. Cambridge University Press,
  Cambridge.

\bibitem[{Turcotte and Schubert(2014)}]{DLT-GS:2014}
Turcotte, D.~L., Schubert, G., 2014. Geodynamics, 3rd Edition. Cambridge
  University Press.

\bibitem[{van Keken et~al.(1997)van Keken, King, Schmeling, Christensen,
  Neumeister, and Doin}]{PEVK-SDK-HS-URC-DN-MPD:1997}
van Keken, P.~E., King, S.~D., Schmeling, H., Christensen, U.~R., Neumeister,
  D., Doin, M.-P., 1997. A comparison of methods for the modeling of
  thermochemical convection. Journal of Geophysical Research: Solid Earth
  102~(B10), 22477--22495.

\bibitem[{Wanner and Hairer(1991)}]{GW-EH:1991}
Wanner, G., Hairer, E., 1991. Solving ordinary differential equations II.
  Vol.~14 of Springer Series in Computational Mathematics. Springer-Verlag
  Berlin Heidelberg.

\bibitem[{Zhong(2006)}]{SZ:2006}
Zhong, S., 2006. Constraints on thermochemical convection of the mantle from
  plume heat flux, plume excess temperature, and upper mantle temperature.
  Journal of Geophysical Research 111~(B4).

\bibitem[{Zhong et~al.(2000)Zhong, Zuber, Moresi, and
  Gurnis}]{SZ-MTZ-LM-MG:2000}
Zhong, S., Zuber, M., Moresi, L., Gurnis, M., 2000. Role of
  temperature-dependent viscosity and surface plates in spherical shell models
  of mantle convection. Journal of Geophysical Research: Solid Earth 105~(B5),
  11063--11082.

\end{thebibliography}

\end{document}